\documentclass[11pt]{amsart}

\usepackage{latexsym}
\usepackage{amssymb}
\usepackage{stmaryrd}

\newcommand{\excise}[1]{}

\numberwithin{section}{part}
\newtheorem{thm}{Theorem}[section]
\newtheorem{lemma}[thm]{Lemma}
\newtheorem{claim}[thm]{Claim}
\newtheorem{cor}[thm]{Corollary}
\newtheorem{prop}[thm]{Proposition}

\newtheorem{question}[thm]{Question}

\newtheorem{theorem}{Theorem}

\theoremstyle{definition}
\newtheorem{example}[thm]{Example}
\newtheorem{remark}[thm]{Remark}

\newtheorem{defn}[thm]{Definition}

\newenvironment{sketch}{\begin{trivlist}\item {\it
        Sketch of proof.\,}}{\mbox{}\hfill$\square$\end{trivlist}}

\newenvironment{proofof}[1]{\begin{trivlist}\item {\bf
        Proof of {#1}.\,}}{\mbox{}\hfill$\square$\end{trivlist}}

\newenvironment{rcgraph}{\begin{trivlist}\item\centering\footnotesize$}
                        {$\end{trivlist}}

\def\hln{\\[-.2ex]\hline}

\newenvironment{numbered}%
        {\begin{list}
                {\noindent\makebox[0mm][r]{\arabic{enumi}.}}
                {\leftmargin=5.5ex \usecounter{enumi}}
        }
        {\end{list}}

\newenvironment{romanlist}%
        {\begin{list}
                {\noindent\makebox[0mm][r]{(\roman{enumi})}}
                {\leftmargin=5.5ex \usecounter{enumi}}
        }
        {\end{list}}

\def\bem#1{\textbf{#1}}

\def\<{\langle}
\def\>{\rangle}
\def\0{{\mathbf 0}}
\def\1{{\mathbf 1}}
\def\AA{{\mathbb A}}

\def\CC{{\mathbb C}}
\def\cC{{\mathcal C}}

\def\EE{{\mathcal E}}
\def\FF{{\mathcal F}}
\def\GG{{\mathcal G}}

\def\JJ{{\mathcal J}}
\def\KK{{\mathcal K}}

\def\LL{{\mathcal L}}

\def\NN{{\mathbb N}}
\def\OO{{\mathcal O}}
\def\PP{{\mathbb P}}

\def\SS{{\mathfrak S}}

\def\WW{{\mathbf \Omega}}
\def\XX{{\mathcal X}}
\def\YY{{\mathcal Y}}
\def\ZZ{{\mathbb Z}}
\def\aa{{\mathbf a}}
\def\bb{{\mathbf b}}
\def\cc{{\mathbf c}}

\def\kk{\Bbbk}
\def\mm{{\mathfrak m}}
\def\pp{{\mathfrak p}}
\renewcommand\tt{{\mathbf t}}
\def\xx{{\mathbf x}}
\def\yy{{\mathbf y}}
\def\zz{{\mathbf z}}

\def\id{{\rm id}}

\def\rk{{\rm rk}}

\def\th{{\rm th}}

\def\del{{\rm del}}
\def\Hom{{\rm Hom}}

\def\std{{\rm std}}
\def\sym{{\rm Sym}}

\def\crit{{\rm crit}}

\def\link{{\rm link}}

\def\mult{{\rm mult}}
\def\prom{{\rm prom}}
\def\rank{{\rm rank}}
\def\SPEC{{\bf Spec}}
\def\spec{{\rm Spec}}
\def\supp{{\rm supp}}

\def\west{{\rm west}}

\def\cross{\textrm{`+'} }

\def\start{{\rm start}}

\def\facets{{\rm facets}}
\def\length{{\rm length}}
\def\weight{{\rm wt}}
\def\mitosis{{\rm mitosis\hspace{1pt}}}

\def\K{$K$}
\def\FL{{\mathcal F}{\ell}_n}
\def\FLN{{\mathcal F}{\ell}_N}
\def\IN{\mathsf{in}}
\def\mn{{M_n}}

\def\rc{\mathcal{RP}}
\def\dom{\backslash}

\def\fln{{B \backslash \gln\CC}}

\def\gln{{G_{\!}L_n}}
\def\glN{{G_{\!}L_N}}
\def\too{\longrightarrow}
\def\Znn{\ZZ^{\raisebox{-.1ex}%
	{$\scriptstyle n^{\raisebox{-.25ex}{$\scriptscriptstyle 2$}}$}}}
\def\Znp{\ZZ^{\raisebox{-0ex}%
	{$\scriptstyle d^{\raisebox{-.1ex}{$\scriptscriptstyle \prime$}}$}}}
\def\bgln{B\dom\hspace{-.2ex}\gln}
\def\bglN{B\dom\hspace{-.2ex}\glN}
\def\from{\leftarrow}
\def\hhom{{\mathcal{H}\hspace{-.2ex}\mathit{om}}}
\def\into{\hookrightarrow}
\def\spot{{\hbox{\raisebox{1.7pt}{\large\bf .}}\hspace{-.5pt}}}
\def\onto{\twoheadrightarrow}

\def\adots{{.\hspace{1pt}\raisebox{2pt}{.}\hspace{1pt}\raisebox{4pt}{.}}}
\def\minus{\smallsetminus}
\def\congto{\stackrel{\begin{array}{@{}c@{\;}}
                        \\[-4ex]\scriptstyle\approx\\[-1.6ex]
                      \end{array}}\to}
\def\implies{\Rightarrow}
\def\nothing{\varnothing}
\def\bigcupdot{\makebox[0pt][l]{$\hspace{1.05ex}\cdot$}\textstyle\bigcup}
\def\toplinedots{\hfill\raisebox{-5.1pt}[0pt][0pt]{\ $\ldots$\ }\hfill}
\renewcommand\iff{\Leftrightarrow}

\def\cj#1{\kk[\zz]/J_{#1}}
\def\ci#1{\kk[\zz]/I_{#1}}
\def\ol#1{{\overline {#1}}}
\def\ub#1{\hbox{\underbar{$#1$\hspace{-.25ex}}\hspace{.25ex}}}
\def\wt#1{{\widetilde{#1}}}

\def\dem#1{{\ol \partial_{#1}}}
\def\sub#1#2{_{#1 \times #2}}
\def\ddem#1#2{{\varepsilon_{#1}^{#2}}}

\newcommand{\aoverb}[2]{{\genfrac{}{}{0pt}{1}{#1}{#2}}}

\font\co=lcircle10
\def\petit#1{{\scriptstyle #1}}

\def\jr{\smash{\raise2pt\hbox{\co \rlap{\rlap{\char'005} \char'007}}
               \raise6pt\hbox{\rlap{\vrule height6.5pt}}
               \raise2pt\hbox{\rlap{\hskip4pt \vrule height0.4pt depth0pt
                width7.7pt}}}}
\def\je{\smash{\raise2pt\hbox{\co \rlap{\rlap{\char'005}
                \phantom{\char'007}}}\raise6pt\hbox{\rlap{\vrule height6pt}}}}
\def\+{\smash{\lower2pt\hbox{\rlap{\vrule height14pt}}
                \raise2pt\hbox{\rlap{\hskip-3pt \vrule height.4pt depth0pt
                width14.7pt}}}}
\def\perm#1#2{\hbox{\rlap{$\petit {#1}_{\scriptscriptstyle #2}$}}%
                \phantom{\petit 1}}
\def\phperm{\phantom{\perm w3}}
\def\mcc{\multicolumn{1}{@{}c@{}}}

\def\ld#1#2{\raisebox{#1 pt}%
	{${\scriptstyle\partial}_{\scriptscriptstyle#2}\!\!\!$}}
\def\rd#1#2{\raisebox{#1 pt}%
	{$\!\!\!{\scriptstyle\partial}_{\scriptscriptstyle#2}$}}
\def\textcross{\ \smash{\lower4pt\hbox{\rlap{\hskip4.15pt\vrule height14pt}}
                \raise2.8pt\hbox{\rlap{\hskip-3pt \vrule height.4pt depth0pt
                width14.7pt}}}\hskip12.7pt}
\def\textelbow{\ \hskip.1pt\smash{\raise2.8pt%
                \hbox{\co \hskip 4.15pt\rlap{\rlap{\char'005} \char'007}
                \lower6.8pt\rlap{\vrule height3.5pt}
                \raise3.6pt\rlap{\vrule height3.5pt}}
                \raise2.8pt\hbox{%
                  \rlap{\hskip-7.15pt \vrule height.4pt depth0pt width3.5pt}%
                  \rlap{\hskip4.05pt \vrule height.4pt depth0pt width3.5pt}}}
                \hskip8.7pt}

\begin{document}

\title{Gr\"obner geometry of Schubert polynomials}
\author{Allen Knutson}
\thanks{AK was partly supported by the Clay Mathematics Institute, Sloan
	Foundation, and NSF}
\address{Mathematics Department\\ UC Berkeley\\ Berkeley, California}
\email{allenk@math.berkeley.edu}
\author{Ezra Miller}
\thanks{EM was supported by the Sloan Foundation and NSF}
\address{Mathematical Sciences Research Institute\\Berkeley, CA}
\email{ezra@math.umn.edu}
\date{9 September 2003}

\begin{abstract}
\noindent
Given a permutation $w \in S_n$, we consider a determinantal ideal
$I_w$ whose generators are certain minors in the generic $n \times n$
matrix (filled with independent variables).  Using `multidegrees' as
simple algebraic substitutes for torus-equivariant cohomology classes
on vector spaces, our main theorems describe, for each ideal~$I_w$:
\begin{itemize}
\item
variously graded multidegrees and Hilbert series in terms of ordinary
and double Schubert and Grothendieck polynomials;

\item
a Gr\"obner basis consisting of minors in the generic $n \times n$
matrix;

\item
the Stanley--Reisner simplicial complex of the initial ideal in terms of
known combinatorial diagrams \cite{FKyangBax,BB} associated to
permutations in~$S_n$; and

\item
a procedure inductive on weak Bruhat order for listing the facets of
this complex.
\end{itemize}
We show that the initial ideal is Cohen--Macaulay, by identifying the
Stanley--Reisner complex as a special kind of ``subword complex
in~$S_n$'', which we define generally for arbitrary Coxeter groups, and
prove to be shellable by giving an explicit vertex decomposition.
We~also prove geometrically a general positivity statement for
multidegrees of subschemes.

Our main theorems provide a geometric explanation for the naturality of
Schubert polynomials and their associated combinatorics.  More
precisely, we \mbox{apply these theorems to}:
\begin{itemize}
\item
define a single geometric setting in which polynomial representatives
for Schubert classes in the integral cohomology ring of the flag
manifold are determined uniquely, and have positive coefficients for
geometric reasons;


\item
rederive from a topological perspective Fulton's Schubert polynomial
formula for universal cohomology classes of degeneracy loci of maps
between flagged vector bundles;

\item
supply new proofs that Schubert and Grothendieck polynomials represent
cohomology and $K$-theory classes on the flag manifold; and

\item
provide determinantal formulae for the multidegrees of ladder
determinantal rings.
\end{itemize}

The proof of the main theorems introduces the technique of ``Bruhat
induction'', consisting of a collection of geometric, algebraic, and
combinatorial tools, based on divided and isobaric divided differences,
that allow one to prove statements about determinantal ideals by
induction on weak Bruhat order.
\end{abstract}

\maketitle

{}

\setcounter{tocdepth}{1}
\tableofcontents

\part*{Introduction}


\noindent
The manifold~$\FL$ of complete flags (chains of vector subspaces) in the
vector space~$\CC^n$ over the complex numbers has historically been a
focal point for a number of distinct fields within mathematics.
By definition, $\FL$ is an object at the intersection of algebra and
geometry.  The fact that $\FL$ can be expressed as the quotient $\bgln$
of all invertible $n \times n$ matrices by its subgroup of lower
triangular matrices places it within the realm of Lie group theory, and
explains its appearance in representation theory.  In topology, flag
manifolds arise as fibers of certain bundles constructed universally
from complex vector bundles, and in that context the cohomology ring
$H^*(\FL) = H^*(\FL;\ZZ)$ with integer coefficients~$\ZZ$ plays an
important role.  Combinatorics, especially related to permutations of a
set of cardinality~$n$,
aids in understanding the topology of~$\FL$ in a geometric manner.

To be more precise, the cohomology ring $H^*(\FL)$ equals---in a
canonical way---the quotient of a polynomial ring $\ZZ[x_1,\ldots,x_n]$
modulo the ideal generated by all nonconstant homogeneous functions
invariant under permutation of the indices $1,\ldots,n$ \cite{Borel53}.
This quotient is a free abelian group of rank~$n!$ and has a basis given
by monomials
dividing $\prod_{i=1}^{n-1} x_i^{n-i}$.  This algebraic basis does not
reflect the geometry of flag manifolds as well as the basis of {\em
Schubert classes}, which are the cohomology classes of Schubert
varieties~$X_w$, indexed by permutations $w \in S_n$ \cite{Ehr}.  The
{\em Schubert variety}\/~$X_w$ consists of flags \mbox{$V_0 \subset
V_1 \subset \cdots \subset V_{n-1} \subset V_n$} whose intersections
$V_i \cap \CC^j$ have dimensions determined in a certain way by~$w$,
where $\CC^j$ is spanned by the first~$j$ basis vectors of~$\CC^n$.

A great deal of research has grown out of attempts to understand the
connection between the algebraic basis of monomials and the geometric
basis of Schubert classes $[X_w]$ in the cohomology ring~$H^*(\FL)$.
For this purpose, Lascoux and Sch\"utzenberger singled out Schubert
polynomials $\SS_w \in \ZZ[x_1,\ldots,x_n]$ as
representatives for Schubert classes \cite{LSpolySchub}, relying in
large part on earlier work of Demazure \cite{Dem} and
Bernstein--Gel$'$fand--Gel$'$fand \cite{BGG}.  Lascoux and
Sch\"utzenberger justified their choices with algebra and combinatorics,
whereas the earlier work had been in the context of geometry.  This
paper bridges the
algebra and combinatorics of Schubert polynomials on the one hand, with
geometry of Schubert varieties on the other.  In the process, it brings
a new perspective to problems in commutative algebra concerning ideals
generated by minors of generic matrices.

Combinatorialists have in fact recognized the intrinsic interest of
Schubert polynomials~$\SS_w$ for some time, and have therefore produced
a wealth of interpretations for their coefficients.  For example, see
\cite{BergeronCombConstSchub}, \cite[Appendix to Chapter IV, by N.\
Bergeron]{NoteSchubPoly}, \cite{BJS}, \cite{FKyangBax},
\cite{FSnilCoxeter}, \cite{Kohnert}, and \cite{Winkel}.  Geometers,
on the other hand, who take for granted Schubert {\em classes}\/~$[X_w]$
in cohomology of flag manifold~$\FL$, generally remain less convinced of
the naturality of Schubert {\em polynomials}, even though these
polynomials arise in certain universal geometric contexts
\cite{FulDegLoc}, and there are geometric proofs of positivity for their
\mbox{coefficients} \mbox{\cite{BerSot02, KoganThesis}}.

Our primary motivation for undertaking this project was to provide a
geometric context in which both (i)~polynomial representatives for
Schubert classes $[X_w]$ in the integral cohomology ring $H^*(\FL)$ are
uniquely singled out, with no choices other than a Borel subgroup of the
general linear group~$\gln\CC$; and (ii)~it is geometrically obvious
that these representatives have nonnegative coefficients.  That our
polynomials turn out to be the Schubert polynomials is a testament to
the naturality of Schubert polynomials; that our geometrically positive
formulae turn out to reproduce
known combinatorial structures is a testament to the naturality of the
combinatorics previously unconvincing to geometers.

The kernel of our idea was to translate ordinary cohomological
statements concerning Borel orbit closures on the flag manifold~$\FL$
into equivariant-cohomological statements concerning double Borel orbit
closures on the $n \times n$ matrices~$\mn$.  Briefly, the preimage
\mbox{$\tilde X_w \subseteq \gln$} of a Schubert variety $X_w \subseteq
\FL = \bgln$ is an orbit closure for the action of $B \times B^+$, where
$B$ and $B^+$ are the lower and upper triangular Borel subgroups
of~$\gln$ acting by multiplication on the left and right.  Letting $\ol
X_w \subseteq \mn$ be the closure of~$\tilde X_w$ and $T$ the torus
in~$B$, the $T$-equivariant cohomology class $[\ol X_w]_T \in H^*_T(\mn)
= \ZZ[x_1,\ldots,x_n]$ is our polynomial representative.
It has positive coefficients because there is a $T$-equivariant flat
(Gr\"obner) degeneration $\ol X_w \rightsquigarrow \LL_w$ to a union of
coordinate subspaces $L \subseteq \mn$.  Each subspace $L \subseteq
\LL_w$ has equivariant cohomology class $[L]_T \in H^*_T(\mn)$ that is a
monomial in $x_1,\ldots,x_n$, and the sum of these is $[\ol X_w]_T$.
Our obviously positive formula is thus simply
\begin{eqnarray} \label{pos}
  [\ol X_w]_T &=& [\LL_w]_T\:\ = \:\ \sum_{L \in \LL_w} [L]_T.
\end{eqnarray}

In fact, one need not actually {\em produce}\/ a degeneration of~$\ol
X_w$ to a union of coordinate subspaces: mere {\em existence}\/ of such
a degeneration is enough to conclude positivity of the cohomology
class~$[\ol X_w]_T$, although if the limit is nonreduced then subspaces
must be counted according to their (positive) multiplicities.  This
positivity holds quite generally for sheaves on vector spaces with torus
actions, because existence of degenerations is a standard consequence of
Gr\"obner basis theory.
That being said, in
our main results we identify a particularly natural
degeneration of the matrix Schubert variety~$\ol X_w$, with reduced and
Cohen--Macaulay limit~$\LL_w$, in which the subspaces have combinatorial
interpretations, and~\eqref{pos} coincides with the known combinatorial
formula \cite{BJS,FSnilCoxeter} for Schubert polynomials.

The above argument, as presented, requires equivariant cohomology
classes associated to closed subvarieties of noncompact spaces such
as~$\mn$, the subtleties of which might be considered unpalatable, and
certainly require characteristic zero.  Therefore we instead develop our
theory in the context of \bem{multidegrees},
which are algebraically defined substitutes.
In this setting, equivariant considerations for \bem{matrix Schubert
varieties} $\ol X_w \subseteq \mn$ guide our path directly toward
multigraded commutative algebra for the \bem{Schubert determinantal
ideals}~$I_w$ cutting out the varieties~$\ol X_w$.

\subsection*{Example}
Let $w=2143$ be the permutation in the symmetric group~$S_4$ sending $1
\mapsto 2$, $2 \mapsto 1$, $3 \mapsto 4$ and $4 \mapsto 3$.  The matrix
Schubert variety $\ol X_{2143}$ is the set of $4\times 4$ matrices $Z =
(z_{ij})$ whose upper-left entry is zero, and whose upper-left $3\times
3$ block has rank at most two.  The equations defining $\ol X_{2143}$
are the vanishing of the determinants
$$
\Big\< z_{11}, \quad
\left|\begin{array}{ccc}z_{11}&z_{12}&z_{13}\\
                        z_{21}&z_{22}&z_{23}\\
                        z_{31}&z_{32}&z_{33}\end{array}\right|
= -z_{13} z_{22} z_{31} + \ldots
\Big\>.
$$
When we Gr\"obner-degenerate the matrix Schubert variety to the scheme
defined by the initial ideal $\<z_{11}, -z_{13} z_{22} z_{31}\>$, we get
a union $\LL_{2143}$ of three coordinate subspaces
$$
  L_{11,13},\ L_{11,22}, \hbox{ and } L_{11,31}, \quad\hbox{with
  ideals}\quad \<z_{11},z_{13}\>,\ \<z_{11},z_{22}\>, \hbox{ and }
  \<z_{11},z_{31}\>.
$$
In the $\ZZ^n$-grading where $z_{ij}$ has weight $x_i$, the multidegree
of $L_{i_1j_1,i_2j_2}$ equals $x_{i_1}x_{i_2}$.  Our ``obviously
positive'' formula~\eqref{pos} for $\SS_{2143}(\xx)$ says that $[\ol
X_{2143}]_T = x_1^2 + x_1x_2 + x_1x_3$.

Pictorially, we represent the subspaces $L_{11,13}$, $L_{11,22}$, and
$L_{11,31}$ inside $\LL_{2143}$ as subsets
\begin{rcgraph}
\hbox{\normalsize{$\<z_{11},z_{13}\>\ =\ $}}\
\begin{array}{|c|c|c|c|}
\multicolumn{4}{c}{}
\\[-2.5ex]
  \hline\!+\!&     &\!+\!&     
\\\hline     &\ \, &     &\ \, 
\\\hline     &     &     &     
\\\hline     &     &     &     
\\\hline
\end{array}
\:
\hbox{\normalsize{$,\quad\<z_{11},z_{22}\>\ =\ $}}\
\begin{array}{|c|c|c|c|}
\multicolumn{4}{c}{}
\\[-2.5ex]
  \hline\!+\!&     &\ \, &\ \, 
\\\hline     &\!+\!&     &     
\\\hline     &     &     &     
\\\hline     &     &     &     
\\\hline
\end{array}
\:
\hbox{\normalsize{$,\quad\<z_{11},z_{31}\>\ =\ $}}\
\begin{array}{|c|c|c|c|}
\multicolumn{4}{c}{}
\\[-2.5ex]
  \hline\!+\!&\ \, &\ \, &\ \,
\\\hline     &     &     &     
\\\hline\!+\!&     &     &     
\\\hline     &     &     &     
\\\hline
\end{array}
\end{rcgraph}
of the $4 \times 4$ grid, or equivalently as ``pipe dreams'' with
crosses~$\textcross$ and ``elbow joints''~$\textelbow$ instead of boxes
with $+$ or nothing, respectively (imagine~$\textelbow$ filling the
lower right
corners):
$$
\qquad\qquad\
\begin{array}{ccccc}
       &\perm1{}&\perm2{}&\perm3{}&\perm4{}\\
\petit2&   \+   &   \jr  &   \+   &  \je   \\
\petit1&   \jr  &   \jr  &   \je  &\\
\petit4&   \jr  &   \je  &        &\\
\petit3&   \je  &        &        &\\
\end{array}
\qquad\qquad\hbox{}\qquad\
\begin{array}{ccccc}
       &\perm1{}&\perm2{}&\perm3{}&\perm4{}\\
\petit2&   \+   &   \jr  &   \jr  &  \je   \\
\petit1&   \jr  &   \+   &   \je  &\\
\petit4&   \jr  &   \je  &        &\\
\petit3&   \je  &        &        &\\
\end{array}
\qquad\qquad\hbox{}\qquad\
\begin{array}{ccccc}
       &\perm1{}&\perm2{}&\perm3{}&\perm4{}\\
\petit2&   \+   &   \jr  &   \jr  &  \je   \\
\petit1&   \jr  &   \jr  &   \je  &\\
\petit4&   \+   &   \je  &        &\\
\petit3&   \je  &        &        &\\
\end{array}
$$
These are the three ``reduced pipe dreams'', or ``planar histories'',
for
$w = 2143$
\cite{FKyangBax}, so we recover the combinatorial formula for
$\SS_w(\xx)$ from \cite{BJS,FSnilCoxeter}.

\medskip
Our main `Gr\"obner geometry' theorems describe, for every matrix
Schubert variety~$\ol X_w$:
\begin{itemize}
\item
its multidegree and Hilbert series, in terms of Schubert and
Grothendieck polynomials (Theorem~\ref{t:formulae});

\item
a Gr\"obner basis consisting of minors in its defining ideal~$I_w$
(Theorem~\ref{t:gb});

\item
the Stanley--Reisner complex $\LL_w$ of its
initial ideal~$J_w$, which we prove is Cohen--Macaulay, in terms of pipe
dreams and combinatorics
of~$S_n$ (Theorem~\ref{t:gb}); and

\item
an inductive irredundant algorithm (`mitosis') on weak Bruhat order for
listing the facets of~$\LL_w$ (Theorem~\ref{t:mitosis}).
\end{itemize}
Gr\"obner geometry of Schubert polynomials thereby provides a geometric
explanation for the naturality of Schubert polynomials and their
associated combinatorics.

The divided and isobaric divided differences used by Lascoux and
Sch\"utzenberger to define Schubert and Gro\-then\-dieck polynomials
inductively \cite{LSpolySchub,LSgrothVarDrap} were originally invented
by virtue of their geometric interpretation by Demazure \cite{Dem} and
Bernstein--Gel$'$fand--Gel$'$fand \cite{BGG}.  The heart of our proof of
the Gr\"obner geometry theorem for Schubert polynomials captures the
divided and isobaric divided differences in their algebraic and
combinatorial manifestations.  Both manifestations are positive: one in
terms of the generators of the initial ideal~$J_w$ and the monomials
outside~$J_w$, and the other in terms of certain combinatorial diagrams
(reduced pipe dreams) associated to permutations by
\mbox{Fomin--Kirillov}~\cite{FKyangBax}.
Taken together, the geometric, algebraic, and combinatorial
interpretations provide a powerful inductive method, which we call
\bem{Bruhat induction}, for working with determinantal ideals and
their initial ideals, as they relate to multigraded cohomological and
combinatorial invariants.  In particular, Bruhat induction applied to
the facets of $\LL_w$ proves a geometrically motivated substitute for
Kohnert's conjecture \cite{Kohnert}.

At present, ``almost all of the approaches one can choose for the
investigation of determinantal rings use standard bitableaux and the
straightening law'' \cite[p.~3]{BruConKRSandDet}, and are thus
intimately tied to the Robinson--Schensted--Knuth correspondence.
Although Bruhat induction as developed here may seem similar in spirit
to RSK, in that both allow one to work directly with vector space bases
in the quotient ring, Bruhat induction contrasts with methods based on
RSK in that
it compares standard monomials of {\em different ideals}\/ inductively
on weak Bruhat order, instead of comparing distinct bases associated to
the same ideal, as RSK does.  Consequently, Bruhat induction encompasses
a substantially larger class of determinantal ideals.

Bruhat induction, as well as the derivation of the main theorems
concerning Gr\"obner geometry of Schubert polynomials from it, relies on
two general results concerning
\begin{itemize}
\item
positivity of multidegrees---that is, positivity of torus-equivariant
cohomology classes represented by subschemes or coherent sheaves on
vector spaces (Theorem~\ref{t:pos}); and

\item
shellability of certain simplicial complexes that reflect the nature of
reduced subwords of words in Coxeter generators for Coxeter groups
(Theorem~\ref{t:subword}).
\end{itemize}
The latter of these allows us to approach the combinatorics of Schubert
and Grothendieck polynomials from a new perspective, namely that of
simplicial topology.  More precisely, our proof of shellability for the
initial complex $\LL_w$ draws on previously unknown combinatorial
topological aspects of reduced expressions in symmetric groups, and more
generally in arbitrary Coxeter groups.  We touch relatively briefly on
this aspect of the story here, only proving what is essential for the
general picture in the present context, and refer the reader to
\cite{KMsubword} for a complete treatment, including applications to
Grothendieck polynomials.

\subsection*{Organization}
Our main results, Theorems~\ref{t:formulae}, \ref{t:gb},
\ref{t:mitosis}, \ref{t:pos}, and~\ref{t:subword}, appear in
Sections~\ref{sec:matSchub}, \ref{sec:grobGeom}, \ref{sec:alg},
\ref{sec:pos}, and~\ref{sec:subword}, respectively.  The sections in
Part~\ref{part:intro} are almost entirely expository in nature, and
serve not merely to define all objects appearing in the central
theorems, but also to provide independent motivation and examples for
the theories they describe.  For each of
Theorems~\ref{t:formulae},~\ref{t:gb}, \ref{t:mitosis},
and~\ref{t:subword}, we develop just enough prerequisites before it to
give a complete statement, while for Theorem~\ref{t:pos} we first
provide a crucial characterization of multidegrees, in
Theorem~\ref{t:multidegs}.

Readers seeing this paper for the first time should note that
Theorems~\ref{t:formulae}, \ref{t:gb}, and~\ref{t:pos} are core results,
not to be overlooked on a first pass through.  Theorems~\ref{t:mitosis}
and~\ref{t:subword} are less essential to understanding the main point
as outlined in the Introduction, but still fundamental for the
combinatorics of Schubert polynomials as derived from geometry via
Bruhat induction (which is used to prove Theorems~\ref{t:formulae}
and~\ref{t:gb}), and for substantiating the naturality of the
degeneration in Theorem~\ref{t:gb}.

The paper is structured logically as follows.
There are no proofs in Sections~\ref{sec:schub}--\ref{sec:alg}
except for a few easy lemmas that serve the exposition.  The complete
proof of Theorems~\ref{t:formulae},~\ref{t:gb}, and~\ref{t:mitosis} must
wait until the last section of Part~\ref{part:gb}
(Section~\ref{sec:proof}), because these results rely on Bruhat
induction.  Section~\ref{sec:proof} indicates which parts of the
theorems from Part~\ref{part:intro} imply the others, while gathering
the
results from Part~\ref{part:gb} to prove those required parts.  In
contrast, the proofs of Theorems~\ref{t:pos} and~\ref{t:subword} in
Sections~\ref{sec:pos} and~\ref{sec:subword} are completely
self-contained, relying on nothing other than definitions.  Results of
Part~\ref{part:intro} are used freely in Part~\ref{part:app} for
applications to consequences not found or only briefly mentioned in
Part~\ref{part:intro}.  The development of Bruhat induction in
Part~\ref{part:gb}
depends only on Section~\ref{sec:pos} and defintions from
Part~\ref{part:intro}.

In terms of content,
Sections~\ref{sec:schub}, \ref{sec:mult}, and~\ref{sec:pipe}, as well as
the first half of Section~\ref{sec:matSchub}, review known definitions,
while the other sections in Part~\ref{part:intro} introduce topics
appearing here for the first time.  In more detail,
Section~\ref{sec:schub} recalls the Schubert and Grothendieck
polynomials of Lascoux and Sch\"utzenberger via divided differences and
their isobaric relatives.  Then Section~\ref{sec:mult}
reviews \K-polynomials and multidegrees, which are rephrased versions of
the equivariant multiplicities in \cite{joseph,rossmann}.
We start Section~\ref{sec:matSchub} by introducing matrix Schubert
varieties and Schubert determinantal ideals, which are due (in different
language) to Fulton \cite{FulDegLoc}.  This discussion culminates in the
statement of Theorem~\ref{t:formulae}, giving the multidegrees and
\K-polynomials of matrix Schubert varieties.

We continue in Section~\ref{sec:pipe} with some combinatorial diagrams
that we call `reduced pipe dreams', associated to permutations.  These
were invented by Fomin and Kirillov and studied by Bergeron and Billey,
who called them `rc-graphs'.  Section~\ref{sec:grobGeom} begins with the
definition of `antidiagonal' squarefree monomial ideals, and proceeds to
state Theorem~\ref{t:gb}, which describes Gr\"obner bases and
initial ideals for matrix Schubert varieties in terms of reduced pipe
dreams.  Section~\ref{sec:alg} defines our combinatorial `mitosis' rule
for manipulating subsets of the $n \times n$ grid, and describes in
Theorem~\ref{t:mitosis} \mbox{how mitosis generates all reduced pipe
dreams}.

Section~\ref{sec:pos} works with multidegrees in the general context of
a positive multigrading, proving the characterization
Theorem~\ref{t:multidegs} and then its consquence, the Positivity
Theorem~\ref{t:pos}.  Also in a general setting---that of arbitrary
Coxeter groups---we define `subword complexes' in
Section~\ref{sec:subword}, and prove their vertex-decomposability in
Theorem~\ref{t:subword}.

Our most important application, in Section~\ref{sec:bjs}, consists of
the geometrically positive formulae for Schubert polynomials that
motivated this paper.  Other applications include connections with
Fulton's theory of degeneracy loci in Section~\ref{sec:loci}, relations
between our multidegrees and \K-polynomials on $n \times n$ matrices
with classical cohomological theories on the flag manifold in
Section~\ref{sec:flag}, and comparisons in Section~\ref{sec:ladder} with
the commutative algebra literature on determinantal~ideals.

Part~\ref{part:gb} demonstrates how the method of Bruhat induction
works geometrically, algebraically, and combinatorially to provide full
proofs of Theorems~\ref{t:formulae},~\ref{t:gb}, and~\ref{t:mitosis}.
We postpone the detailed overview of Part~\ref{part:gb} until
Section~\ref{sec:gb}, although we mention here that the geometric
Section~\ref{sec:multimat} has a rather different flavor than
Sections~\ref{sec:mutation}--\ref{sec:rp}, which deal mostly with the
combinatorial commutative algebra spawned by divided differences, and
Section~\ref{sec:proof}, which collects Part~\ref{part:gb} into a
coherent whole in order to prove Theorems~\ref{t:formulae},~\ref{t:gb},
and~\ref{t:mitosis}.  Generally speaking, the material in
Part~\ref{part:gb} is more technical than earlier parts.

We have tried to make the material here as accessible as possible to
combinatorialists, geometers, and commutative algebraists alike.  In
particular, except for applications in Part~\ref{part:app}, we have
assumed no specific knowledge of the algebra, geometry, or combinatorics
of flag manifolds, Schubert varieties, Schubert polynomials,
Grothendieck polynomials, or determinantal ideals.  Many of our examples
interpret the same underlying data in varying contexts, to highlight and
contrast common themes.  In particular this is true of
Examples~\ref{ex:intro}, \ref{ex:pipe}, \ref{ex:rp}, \ref{ex:gb},
\ref{ex:mitosis}, \ref{ex:t:mitosis},
\ref{ex:mutate}, \ref{ex:prom}, \ref{ex:ddem}, \ref{ex:ii},
\ref{ex:revert}, \ref{ex:offspring}, \ref{ex:mitosis'},
and~\ref{ex:remove}.

\subsection*{Conventions}
Throughout this paper, $\kk$ is an arbitary field.  In particular, we
impose no restrictions on its characteristic.  Furthermore, although
some geometric statements or arguments may seem to require that $\kk$ be
algebraically closed, this hypothesis could be dispensed with formally
by resorting to sufficiently abstruse language.

We consciously chose our notational conventions (with considerable
effort) to mesh with those of \cite{FulDegLoc}, \cite{LSpolySchub},
\cite{FKgrothYangBax}, \cite{HerzogTrung}, and \cite{BB} concerning
permutations ($w^T$ versus~$w$), the indexing on (matrix) Schubert
varieties and polynomials (open orbit corresponds to identity
permutation and smallest orbit corresponds to long word), the placement
of one-sided ladders (in the northwest corner as opposed to the
southwest), and reduced pipe dreams.  These conventions dictated our
seemingly idiosyncratic choices of Borel subgroups as well as the
identification $\FL \cong \bgln$ as the set of right cosets, and
resulted in our use of row vectors in $\kk^n$ instead of the usual
column vectors.  That there even existed consistent conventions came as
a relieving surprise.

\subsection*{Acknowledgements}
The authors are grateful to Bernd Sturmfels, who took part in the
genesis of this project, and to Misha Kogan, as well as to Sara Billey,
Francesco Brenti, Anders Buch, Christian Krattenthaler, Cristian Lenart,
Vic Reiner, Rich\'ard Rim\'anyi, Anne Schilling, Frank Sottile, and
Richard Stanley for inspiring conversations and references.  Nantel
Bergeron kindly provided \LaTeX\ macros for drawing pipe dreams.

\flushbottom

\part{The Gr\"obner geometry theorems}

\label{part:intro}

\section{Schubert and Grothendieck polynomials}\label{sec:schub}

We write all permutations in one-line (not cycle) notation, where
\mbox{$w = w_1 \ldots w_n$} sends $i \mapsto w_i$.  Set $w_0 = n \ldots
3 2 1$ equal to the \bem{long permutation} reversing the order of
$1,\ldots,n$.

\begin{defn} \label{defn:schub}
Let $R$ be a commutative ring, and $\xx = x_1,\ldots,x_n$ independent
variables.  The \bem{$i^\th$ divided difference operator} $\partial_i$
takes each polynomial $f \in R[\xx]$ to
\begin{eqnarray*}
  \partial_i f(x_1, x_2, \ldots) &=& \frac{f(x_1, x_2, \ldots, ) - f(x_1,
        \ldots, x_{i-1}, x_{i+1}, x_i, x_{i+2}, \ldots)}{x_i - x_{i+1}}.
\end{eqnarray*}
The \bem{Schubert polynomial} for~$w \in S_n$ is defined by the recursion
\begin{eqnarray*}
  \SS_{ws_i}(\xx) &=& \partial_i \SS_w(\xx)
\end{eqnarray*}
whenever $\length(ws_i) < \length(w)$, and the initial condition
$\SS_{w_0}(\xx) = \prod_{i=1}^n x_i^{n-i} \in \ZZ[\xx]$.  The
\bem{double Schubert polynomials} $\SS_w(\xx,\yy)$ are defined by the
same recursion, but starting from $\SS_{w_0}(\xx,\yy) = \prod_{i+j \leq
n} (x_i-y_j) \in \ZZ[\yy][\xx]$.
\end{defn}

In the definition of $\SS_w(\xx,\yy)$, the operator $\partial_i$ is to
act only on the $\xx$~variables and not on the $\yy$~variables.  Checking
monomial by monomial
verifies that $x_i-x_{i+1}$ divides the numerator of $\partial_i(f)$,
so $\partial_i(f)$ is again a polynomial, homogeneous of degree $d-1$
if~$f$ is homogeneous of degree~$d$.

\begin{example} \label{ex:schub}
Here are all of the Schubert polynomials for permutations in $S_3$, along
with the rules for applying divided differences.
$$
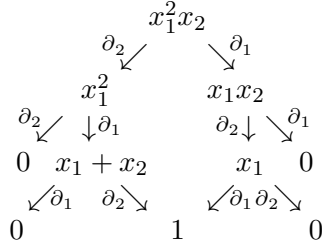
\begin{array}{@{}r@{\hspace{-1ex}}r@{}r@{}c@{}l@{}l@{\hspace{-1ex}}l@{}}
       &        &              &  x_1^2x_2   &              &        &
\\     &        & \ld52\swarrow&             &\searrow\rd51 &        &
\\     &        &x_1^2\ \ \ \  &             & \mcc{x_1x_2} &        &
\\&\ld52\swarrow\!\!&\mcc{\downarrow\;\rd21}&&\mcc{\ld22\;\downarrow}&\searrow\rd51
\\  0  &        &      x_1+x_2 &             &\ \ \ x_1     &        & 0
\\&\swarrow\rd01\!\!\!\!& \ld02\searrow&     &\swarrow\rd01 &\!\!\ld02\searrow&
\\0\: &        &              &      1      &              &        &\  0
\end{array}
$$
\end{example}

The recursion for both single and double Schubert polynomials can be
summarized as
\begin{eqnarray*}
  \SS_w &=& \partial_{i_k} \cdots \partial_{i_1} \SS_{w_0},
\end{eqnarray*}
where $w_0w = s_{i_1} \cdots s_{i_k}$ and\/ $\length(w_0w)=k$.  The
condition $\length(w_0w)=k$ means by definition that $k$ is minimal, so
$w_0w = s_{i_1} \cdots s_{i_k}$ is a \bem{reduced expression} for $w_0
w$.  It is not immediately obvious from Definition~\ref{defn:schub} that
$\SS_w$ is well-defined, but it follows from the fact that divided
differences satisfy the \bem{Coxeter relations}, $\partial_i
\partial_{i+1} \partial_i = \partial_{i+1} \partial_i \partial_{i+1}$
and $\partial_i \partial_{i'} = \partial_{i'} \partial_i$ when $|i-i'|
\geq 2$.

Divided differences arose geometrically in work of Demazure \cite{Dem}
and Bernstein--Gel$'$fand--Gel$'$fand \cite{BGG}, where they reflected a
`Bott--Samelson crank': form a $\PP^1$ bundle over a Schubert variety
and smear it out onto the flag manifold $\FL$ to get a Schubert variety
of dimension $1$~greater.  In their setting, the variables~$\xx$
represented Chern classes of standard line bundles $L_1,\ldots,L_n$
on~$\FL$, where the fiber of~$L_i$ over a flag \mbox{$F_0 \subset \cdots
\subset F_n$} is the dual vector space $(F_i/F_{i-1})^*$.  The divided
differences acted on the cohomology ring $H^*(\FL)$, which is the
quotient of~$\ZZ[\xx]$ modulo the ideal generated by symmetric functions
with no constant term~\cite{Borel53}.  The insight of Lascoux and
Sch\"utzenberger in~\cite{LSpolySchub} was to impose a stability
condition on the collection of polynomials~$\SS_w$ that defines them
uniquely among representatives for the cohomology classes of Schubert
varieties.  More precisely, although Definition~\ref{defn:schub} says
that $w$ lies in~$S_n$, the number~$n$ in fact plays no role: if $w_N
\in S_N$ for $n \geq N$ agrees with $w$ on $1,\ldots,n$ and fixes
$n+1,\ldots,N$, then $\SS_{w_N}(x_1,\ldots,x_N) =
\SS_w(x_1,\ldots,x_n)$.

The `double' versions represent Schubert classes in equivariant
cohomology for the Borel group action on~$\FL$.  As the ordinary
Schubert polynomials are much more common in the literature than double
Schubert polynomials, we have phrased many of our coming results both in
terms of Schubert polynomials as well as double Schubert polynomials.
This choice has the advantage of demonstrating how the notation
simplifies in the single case.

Schubert polynomials have their analogues in $K$-theory of~$\FL$, where
the recurrence uses a ``homogenized'' operator (sometimes called an
\bem{isobaric divided difference operator}):

\begin{defn} \label{defn:dem}
Let $R$ be a
commutative ring.  The \bem{$i^\th$ Demazure operator} $\dem i: R[[\xx]]
\to R[[\xx]]$ sends a power series $f(\xx)$ to
\begin{eqnarray*}
  \frac{x_{i+1}f(x_1, \ldots, x_n) - x_if(x_1, \ldots, x_{i-1}, x_{i+1},
  x_i, x_{i+2}, \ldots, x_n)}{x_{i+1} - x_i} &=& -\partial_i(x_{i+1}f).
\end{eqnarray*}
The \bem{Grothendieck polynomial} $\GG_w(\xx)$ is obtained recursively
from the ``top'' Grothendieck polynomial $\GG_{w_0}(\xx) :=
\prod_{i=1}^n (1-x_i)^{n-i}$ via the recurrence
\begin{eqnarray*}
  \GG_{ws_i}(\xx) &=& \dem i \GG_w(\xx)
\end{eqnarray*}
whenever $\length(ws_i) < \length(w)$.  The \bem{double Grothendieck
polynomials} are defined by the same recurrence, but start from
$\GG_{w_0}(\xx,\yy) := \prod_{i+j \leq n} (1-x_i y_j^{-1})$.
\end{defn}

As with divided differences, one can check directly that Demazure
operators $\dem i$ take power series to power series, and satisfy the
Coxeter relations.  Lascoux and Sch\"utzenberger \cite{LSgrothVarDrap}
showed that Grothendieck polynomials enjoy the same stability property
as do Schubert polynomials; we shall rederive this fact directly from
Theorem~\ref{t:formulae} in Section~\ref{sec:flag}
(Lemma~\ref{lemma:stable}), where we also
construct the bridge from Gr\"obner geometry of Schubert and
Grothendieck polynomials to classical geometry on flag manifolds.

Schubert polynomials represent data that are leading terms for the richer
structure encoded by Grothendieck polynomials.

\begin{lemma} \label{lemma:schubert}
The Schubert polynomial $\SS_w(\xx)$ is the sum of all lowest-degree
terms in $\GG_w(\1 - \xx)$, where $(\1-\xx) = (1-x_1,\ldots,1-x_n)$.
Similarly, the double Schubert polynomial $\SS_w(\xx,\yy)$ is the sum of
all lowest-degree terms in $\GG_w(\1-\xx,\1-\yy)$.
\end{lemma}
\begin{proof}
Assuming $f(\1-\xx)$ is homogeneous, plugging $\1-\xx$ for $\xx$ into the
first displayed equation in Definition~\ref{defn:dem} and taking the
lowest degree terms yields $\partial_i f(\1-\xx)$.  Since $\SS_{w_0}$ is
homogeneous, the result follows by induction on $\length(w_0w)$.
\end{proof}

Although the Demazure operators are usually applied only to polynomials
in~$\xx$, it will be crucial in our applications to use them on power
series in~$\xx$.
We shall also use the fact that, since the standard denominator $f(\xx) =
\prod_{i=1}^n(1-x_i)^n$ for $\ZZ^n$-graded Hilbert series over $\kk[\zz]$
is symmetric in $x_1,\ldots,x_n$, applying $\dem i$ to a Hilbert series
$g/\hspace{-.3ex}f$ simplifies: $\dem i(g/\hspace{-.3ex}f) = (\dem i
g)/\hspace{-.3ex}f$.  This can easily be checked directly.  The same
comment applies when $f(\xx) = \prod_{i,j=1}^n(1-x_i/y_j)$ is the
standard denominator for $\ZZ^{2n}$-graded Hilbert series.

\section{Multidegrees and \K-polynomials}\label{sec:mult}

Our first main theorem concerns cohomological and \K-theoretic
invariants of matrix Schubert varieties, which are given by
multidegrees and \K-polynomials, respectively.  We
work with these here in the setting of a polynomial ring~$\kk[\zz]$ in
$m$ variables $\zz = z_1,\ldots,z_m$, with a grading by $\ZZ^d$ in which
each variable $z_i$ has \bem{exponential weight} $\weight(z_i) =
\tt^{\aa_i}$ for some vector $\aa_i = (a_{i1},\ldots,a_{id}) \in \ZZ^d$,
where $\tt = t_1,\ldots,t_d$.  We~call $\aa_i$ the \bem{ordinary weight}
of~$z_i$, and sometimes write $\aa_i = \deg(z_i) = a_{i1} t_1 + \cdots +
a_{id}t_d$.  It can be useful to think of this as the logarithm of the
Laurent monomial~$\tt^{\aa_i}$.

\begin{example} \label{ex:exp}
Our primary concern is the case $\zz = (z_{ij})_{i,j=1}^n$ with various
gradings, in which the different kinds of weights are:
$$
\begin{array}{r|c@{\quad}c@{\quad}c@{\quad}c}
                      \hbox{grading} &\ZZ&\ZZ^n&\ZZ^{2n} &  \Znn
\\[0ex]\hline
\hbox{exponential weight of } z_{ij} & t & x_i & x_i/y_j &\!\!z_{ij}
\\[0ex]\hline
\hbox{ordinary weight of }    z_{ij} & t & x_i & x_i-y_j &\!\!z_{ij}
\end{array}
$$
The exponential weights are Laurent monomials that we treat as elements
in the group rings $\ZZ[t^{\pm 1}]$, $\ZZ[\xx^{\pm 1}]$, $\ZZ[\xx^{\pm
1},\yy^{\pm 1}]$, $\ZZ[\zz^{\pm 1}]$ of the grading groups.  The
ordinary weights are linear forms that we treat as elements in the
integral symmetric algebras $\ZZ[t] = \sym^\spot_\ZZ(\ZZ)$, $\ZZ[\xx] =
\sym^\spot_\ZZ(\ZZ^n)$, $\ZZ[\xx,\yy] = \sym^\spot_\ZZ(\ZZ^{2n})$,
$\ZZ[\zz] = \sym^\spot_\ZZ(\Znn)$ of the grading groups.%
%
\end{example}

Every finitely generated $\ZZ^d$-graded module $\Gamma = \bigoplus_{\aa
\in \ZZ^d} \Gamma_\aa$ over $\kk[\zz]$ has a free resolution
$$
  \EE_\spot: 0 \from \EE_0 \from \EE_1 \from \cdots \from \EE_m \from 0,
  \qquad\hbox{where}\qquad
  \EE_i = \bigoplus_{j=1}^{\beta_i} \kk[\zz](-\bb_{ij})
$$
is graded, with the $j^\th$ summand of $\EE_i$ generated in
$\ZZ^d$-graded degree $\bb_{ij}$.

\begin{defn} \label{d:Kpoly}
The \bem{\K-polynomial} of $\Gamma$ is $\KK(\Gamma;\tt) = \sum_i (-1)^i
\sum_j \tt^{\bb_{ij}}$.
\end{defn}

Geometrically, the $K$-polynomial of $\Gamma$ represents the class of the
sheaf $\tilde \Gamma$ on $\kk^m$ in equivariant $K$-theory for the action
of the $d$-torus whose weight lattice is $\ZZ^d$.
Algebraically, when the $\ZZ^d$-grading is \bem{positive}, meaning that
the ordinary weights $\aa_1,\ldots,\aa_d$ lie in a single open
half-space in $\ZZ^d$, the vector space dimensions
$\dim_\kk(\Gamma_\aa)$ are finite for all $\aa \in \ZZ^d$, and the
$K$-polynomial of~$\Gamma$ is the numerator of its $\ZZ^d$-graded
Hilbert series $H(\Gamma;\tt)$:
\begin{eqnarray*}
  H(\Gamma;\tt)\:\ :=\:\ \sum_{\aa \in \ZZ^d}
  \dim_\kk(\Gamma_\aa)\cdot\tt^\aa &=&
  \frac{\KK(\Gamma;\tt)}{\prod_{i=1}^m (1-\weight(z_i))}.
\end{eqnarray*}
We shall only have a need to consider positive multigradings in this
paper.

Given any Laurent monomial $\tt^\aa = t_1^{a_1} \cdots t_d^{a_d}$, the
rational function $\prod_{j=1}^d(1-t_j)^{a_j}$ can be expanded as a
well-defined (that is, convergent in the $\tt$-adic topology) formal
power series \mbox{$\prod_{j=1}^d(1 - a_jx_j + \cdots)$} in~$\tt$.
Doing the same for each monomial in an arbitrary Laurent polynomial
$\KK(\tt)$ results in a power series denoted by $\KK(\1-\tt)$.

\begin{defn} \label{defn:multideg}
The \bem{multidegree} of a $\ZZ^d$-graded $\kk[\zz]$-module $\Gamma$ is
the sum $\cC(\Gamma;\tt)$ of the lowest degree terms in
$\KK(\Gamma;\1-\tt)$.  If $\Gamma = \kk[\zz]/I$ is the coordinate ring
of a subscheme $X \subseteq \kk^m$, then we may also write $[X]_{\ZZ^d}$
or $\cC(X;\tt)$ to mean $\cC(\Gamma;\tt)$.
\end{defn}

Geometrically, multidegrees are just an algebraic reformulation of
torus-equivariant cohomology of affine space, or equivalently the
equivariant Chow ring \cite{TotChowRing,EGequivInt}.
Multidegrees originated in \cite{joseph}, and are called {\em
equivariant multiplicities}\/ in \cite{rossmann}.

\begin{example} \label{ex:subspace}
Let $n = 2$ in Example~\ref{ex:exp}, and set $\Gamma = \kk\Big[
\raisebox{.5ex}{\footnotesize$\begin{array}{@{}c@{\;}c@{}}z_{11}&z_{12}\\
z_{21}&z_{22}\end{array}$}\Big]/\<z_{11},z_{22}\>$.  Then
$$
  \KK(\Gamma;\zz) = {(1-z_{11})(1-z_{22})}
  \quad\hbox{and}\quad
  \KK(\Gamma;\xx,\yy) = {(1-x_1/y_1)(1-x_2/y_2)}
$$
because of the Koszul resolution.  Thus $\KK(\Gamma;\1-\zz) =
z_{11}z_{22} = \cC(\Gamma;\zz)$, and
\begin{eqnarray*}
  \KK(\Gamma;\1-\xx,\1-\yy) &=&
  (x_1 - y_1 + x_1y_1 - y_1^2 + \cdots)
  (x_2 - y_2 + x_2y_2 - y_2^2 + \cdots),
\end{eqnarray*}
whose sum of lowest degree terms is $\cC(\Gamma;\xx,\yy) =
(x_1-y_1)(x_2-y_2)$.%
\end{example}

The letters $\cC$ and $\KK$ stand for `cohomology' and
\hbox{`$K$-theory'}, the relation between them (`take lowest degree
terms') reflecting the Grothendieck--Riemann--Roch transition from
$K$-theory to its associated graded ring.  When $\kk$ is the complex
field~$\CC$, the (Laurent) polynomials denoted by $\cC$ and $\KK$ are
honest torus-equivariant cohomology and \mbox{$K$-classes on~$\CC^m$.}

\section{Matrix Schubert varieties}\label{sec:matSchub}

Let $\mn$ be the variety of $n \times n$ matrices over~$\kk$, with
coordinate ring $\kk[\zz]$ in indeterminates $\{z_{ij}\}_{i,j = 1}^n$.
Throughout the paper, $q$ and $p$ will be integers with $1 \leq q,p \leq
n$, and $Z$ will stand for an $n \times n$ matrix.  Most often, $Z$ will
be the \bem{generic matrix} of variables $(z_{ij})$, although
occasionally $Z$ will be an element of $\mn$.  Denote by $Z\sub qp$ the
northwest $q \times p$ submatrix of~$Z$.  For instance, given a
permutation $w \in S_n$, the permutation matrix $w^T$ with `$1$'~entries
in row~$i$ and column~$w(i)$ has upper-left $q \times p$ submatrix with
rank given by
\begin{eqnarray*}
  \rank(w^T\sub qp) &=& \#\{(i,j) \leq (q,p) \mid w(i) = j\},
\end{eqnarray*}
the number of `$1$' entries in the submatrix
$w^T\sub qp$.

The class of determinantal ideals in the following definition was
identified by Fulton in~\cite{FulDegLoc}, though in slightly different
language.

\begin{defn} \label{defn:Iw}
Let $w \in S_n$ be a permutation.  The \bem{Schubert determinantal ideal
$I_w \subset \kk[\zz]$} is generated by all minors in $Z\sub qp$ of size
$1 + \rank(w^T\sub qp)$ for all $q,p$, where $Z = (z_{ij})$ is the
matrix of variables.
\end{defn}

The subvariety of~$\mn$ cut out by~$I_w$ is the central geometric object
in this paper.

\begin{defn} \label{defn:matrixSchub}
Let $w \in S_n$.  The \bem{matrix Schubert variety} $\ol X_w \subseteq
\mn$ consists of the matrices $Z \in \mn$ such that $\rank(Z\sub qp) \leq
\rank(w^T\sub qp)$ for all $q,p$.
\end{defn}

\begin{example} \label{ex:w0}
The smallest matrix Schubert variety is $\ol X_{w_0}$, where $w_0$ is
the \bem{long permutation} \hbox{$n\,\cdots\,2\,1$} reversing the order
of $1,\ldots,n$.  The variety $\ol X_{w_0}$ is just the linear subspace
of lower-right-triangular matrices; its ideal is $\<z_{ij} \mid i+j \leq
n\>$.
\end{example}

\begin{example} \label{ex:s3}
Five of the six $3\times 3$ matrix Schubert varieties are linear
subspaces:
$$
\begin{array}{r@{\:\ =\:\ }l@{\qquad\qquad}r@{\:\ =\:\ }l}
  I_{123} & 0
  &
  \ol X_{123} & M_3
\\
  I_{213} & \<z_{11}\>
  &
  \ol X_{213} & \{Z \in M_3 \mid z_{11} = 0\}
\\
  I_{231} & \<z_{11},z_{12}\>
  &
  \ol X_{231} & \{Z \in M_3 \mid z_{11} = z_{12} = 0\}
\\
  I_{231} & \<z_{11},z_{21}\>
  &
  \ol X_{312} & \{Z \in M_3 \mid z_{11} = z_{21} = 0\}
\\
  I_{321} & \<z_{11},z_{12},z_{21}\>
  &
  \ol X_{321} & \{Z \in M_3 \mid z_{11} = z_{12} = z_{21} = 0\}
\end{array}
$$
The remaining permutation, $w = 132$, has
$$
  \;
  I_{132}\ =\ \<z_{11}z_{22} - z_{12}z_{21}\>
  \qquad\;
  \ol X_{132}\ =\ \{Z \in M_3 \mid \rank(Z\sub 22) \leq 1\},\quad\ \:
$$
so $\ol X_{132}$ is the set of matrices whose upper-left $2\times 2$
block is singular.
\end{example}

\begin{example} \label{ex:intro}
Let $w = 13865742$,
so that $w^T$ is given by replacing each $*$ by~$1$ in the left matrix
below.
$$
\begin{tinyrc}{
\begin{array}{@{}|@{\:}c@{\:}|@{\:}c@{\:}|@{\:}c@{\:}|@{\:}c@{\:}
		 |@{\:}c@{\:}|@{\:}c@{\:}|@{\:}c@{\:}|@{\:}c@{\:}|@{}}
\hline  *  &     &     &     &     &     &     &
\hln       &     &  *  &     &     &     &     &
\hln       &     &     &     &     &     &     &  *
\hln       &     &     &     &     &  *  &     &
\hln       &     &     &     &  *  &     &     &
\hln       &     &     &     &     &     &  *  &
\hln       &     &     &  *  &     &     &     &
\hln       &  *  &     &     &     &     &     &
\\\hline
\end{array}
\hbox{\normalsize$ \ \ \implies\ \ $}
\begin{array}{@{}|@{\:}c@{\:}|@{\:}c@{\:}|@{\:}c@{\:}|@{\:}c@{\:}
		 |@{\:}c@{\:}|@{\:}c@{\:}|@{\:}c@{\:}|@{\:}c@{\:}|@{}}
\hline  1  &  1  &  1  &  1  &  1  &  1  &  1  &  1
\hln    1  &  1  &     &     &     &     &     &
\hln    1  &  1  &     &     &     &     &     &
\hln    1  &  1  &     &     &     &     &     &
\hln    1  &  1  &     &     &     &     &     &
\hln    1  &  1  &     &     &     &     &     &
\hln    1  &  1  &     &     &     &     &     &
\hln    1  &     &     &     &     &     &     &
\\\hline
\end{array}
\hbox{\normalsize$ \ ,\ $}
\begin{array}{@{}|@{\:}c@{\:}|@{\:}c@{\:}|@{\:}c@{\:}|@{\:}c@{\:}
  		 |@{\:}c@{\:}|@{\:}c@{\:}|@{\:}c@{\:}|@{\:}c@{\:}|@{}}
\hline  2  &  2  &  2  &  2  &  2  &  2  &  2  &  2
\hln    2  &  2  &  2  &  2  &  2  &  2  &  2  &  2
\hln    2  &  2  &  2  &  2  &  2  &  2  &  2  &
\hln    2  &  2  &  2  &  2  &  2  &     &     &
\hln    2  &  2  &  2  &  2  &     &     &     &
\hln    2  &  2  &  2  &  2  &     &     &     &
\hln    2  &  2  &  2  &     &     &     &     &
\hln    2  &  2  &     &     &     &     &     &
\\\hline
\end{array}
\hbox{\normalsize$ \ ,\ $}
\begin{array}{@{}|@{\:}c@{\:}|@{\:}c@{\:}|@{\:}c@{\:}|@{\:}c@{\:}
  		 |@{\:}c@{\:}|@{\:}c@{\:}|@{\:}c@{\:}|@{\:}c@{\:}|@{}}
\hline  3  &  3  &  3  &  3  &  3  &  3  &  3  &  3
\hln    3  &  3  &  3  &  3  &  3  &  3  &  3  &  3
\hln    3  &  3  &  3  &  3  &  3  &  3  &  3  &  3
\hln    3  &  3  &  3  &  3  &  3  &  3  &  3  &
\hln    3  &  3  &  3  &  3  &  3  &     &     &
\hln    3  &  3  &  3  &  3  &  3  &     &     &
\hln    3  &  3  &  3  &  3  &     &     &     &
\hln    3  &  3  &  3  &     &     &     &     &
\\\hline
\end{array}
\hbox{\normalsize$ \ ,\ $}
\begin{array}{@{}|@{\:}c@{\:}|@{\:}c@{\:}|@{\:}c@{\:}|@{\:}c@{\:}
  		 |@{\:}c@{\:}|@{\:}c@{\:}|@{\:}c@{\:}|@{\:}c@{\:}|@{}}
\hline  4  &  4  &  4  &  4  &  4  &  4  &  4  &  4
\hln    4  &  4  &  4  &  4  &  4  &  4  &  4  &  4
\hln    4  &  4  &  4  &  4  &  4  &  4  &  4  &  4
\hln    4  &  4  &  4  &  4  &  4  &  4  &  4  &  4
\hln    4  &  4  &  4  &  4  &  4  &  4  &  4  &
\hln    4  &  4  &  4  &  4  &  4  &  4  &     &
\hln    4  &  4  &  4  &  4  &  4  &     &     &
\hln    4  &  4  &  4  &  4  &     &     &     &
\\\hline
\end{array}
\hbox{\normalsize$ \ ,\ \ldots$}
}\end{tinyrc}
$$
Each matrix in $\ol X_w \subseteq \mn$ has the property that every
rectangular submatrix contained in the region filled with $1$'s has rank
$\leq 1$, and every rectangular submatrix contained in the region filled
with $2$'s has rank $\leq 2$, and so on.  The ideal $I_w$ therefore
contains the $21$ minors of size $2 \times 2$ in the first region and the
$144$
\begin{excise}%
{$$
\textstyle
  [{7 \choose 3} \cdot 1 + {6 \choose 3} \cdot {4 \choose 3} + {4 \choose
  3} \cdot {5 \choose 3} + 1 \cdot {7 \choose 3}] - [{6 \choose 3} \cdot
  1 + {4 \choose 3} \cdot {4 \choose 3} + 1 \cdot {5 \choose 3}]
$$
$$
\textstyle
  [35 + 20 \cdot 4 + 4 \cdot 10 + 35] - [20 + 4 \cdot 4 + 10] =
  [35 + 80 + 40 + 35] - [20 + 16 + 10] = 190 - 46 = 144
$$
$$
\textstyle
  [{7 \choose 3}] + [{6 \choose 3} \cdot {3 \choose 2}] + [{4 \choose 3}
  \cdot {4 \choose 2}] + [2 \cdot 1 \cdot {5 \choose 2} + 1 \cdot 1
  \cdot {5 \choose 1}]
$$
$$
\textstyle
  [35] + [20 \cdot 3] + [4 \cdot 6] + [20 + 5] = 35 + 60 + 24 + 25 = 144
$$}%
\end{excise}%
minors of size $3 \times 3$ in the second region.  These 165 minors in
fact generate~$I_w$, as can be checked either directly by Laplace
expansion of each determinant in~$I_w$ along its last row(s) or
column(s), or indirectly using Fulton's notion of `essential set'
\cite{FulDegLoc}.
See also Example~\ref{ex:gb}.%
\end{example}

Our first main theorem provides a straightforward geometric explanation
for the naturality of Schubert and Grothendieck polynomials.  More
precisely, our context automatically makes them well-defined as
(Laurent) polynomials, as opposed to being identified as (particularly
nice) representatives for classes in some quotient of a polynomial ring.

\begin{theorem} \label{t:formulae}
The Schubert determinantal ideal~$I_w$ is prime, so $I_w$ is the ideal
$I(\ol X_w)$ of the matrix Schubert variety~$\ol X_w$.  The
$\ZZ^n$-graded and $\ZZ^{2n}$-graded \K-polynomials
of\/~%
$\ol X_w$ are the Grothendieck and double Grothendieck polynomials
for~$w$, respectively:
\begin{eqnarray*}
  \ \;
  \KK(\ol X_w; \xx)\:=\:\GG_w(\xx)
  \:&{\rm and}&\:
  \KK(\ol X_w;\xx,\yy)\:=\:\GG_w(\xx,\yy).
\end{eqnarray*}
The $\ZZ^n$-graded and $\ZZ^{2n}$-graded multidegrees of $\ol X_w$ are
the Schubert and double Schubert polynomials for~$w$, respectively:
\begin{eqnarray*}
  [\ol X_w]_{\ZZ^n}\:=\:\SS_w(\xx)
  \:&{\rm and}&\:
  [\ol X_w]_{\ZZ^{2n}}\:=\:\SS_w(\xx,\yy),
\end{eqnarray*}
\end{theorem}

Primality of~$I_w$ was proved by Fulton \cite{FulDegLoc}, but we shall
not assume it in our proofs.

\begin{example} \label{ex:formulae}
Let $w = 2143$ as in the example from the Introduction.  Computing the
\K-polynomial of the complete intersection $\ci{2143}$ yields (in the
$\ZZ^n$-grading for simplicity)
$$
  (1-x_1)(1-x_1x_2x_3)\:\ =\:\ \GG_{2143}(\xx)\:\ =\:\ \dem2 \dem 1 \dem3
  \dem2\big((1-x_1)^3(1-x_2)^2(1-x_3)\big),
$$
the latter equality by Theorem~\ref{t:formulae}.
Substituting $\xx \mapsto \1-\xx$ in $\GG_{2143}(\xx)$ yields
\begin{eqnarray*}
  \GG_{2143}(\1-\xx) &=& x_1(x_1+x_2+x_3-x_1x_2-x_2x_3-x_1x_3+x_1x_2x_3),
\end{eqnarray*}
whose sum of lowest degree terms equals the multidegree $\cC(\ol
X_{2143};\xx)$ by definition.  This agrees with the Schubert polynomial
$\SS_{2143}(\xx) = x_1^2 + x_1x_2 + x_1x_3$.%
\end{example}

That Schubert and Grothendieck polynomials represent cohomology and
\K-theory classes of Schubert varieties in flag manifolds will be shown
in Section~\ref{sec:flag} to follow from Theorem~\ref{t:formulae}.

\section{Pipe dreams}\label{sec:pipe}

In this section we introduce the set $\rc(w)$ of reduced pipe dreams%
	\footnote{In the game Pipe Dream, the player is supposed to
	guide water flowing out of a spigot at one edge of the game
	board to its destination at another edge by laying down given
	square tiles with pipes going through them; see
	Definition~\ref{defn:pipe}.  The spigot placements and
	destinations are interpreted in Definition~\ref{defn:rp}.}
for a permutation $w \in S_n$.  Each diagram $D \in \rc(w)$ is a subset
of the $n \times n$ grid $[n]^2$ that represents an example of the curve
diagrams invented by Fomin and Kirillov \cite{FKyangBax}, though our
notation follows Bergeron and Billey~\cite{BB} in this regard.%
	\footnote{The corresponding objects in \cite{FKyangBax} look
	like reduced pipe dreams rotated by $135^\circ$.}
Besides being attractive ways to draw permutations, reduced pipe dreams
generalize to flag manifolds the semistandard Young tableaux for
Grassmannians.  Indeed, there is even a natural bijection between
tableaux and reduced pipe dreams for Grassmannian permutations (see
\cite{KoganThesis}, for instance).

Consider a square grid $\ZZ_{> 0} \times \ZZ_{> 0}$ extending infinitely
south and east, with the box in row~$i$ and column~$j$ labeled $(i,j)$,
as in an $\infty \times \infty$ matrix.  If each box in the grid is
covered with a square tile containing either $\textcross$ or
$\textelbow$, then one can think of the tiled grid as a network of pipes.

\begin{defn} \label{defn:pipe}
A \bem{pipe dream}
is a finite subset of $\ZZ_{> 0} \times \ZZ_{> 0}$, identified as the set
of crosses in a tiling by \bem{crosses} $\textcross$ and \bem{elbow
joints} $\textelbow$.%
\end{defn}

Whenever we draw pipe dreams, we fill the boxes with crossing tiles by
\cross.  However, we often leave the elbow tiles blank, or denote them by
dots for ease of notation.  The pipe dreams we consider all represent
subsets of the pipe dream~$D_0$ that has crosses in the triangular region
strictly above the main antidiagonal (in spots $(i,j)$ with $i+j \leq n$)
and elbow joints elsewhere.  Thus we can safely limit ourselves to
drawing inside $n \times n$ grids.

\begin{example} \label{ex:pipe}
Here are two rather arbitrary pipe dreams with $n=5$:
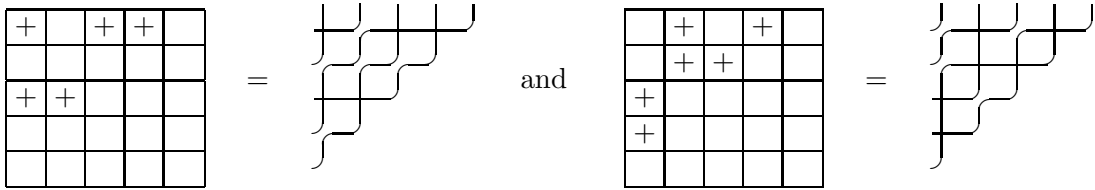
\vspace{-2ex}
$$
\begin{array}{|c|c|c|c|c|}
\multicolumn{5}{c}{}
\\\hline\!+\!&     &\!+\!&\!+\!&\ \,
\\\hline     &     &     &     &
\\\hline\!+\!&\!+\!&     &     &
\\\hline     &     &     &     &
\\\hline     &     &     &     &
\\\hline
\end{array}
\,\quad =
\begin{array}{cccccc}
&\phperm &\phperm &\phperm &\phperm &\phperm \\
&   \+   &   \jr  &   \+   &   \+   &  \je   \\
&   \jr  &   \jr  &   \jr  &   \je  &\\
&   \+   &   \+   &   \je  &        &\\
&   \jr  &   \je  &        &        &\\
&   \je  &        &        &        &
\end{array}
\quad \hbox{and} \qquad
\begin{array}{|c|c|c|c|c|}
\multicolumn{5}{c}{}
\\\hline     &\!+\!&     &\!+\!&\ \,
\\\hline     &\!+\!&\!+\!&     &
\\\hline\!+\!&     &     &     &
\\\hline\!+\!&     &     &     &
\\\hline     &     &     &     &
\\\hline
\end{array}
\,\quad =
\begin{array}{cccccc}
&\phperm &\phperm &\phperm &\phperm &\phperm \\
&   \jr  &   \+   &   \jr  &   \+   &  \je   \\
&   \jr  &   \+   &   \+   &   \je  &\\
&   \+   &   \jr  &   \je  &        &\\
&   \+   &   \je  &        &        &\\
&   \je  &        &        &        &
\end{array}
$$
Another (slightly less arbitrary) example, with $n=8$, is the pipe
dream~$D$ in Fig.~\ref{fig:pipe}.
\begin{figure}[ht]
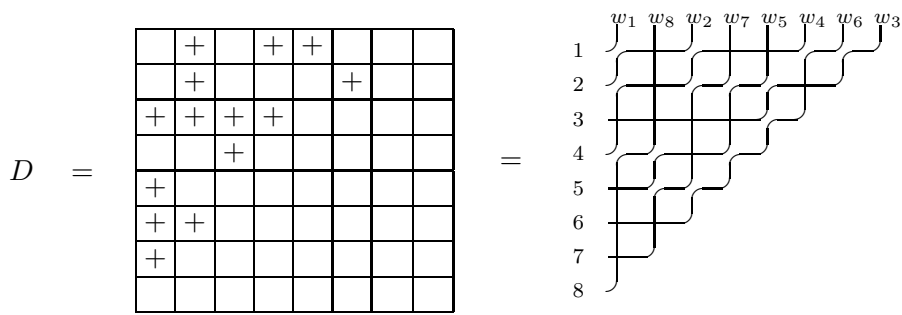

$$
\begin{array}{c}\\[-4ex]
	\\\\\\\\\\\raisebox{2ex}{$D\quad=\quad$}\\\\\\\\\end{array}
\begin{array}{|c|c|c|c|c|c|c|c|}
\multicolumn{8}{c}{}\\[-4ex]
\multicolumn{8}{c}{}
\\\hline     &\!+\!&     &\!+\!&\!+\!&     &\ \, &\ \,
\\\hline     &\!+\!&     &     &     &\!+\!&     &
\\\hline\!+\!&\!+\!&\!+\!&\!+\!&     &     &     &
\\\hline     &     &\!+\!&     &     &     &     &
\\\hline\!+\!&     &     &     &     &     &     &
\\\hline\!+\!&\!+\!&     &     &     &     &     &
\\\hline\!+\!&     &     &     &     &     &     &
\\\hline     &     &     &     &     &     &     &
\\\hline
\end{array}
\ \quad = \quad
\begin{array}{ccccccccc}
\\[-4ex]
     &\perm w1&\perm w8&\perm w2&\perm w7&\perm w5&\perm w4&\perm w6&\perm w3\\
\petit1&  \jr &   \+   &   \jr  &   \+   &   \+   &   \jr  &   \jr  &  \je\\
\petit2&  \jr &   \+   &   \jr  &   \jr  &   \jr  &   \+   &   \je  &\\
\petit3&  \+  &   \+   &   \+   &   \+   &   \jr  &   \je  &        &\\
\petit4&  \jr &   \jr  &   \+   &   \jr  &   \je  &        &        &\\
\petit5&  \+  &   \jr  &   \jr  &   \je  &        &        &        &\\
\petit6&  \+  &   \+   &   \je  &        &        &        &        &\\
\petit7&  \+  &   \je  &        &        &        &        &        &\\
\petit8&  \je &        &        &        &        &        &        &
\end{array}
$$
\caption{A pipe dream with $n=8$} \label{fig:pipe}
\end{figure}
The first diagram represents $D$ as a subset of $[8]^2$, whereas the
second demonstrates how the tiles fit together.  Since no cross in~$D$
occurs on or below the $8^\th$ antidiagonal, the pipe entering row~$i$
exits column~$w_i = w(i)$ for some permutation $w \in S_8$.  In this
case, $w = 13865742$ is the permutation from Example~\ref{ex:intro}.
For clarity, we omit the square tile boundaries as well as the wavy
``sea'' of elbows below the main antidiagonal in the right pipe dream.
We also use the thinner symbol $w_i$ instead of $w(i)$ to make the
column widths come out right.%
\end{example}

\begin{defn}\label{defn:rp}
A pipe dream is \bem{reduced} if each pair of pipes crosses at most
once.  The set $\rc(w)$ of reduced pipe dreams for the permutation $w
\in S_n$ is the set of reduced pipe dreams~$D$ such that the pipe
entering row~$i$ exits from column $w(i)$.
\end{defn}

We shall give some idea of what it means for a pipe dream to be reduced,
in Lemma~\ref{lemma:subword}, below.  For notation, we say that a \cross
at $(q,p)$ in a pipe dream $D$ sits on the \bem{$i^\th$ antidiagonal} if
$q+p-1 = i$.  Let $Q(D)$ be the ordered sequence of simple
reflections~$s_i$ corresponding to the antidiagonals on which the
crosses sit, starting from the northeast corner of $D$ and reading {\em
right to left}\/ in each row, snaking down to the southwest corner.%
	\footnote{The term `rc-graph' was used in \cite{BB} for what we
	call reduced pipe dreams.  The letters~`rc' stand for
	``reduced-compatible''.  The ordered list of row indices for the
	crosses in~$D$, taken in the same order as before, is called in
	\cite{BJS} a ``compatible sequence'' for the expression~$Q(D)$;
	we shall not need this concept.}

\begin{example} \label{ex:subword}
The pipe dream $D_0$ corresponds to the ordered sequence 
\begin{eqnarray*}
  Q(D_0)\ \:=\:\ Q_0 &:=& s_{n-1} \cdots s_2 s_1\ s_{n-1} \cdots s_3
  s_2\ \cdots\cdots\ s_{n-1} s_{n-2} s_{n-1},
\end{eqnarray*}
the \bem{triangular} reduced expression for the long permutation~$w_0 =
n \cdots 321$.  Thus $Q_0 = s_3s_2s_1s_3s_2s_3$ when $n=4$.  For another
example, the first pipe dream in Example~\ref{ex:pipe} yields the
ordered sequence $s_4s_3s_1s_5s_4$.
\end{example}

\begin{lemma} \label{lemma:subword}
If $D$ is a pipe dream, then multiplying the reflections in $Q(D)$
yields the permutation $w$ such that the pipe entering row~$i$ exits
column~$w(i)$.  Furthermore, the number of crossing tiles in~$D$ is at
least\/ $\length(w)$, with equality if and only if $D \in \rc(w)$.
\end{lemma}
\begin{proof}
For the first statement, use induction on the number of crosses: adding
a \cross in the $i^\th$ antidiagonal at the end of the list switches the
destinations of the pipes beginning in rows~$i$ and~$i+1$.  Each
inversion in~$w$ contributes at least one crossing in $D$, whence the
number of crossing tiles is at least $\length(w)$.  The expression
$Q(D)$~is reduced when $D$~is reduced because each inversion in~$w$
contributes at most one crossing tile to~$D$.
\end{proof}

In other words, pipe dreams with no crossing tiles on or below the main
antidiagonal in~$[n]^2$ are naturally `subwords' of~$Q(D_0)$, while
reduced pipe dreams are naturally {\em reduced}\/ subwords.  This point
of view
takes center \mbox{stage in Section~\ref{sec:subword}}.

\begin{example} \label{ex:rp}
The upper-left triangular pipe dream $D_0 \subset [n]^2$ is the unique
pipe dream in~$\rc(w_0)$.  The $8 \times 8$ pipe dream~$D$ in
Example~\ref{ex:pipe} lies in~$\rc(13865742)$.
\end{example}

\section{Gr\"obner geometry}\label{sec:grobGeom}

Using Gr\"obner bases, we next degenerate matrix Schubert varieties into
unions of vector subspaces of~$\mn$ corresponding to reduced pipe
dreams.  A~total order~`$>$' on monomials in~$\kk[\zz]$ is a \bem{term
order} if $1 \leq m$ for all monomials $m \in \kk[\zz]$, and $m \cdot m'
< m \cdot m''$ whenever $m' < m''$.  When a term order~`$>$' is fixed,
the largest monomial~$\IN(f)$ appearing with nonzero coefficient in a
polynomial~$f$ is its \bem{initial term}, and the \bem{initial ideal} of
a given ideal~$I$ is generated by the initial terms of all polynomials
$f \in I$.  A~set $\{f_1,\ldots,f_n\}$ is a \bem{Gr\"obner basis} if
$\IN(I) = \<\IN(f_1),\ldots,\IN(f_n)\>$.  See \cite[Chapter~15]{Eis} for
background on term orders and Gr\"obner bases, including geometric
interpretations in terms of flat families.

\begin{defn}
The \bem{antidiagonal ideal} $J_w$ is generated by the antidiagonals of
the minors of $Z = (z_{ij})$ generating $I_w$.  Here, the
\bem{antidiagonal} of a square matrix or a minor is the product of the
entries on the main antidiagonal.
\end{defn}

There exist numerous \bem{antidiagonal term orders} on $\kk[\zz]$, which
by definition pick off from each minor its antidiagonal term, including:
\begin{itemize}
\item
the reverse lexicographic term order that snakes its way from the
northwest corner to the southeast corner, $z_{11} > z_{12} > \cdots >
z_{1n} > z_{21} > \cdots > z_{nn}$; and

\item
the lexicographic term order that snakes its way from northeast corner
to the southwest corner, $z_{1n} > \cdots > z_{nn} > \cdots > z_{2n} >
z_{11} > \cdots > z_{n1}$.
\end{itemize}
The initial ideal $\IN(I_w)$ for any antidiagonal term order
contains~$J_w$ by definition, and our first point in Theorem~\ref{t:gb}
will be equality of these two monomial ideals.

Our remaining points in Theorem~\ref{t:gb} concern the combinatorics
of~$J_w$.  Being a squarefree monomial ideal, it is by definition the
\bem{Stanley--Reisner ideal} of some simplicial complex~$\LL_w$ with
vertex set $[n]^2 = \{(q,p) \mid 1 \leq q,p \leq n\}$.  That is, $\LL_w$
consists of the subsets of~$[n]^2$ containing no antidiagonal in~$J_w$.
Faces of~$\LL_w$ (or any simplicial complex with~$[n]^2$ for vertex set)
may be identified with coordinate subspaces in~%
$\mn$ as follows.  Let $E_{qp}$ denote the elementary matrix whose only
nonzero entry lies in row~$q$ and column~$p$, and identify vertices
in~$[n]^2$ with variables $z_{qp}$ in the generic matrix~$Z$.  Letting
$D_L = [n]^2 \minus L$ be the pipe dream complementary to~$L$, each
face~$L$ is identified with the coordinate subspace
\begin{eqnarray*}
  L &=& \{z_{qp} = 0 \mid (q,p) \in D_L\}\ \:=\ \:
  {\rm span}(E_{qp} \mid (q,p) \not\in D_L).
\end{eqnarray*}
Thus, considering $D_L$ as a pipe dream, its crosses $\textcross$ lie in
the spots where $L$ is zero.  For instance, the three pipe dreams in the
example from the Introduction
are pipe dreams for the subspaces $L_{11,13}$, $L_{11,22}$,
and~$L_{11,31}$.

The term \bem{facet} means `maximal face', and
Definition~\ref{defn:vertex} gives the meaning of~`shellable'.\!

\begin{theorem} \label{t:gb}
The minors of size $1+\rank(w^T\sub qp)$ in $Z\sub qp$ for all $q,p$
constitute a Gr\"obner basis for any antidiagonal term order;
equivalently, $\IN(I_w) = J_w$ for any such term order.
The Stanley--Reisner complex~$\LL_w$ of~$J_w$ is shellable, and hence
Cohen--Macaulay.  In addition,
\begin{eqnarray*}
  \{D_L \mid L \hbox{ is a facet of } \LL_w\} &=& \rc(w)
\end{eqnarray*}
places the set of reduced pipe dreams for~$w$ in canonical bijection
with the facets of~$\LL_w$.
\end{theorem}

The displayed equation is equivalent to~$J_w$ having the prime
decomposition
\begin{eqnarray*}
  J_w &=& \!\bigcap_{D \in \rc(w)} \<z_{ij} \mid (i,j) \in D\>.
\end{eqnarray*}
	
Geometrically, Theorem~\ref{t:gb} says that the matrix Schubert
variety~$\ol X_w$ has a flat degeneration whose limit is both reduced
and Cohen--Macaulay, and whose components are in natural bijection with
reduced pipe dreams.  On its own, Theorem~\ref{t:gb} therefore ascribes
a truly geometric origin to reduced pipe dreams.
Taken together with Theorem~\ref{t:formulae}, it provides in addition a
natural geometric explanation for the combinatorial formulae writing
Schubert
polynomials in terms of pipe dreams: interpret in equivariant cohomology
the decomposition of~$\LL_w$ into irreducible components.  This
procedure is carried out in Section~\ref{sec:bjs} using multidegrees,
for which the required technology is developed in Section~\ref{sec:pos}.
The analogous \K-theoretic formula, which additionally involves {\em
nonreduced}\/ pipe dreams, requires more detailed analysis of subword
complexes (Definition~\ref{defn:subword}), and therefore appears
in~\cite{KMsubword}.

\begin{example} \label{ex:2143'}
Let $w = 2143$ as in the example from the Introduction and
Example~\ref{ex:formulae}.
The term orders that interest us pick out the antidiagonal term $-z_{13}
z_{22} z_{31}$ from the northwest $3 \times 3$ minor.  For $I_{2143}$,
this causes the initial terms of its two generating minors to be
relatively prime, so the minors form a Gr\"obner basis as in
Theorem~\ref{t:gb}.  Observe that the minors generating $I_w$ do not
form a Gr\"obner basis with respect to term orders that pick out the
{\em diagonal}\/ term $z_{11} z_{22} z_{33}$ of the $3 \times 3$ minor,
because $z_{11}$ divides that.

The initial complex $\LL_{2143}$ is shellable, being a cone over the
boundary of a triangle, and as mentioned in the Introduction,
its facets correspond to the reduced pipe dreams for~$2143$.%
\end{example}

\begin{example} \label{ex:gb}
A direct check reveals that every antidiagonal in $J_w$ for $w =
13865742$ stipulated by Definition~\ref{defn:Iw} is divisible by an
antidiagonal of some $2$- or $3$-minor from Example~\ref{ex:intro}.
Hence the $165$ minors of size $2 \times 2$ and $3 \times 3$ in $I_w$
\mbox{form a Gr\"obner basis for~$I_w$}.%
\end{example}

\begin{remark}
M.\thinspace{}Kogan also has a geometric interpretation for reduced pipe
dreams, identifiying them in~\cite{KoganThesis} as subsets of the flag
manifold mapping to corresponding faces of the Gel$'$fand--Cetlin
polytope.  These subsets are not cycles, so they do not individually
determine cohomology classes whose sum is the Schubert class;
nonetheless, their union is a cycle, and its class is the Schubert
class.  See also \cite{KoganM}.%
\end{remark}

\begin{remark}
Theorem~\ref{t:gb} says that every antidiagonal shares at least one
cross with every reduced pipe dream, and moreover, that each
antidiagonal and reduced pipe dream is minimal with this property.
Loosely, antidiagonals and reduced pipe dreams `minimally poison' each
other.  Our proof of this purely combinatorial statement in
Sections~\ref{sec:facets} and~\ref{sec:rp} is indeed essentially
combinatorial, but rather roundabout; we know of no simple
\mbox{reason for~it}.%
\end{remark}

\begin{remark}
The Gr\"obner basis in Theorem~\ref{t:gb} defines a flat degeneration
over any ring, because all of the coefficients of the minors in $I_w$
are integers, and the leading coefficients are all $\pm 1$.  Indeed,
each loop of the division algorithm in Buchberger's criterion
\cite[Theorem~15.8]{Eis} works over $\ZZ$, and therefore over any ring.
\end{remark}

\section{Mitosis algorithm}\label{sec:alg}

Next we introduce a simple combinatorial rule, called `mitosis',%
        \footnote{The term \bem{mitosis} is biological lingo for cell
        division in multicellular organisms.}
%
that creates from each pipe dream a number of new pipe dreams called its
`offspring'.  Mitosis serves as a geometrically motivated improvement on
Kohnert's rule \cite{Kohnert, NoteSchubPoly, Winkel}, which acts on
other subsets of $[n]^2$ derived from
permutation matrices.  In addition to its independent interest from a
combinatorial standpoint, our forthcoming Theorem~\ref{t:mitosis} falls
out of Bruhat induction with no extra work, and in fact the mitosis
operation plays a vital role in Bruhat induction,
toward the end of Part~\ref{part:gb}.

Given a pipe dream in $[n] \times [n]$, define
\begin{eqnarray} \label{eq:pipestart}
  \start_i(D) &=& \hbox{column index of leftmost empty box in row } i
\\[-.5ex]\nonumber
  &=& \min(\{j \mid (i,j) \not\in D\} \cup \{n+1\}).
\end{eqnarray}
Thus in the region to the left of $\start_i(D)$, the $i^\th$ row of $D$
is filled solidly with crosses.  Let
\begin{eqnarray*}
  \JJ_i(D) &=& \{\hbox{columns $j$ strictly to the left of } \start_i(D)
  \mid (i+1,j) \hbox{ has no cross in } D\}.
\end{eqnarray*}
For $p \in \JJ_i(D)$, construct the \bem{offspring} $D_p(i)$ as follows.
First delete the cross at $(i,p)$ from~$D$.  Then take all crosses in
row~$i$ of $\JJ_i(D)$ that are to the left of column~$p$, and move each
one down to the empty box below it in row~$i+1$.

\begin{defn} \label{defn:mitosis}
The $i^\th$ \bem{mitosis} operator sends a pipe dream $D$ to
\begin{eqnarray*}
  \mitosis_i(D) &=& \{D_p(i) \mid p \in \JJ_i(D)\}.
\end{eqnarray*}
Thus all the action takes place in rows~$i$ and~$i+1$, and
$\mitosis_i(D)$ is an empty set if~$\JJ_i(D)$~is empty.  Write
$\mitosis_i({\mathcal P}) = \bigcup_{D \in {\mathcal P}}
\mitosis_i({\mathcal D})$ whenever ${\mathcal P}$ is a set of pipe
dreams.%
\end{defn}


\begin{example} \label{ex:mitosis}
The left diagram $D$ below is the reduced pipe dream for $w = 13865742$
from Example~\ref{ex:pipe} (the pipe dream in Fig.~\ref{fig:pipe})
and Example~\ref{ex:rp}:
$$
\begin{tinyrc}{
\begin{array}{@{}cc@{}}{}\\\\3&\\4\\\\\\\\\\\\\\\end{array}
\begin{array}{@{}|@{\,}c@{\,}|@{\,}c@{\,}|@{\,}c@{\,}|@{\,}c@{\,}
		 |@{\,}c@{\,}|@{\,}c@{\,}|@{\,}c@{\,}|@{\,}c@{\,}|@{}}
\hline     &  +  &     &  +  &  +  &     &     &
\hln       &  +  &     &     &     &  +  &     &\phantom{+}
\hln    +  &  +  &  +  &  +  &     &     &     &
\hln       &     &  +  &     &     &     &     &
\hln    +  &     &     &     &     &     &     &
\hln    +  &  +  &     &     &     &     &     &
\hln    +  &     &     &     &     &     &     &
\hln       &     &     &     &     &     &\phantom{+}&
\\\hline\multicolumn{4}{@{}c@{}}{}&\multicolumn{1}{@{}c@{}}{\!\!\uparrow}
\\\multicolumn{4}{c}{}&\multicolumn{1}{@{}c@{}}{\makebox[0pt]{$\start_3$}}
\end{array}
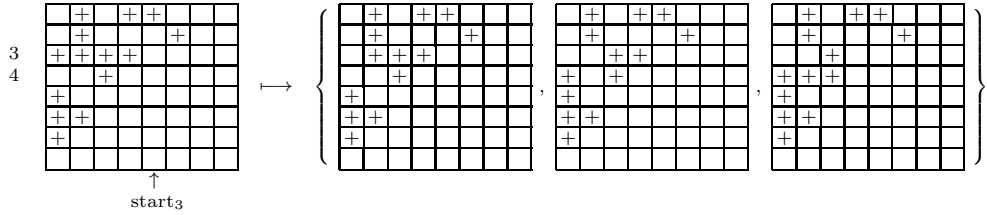
\begin{array}{@{\quad}c@{\ }}
  \longmapsto\\\mbox{}\\\mbox{}
\end{array}
\begin{array}{c}
  \left\{\ 
  \begin{array}{@{}|@{\,}c@{\,}|@{\,}c@{\,}|@{\,}c@{\,}|@{\,}c@{\,}
  		 |@{\,}c@{\,}|@{\,}c@{\,}|@{\,}c@{\,}|@{\,}c@{\,}|@{}}
  \hline     &  +  &     &  +  &  +  &     &     &    
  \hln       &  +  &     &     &     &  +  &     &\phantom{+}
  \hln       &  +  &  +  &  +  &     &     &     &
  \hln       &     &  +  &     &     &     &     &
  \hln    +  &     &     &     &     &     &     &
  \hln    +  &  +  &     &     &     &     &     &
  \hln    +  &     &     &     &     &     &     &
  \hln       &     &     &     &     &     &\phantom{+}&
  \\\hline
  \end{array}
  \ ,\ 
  \begin{array}{@{}|@{\,}c@{\,}|@{\,}c@{\,}|@{\,}c@{\,}|@{\,}c@{\,}
  		 |@{\,}c@{\,}|@{\,}c@{\,}|@{\,}c@{\,}|@{\,}c@{\,}|@{}}
  \hline     &  +  &     &  +  &  +  &     &     &    
  \hln       &  +  &     &     &     &  +  &     &\phantom{+}
  \hln       &     &  +  &  +  &     &     &     &
  \hln    +  &     &  +  &     &     &     &     &
  \hln    +  &     &     &     &     &     &     &
  \hln    +  &  +  &     &     &     &     &     &
  \hln    +  &     &     &     &     &     &     &
  \hln       &     &     &     &     &     &\phantom{+}&
  \\\hline
  \end{array}
  \ ,\ 
  \begin{array}{@{}|@{\,}c@{\,}|@{\,}c@{\,}|@{\,}c@{\,}|@{\,}c@{\,}
  		 |@{\,}c@{\,}|@{\,}c@{\,}|@{\,}c@{\,}|@{\,}c@{\,}|@{}}
  \hline     &  +  &     &  +  &  +  &     &     &    
  \hln       &  +  &     &     &     &  +  &     &\phantom{+}
  \hln       &     &  +  &     &     &     &     &
  \hln    +  &  +  &  +  &     &     &     &     &
  \hln    +  &     &     &     &     &     &     &
  \hln    +  &  +  &     &     &     &     &     &
  \hln    +  &     &     &     &     &     &     &
  \hln       &     &     &     &     &     &\phantom{+}&
  \\\hline
  \end{array}
  \ \right\}
\\
  \begin{array}{c}
  \mbox{}\\
  \mbox{}\\
  \end{array}
\end{array}
}\end{tinyrc}
$$
The set of three pipe dreams on the right is obtained by applying
$\mitosis_3$, since $\JJ_3(D)$ consists of columns $1$, $2$, and
$4$.%
\end{example}

\begin{theorem} \label{t:mitosis}
If\/ $\length(ws_i) < \length(w)$, then $\rc(ws_i)$ is equal to the
disjoint union $\bigcupdot_{D \in \rc(w)} \mitosis_i(D)$.  Thus if
$s_{i_1} \cdots s_{i_k}$ is a reduced expression for~$w_0w$, and $D_0$
is the unique reduced pipe dream for~$w_0$, in which every entry above
the antidiagonal is a~{\rm\cross\!\!},~then
\begin{eqnarray*}
  \rc(w) &=& \mitosis_{i_k} \cdots \mitosis_{i_1}(D_0).
\end{eqnarray*}
\end{theorem}

Readers wishing a simple and purely combinatorial proof that avoids
Bruhat induction as in Part~\ref{part:gb} should
consult~\cite{mitosis}; the proof there uses only definitions and the
statement of Corollary~\ref{cor:BJS}, below, which has elementary
combinatorial proofs.  However, granting Theorem~\ref{t:mitosis} does
not by itself simplify the arguments in Part~\ref{part:gb} here: we
still need the `lifted Demazure operators' from Section~\ref{sec:lift},
of which mitosis is a distilled residue.

\begin{example} \label{ex:t:mitosis}
The left pipe dream in Example~\ref{ex:mitosis} lies in $\rc(13865742)$.
Therefore the three diagrams on the right hand side of
Example~\ref{ex:mitosis} are reduced pipe dreams for $13685742 =
13865742 \cdot s_3$ by Theorem~\ref{t:mitosis}, as can also be checked
directly.
\end{example}

Like Kohnert's rule, mitosis is inductive on weak Bruhat order, starts
with subsets of~$[n]^2$ naturally associated to the permutations
in~$S_n$,
and produces more subsets of~$[n]^2$.  Unlike Kohnert's rule, however,
the offspring of mitosis still lie in the same natural
set~as~the~parent, and the algorithm in Theorem~\ref{t:mitosis} for
generating $\rc(w)$ is irredundant, in the sense that each reduced pipe
dream appears exactly once in the implicit union on the right hand side
of the equation in Theorem~\ref{t:mitosis}.
See~\cite{mitosis} for more on properties of the mitosis recursion and
structures on the set of reduced pipe dreams, as well as background on
other
combinatorial algorithms for coefficients of Schubert polynomials.

\section{Positivity of multidegrees}\label{sec:pos}

The key to our view of positivity, which we state in
Theorem~\ref{t:pos}, lies in three properties of multidegrees
(Theorem~\ref{t:multidegs}) that characterize them uniquely among
functions on multigraded modules.  Since the multigradings considered
here are \bem{positive}, meaning that every graded piece of~$\kk[\zz]$
(and hence every graded piece of every finitely generated graded module)
has finite dimension as a vector space over the field~$\kk$, we are able
to present short complete proofs of the required assertions.

In this section we resume the generality and notation concerning
multigradings from Section~\ref{sec:mult}.  Given a (reduced and
irreducible) variety~$X$ and a module~$\Gamma$ over $\kk[\zz]$, let
$\mult_X(\Gamma)$ denote the \bem{multiplicity} of $\Gamma$ along~$X$,
which by definition equals the length of the largest finite-length
submodule in the localization of~$\Gamma$ at the prime ideal of~$X$.
The \bem{support} of~$\Gamma$ consists of those points at which the
localization of~$\Gamma$ is nonzero.

\begin{thm} \label{t:multidegs}
The multidegree \mbox{$\Gamma \mapsto \cC(\Gamma;\tt)$} is uniquely
characterized among functions from the class of\/ finitely generated
$\ZZ^d$-graded modules to\/~$\ZZ[\tt]$ by the following.
\begin{itemize}
\item
{\rm Additivity:} The (automatically $\ZZ^d$-graded) irreducible
components $X_1,\ldots,X_r$ of maximal dimension in the support of a
module\/~$\Gamma$
satisfy
\begin{eqnarray*}
  \cC(\Gamma;\tt) &=& \sum_{\ell=1}^r\mult_{X_\ell}(\Gamma) \cdot
  \cC(X_\ell,\tt).
\end{eqnarray*}

\item
{\rm Degeneration:} Let~$u$ be a variable of ordinary weight zero.  If a
finitely generated $\ZZ^d$-graded module over $\kk[\zz][u]$ is flat
over~$\kk[u]$ and has $u = 1$ fiber isomorphic to\/~$\Gamma$, then its
$u = 0$ fiber~$\Gamma'$ has the same multidegree as~$\Gamma$ does:
\begin{eqnarray*}
  \cC(\Gamma;\tt) &=& \cC(\Gamma';\tt).
\end{eqnarray*}

\item
{\rm Normalization:} If\/ $\Gamma = \kk[\zz]/\!\<z_i \mid i \in D\>$ is
the coordinate ring of a coordinate subspace of\/~$\kk^m$ for some
subset $D \subseteq \{1,\ldots,m\}$, then
\begin{eqnarray*}
  \cC(\Gamma;\tt) &=& \prod_{i\in D}\Big(\sum_{j=1}^d a_{ij}t_j\Big)
\end{eqnarray*}
is the corresponding product of ordinary weights in\/~$\ZZ[\tt] =
\sym^\spot_\ZZ(\ZZ^d)$.
\end{itemize}
\end{thm}
\begin{proof}
For uniqueness, first observe that every finitely generated
$\ZZ^d$-graded module $\Gamma$ can be degenerated via Gr\"obner bases to
a module $\Gamma'$ supported on a union of coordinate subspaces
\cite[Chapter~15]{Eis}.  By degeneration the module~$\Gamma'$ has the
same multidegree; by additivity the multidegree of~$\Gamma'$ is
determined by the multidegrees of coordinate subpaces; and by
normalization the multidegrees of coordinate subpaces are fixed.

Now we must prove that multidegrees satisfy the three conditions.
Degeneration is easy: since we have assumed the grading to be positive,
$\ZZ^d$-graded modules have $\ZZ^d$-graded Hilbert series, which are
constant in flat families of multigraded modules.

Normalization involves a bit of calculation.  Using the Koszul complex,
the \K-polynomial of $\kk[\zz]/\!\<z_i \mid i \in D\>$ is computed to be
$\prod_{i \in D} (1-\tt^{\aa_i})$.  Thus it suffices to show that if
$K(\tt) = 1-\tt^\bb = 1 - t_1^{b_1} \cdots t_d^{b_d}$, then substituting
$1-t_j$ for each occurrence of $t_j$ yields $K(\1-\tt) = b_1 t_1 +
\cdots + b_d t_d + O(\tt^2)$, where $O(\tt^e)$ denotes a sum of terms
each of which has total degree at least~$e$.  Indeed, then we can
conclude that
\begin{eqnarray*}
  \KK(\kk[\zz]/\!\<z_i \mid i \in D\>;\1-\tt) &=&
  \Big(\prod_{i \in D} \aa_i \Big) + O(\tt^{r+1}),
\end{eqnarray*}
where $r$ is the size of~$D$.  Calculating $K(\1-\tt)$ yields
\begin{eqnarray*}
  1 - \prod_{j=1}^d (1-t_j)^{b_j} &=& 1 - \prod_{j=1}^d \big(1 - b_j t_j
  + O(t_j^2)\big),
\end{eqnarray*}
from which we get the desired formula
\begin{eqnarray*}
  \displaystyle 1 - \Big(1 - \sum_{j=1}^d (b_j t_j) + O(\tt^2)\Big) \
  \,=\,\ \Big(\sum_{j=1}^d b_j t_j\Big) + O(\tt^2).
\end{eqnarray*}

All that remains is additivity.  Every associated prime of~$\Gamma$ is
$\ZZ^d$-graded by \cite[Exercise~3.5]{Eis}.  Choose by noetherian
induction a filtration $\Gamma = \Gamma_\ell \supset \Gamma_{\ell-1}
\supset \cdots \supset \Gamma_1 \supset \Gamma_0 = 0$ in which
$\Gamma_j/\Gamma_{j-1} \cong (\kk[\zz]/\pp_j)(-\bb_j)$ for multigraded
primes~$\pp_j$ and vectors \mbox{$\bb_j \in \ZZ^d$}.  Additivity of
\K-polynomials on short exact sequences implies that $\KK(\Gamma;\tt) =
\sum_{j=1}^\ell \KK(\Gamma_j/\Gamma_{j-1};\tt)$.

The variety of~$\pp_j$ is contained inside the support of~$\Gamma$, and
if $\pp$ has dimension exactly $\dim(\Gamma)$, then $\pp$ equals the
prime ideal of some top-dimensional component $X \in \{X_1,\ldots,X_r\}$
for exactly $\mult_X(\Gamma)$ values~of~$j$ (localize the filtration
at~$\pp$ to see this).

Assume for the moment that~$\Gamma$ is a direct sum of multigraded
shifts of quotients of~$\kk[\zz]$ by monomial ideals.
The filtration can be chosen so that all the primes~$\pp_j$ are of the
form $\<z_i \mid i \in D\>$.
By normalization and the obvious equality $\KK(\Gamma'(\bb);\tt) =
\tt^\bb\KK(\Gamma';\tt)$ for any $\ZZ^d$-graded module~$\Gamma'$, the
only power series $\KK(\Gamma_j/\Gamma_{j-1};\1-\tt)$ contributing terms
to $\KK(\Gamma;\1-\tt)$ are those for which $\Gamma_j/\Gamma_{j-1}$ has
maximal dimension.  Therefore the theorem holds for
direct sums of shifts of monomial quotients.

By Gr\"obner degenration, a general module~$\Gamma$ of codimension~$r$
has the same multidegree as a direct sum of shifts of monomial
quotients.
Using the filtration for this general~$\Gamma$, it follows from the
previous paragraph that $\KK(\Gamma_j/\Gamma_{j-1};\1-\tt) =
\cC(\Gamma_j/\Gamma_{j-1};\tt) + O(\tt^{r+1})$.  Therefore the last two
sentences of the previous paragraph work also for the general
module~$\Gamma$.%
\end{proof}

Our general view of positivity proceeds thus: Multidegrees, like
ordinary degrees, are additive on unions of schemes with equal dimension
and no common components.  Additivity under unions becomes quite useful
for monomial ideals, because their irreducible components are coordinate
subspaces, whose multidegrees are simple.  Knowing explicitly the
multidegrees of monomial subschemes of $\kk^m$ yields formulae for
multidegrees of arbitrary subschemes because multidegrees are constant
in flat families.

\begin{theorem}
\label{t:pos}
The multidegree of any module of dimension $d-r$ over a positively
$\ZZ^d$-graded ring $\kk[\zz]$ is a positive sum of
terms of the form $\aa_{i_1}\!\cdots \aa_{i_r} \in \sym^r_\ZZ(\ZZ^d)$,
where $i_1 < \cdots < i_r$.
\end{theorem}
\begin{proof}
The special fiber of any Gr\"obner degeneration of the module~$\Gamma$
has support equal to a union of coordinate subspaces.  Now use
Theorem~\ref{t:multidegs}.
\end{proof}

The products $\aa_{i_1}\!\cdots \aa_{i_r}$ are all nonzero, and all lie
in a single polyhedral cone containing no linear subspace (a semigroup
with no units) inside $\sym^r_\ZZ(\ZZ^d)$, by positivity.  Thus, when we
say ``positive sum'' in Theorem~\ref{t:pos}, we mean in particular that
the sum is nonzero.

Although the indices on $\aa_{i_1}, \ldots, \aa_{i_r}$ are distinct,
some of the weights themselves might be equal.  This occurs when $\kk^m
= \mn$ and $\weight(z_i) = x_i$, for example: any monomial of degree at
most~$n$ in each~$x_i$ is attainable.
Theorem~\ref{t:pos} implies in this case that a polynomial expressible
as the $\ZZ^n$-graded multidegree of some subscheme of~$\mn$ has
positive coefficients.  In fact, the coefficients count geometric
objects, namely subspaces (with multiplicity) in {\em any}\/ Gr\"obner
degeneration.  Therefore Theorem~\ref{t:pos} completes our second
goal~(ii) from the introduction, that of proving positivity of Schubert
polynomials in a natural geometric setting, in view of
Theorem~\ref{t:formulae}, which completed the first goal~(i).

The conditions in Theorem~\ref{t:multidegs}
overdetermine the multidegree function: there is usually no {\em
single}\/ best way to write a multidegree as a positive sum in
Theorem~\ref{t:pos}.
It happens that antidiagonal degenerations of matrix Schubert varieties
as in Theorem~\ref{t:gb} give particularly nice multiplicity~$1$
formulae, where the geometric objects have combinatorial significance as
in Theorems~\ref{t:gb} and~\ref{t:mitosis}.  The details of this story
are fleshed out in Section~\ref{sec:bjs}.%

\begin{example} \label{ex:s3deg}
Five of the six $3\times 3$ matrix Schubert varieties in
Example~\ref{ex:s3} have $\ZZ^{2n}$-graded multidegrees that are
products of expressions having the form $x_i-y_j$ by the normalization
condition in Theorem~\ref{t:multidegs}:
\begin{eqnarray*}
  [\ol X_{123}]_{\ZZ^{2n}} &=& 1
\\{}
  [\ol X_{213}]_{\ZZ^{2n}} &=& x_1 - y_1
\\{}
  [\ol X_{231}]_{\ZZ^{2n}} &=& (x_1 - y_1)(x_1-y_2)
\\{}
  [\ol X_{312}]_{\ZZ^{2n}} &=& (x_1 - y_1)(x_2-y_1)
\\{}
  [\ol X_{321}]_{\ZZ^{2n}} &=& (x_1 - y_1)(x_1-y_2)(x_2-y_1)
\end{eqnarray*}
The last one, $\ol X_{132}$, has multidegree
\begin{eqnarray*}
  [\ol X_{132}]_{\ZZ^{2n}} &=&  x_1 + x_2 - y_1 - y_2 \qquad\qquad\,
\end{eqnarray*}
that can be written as a sum of expressions $(x_i-y_j)$ in two different
ways.  To see how, pick term orders that choose different leading
monomials for $z_{11}z_{22} - z_{12}z_{21}$.  Geometrically, these
degenerate $\ol X_{132}$ to either the scheme defined by $z_{11}z_{22}$
or the scheme defined by $z_{12}z_{21}$, while preserving the
multidegree in both cases.  The degenerate limits break up as unions
\begin{eqnarray*}
  \ol X_{132}'
     &=&   \{Z \in M_3 \mid z_{11} = 0\}\cup\{Z \in M_3 \mid z_{22} = 0\}
  \ \:=\ \:\{Z \in M_3 \mid z_{11}z_{22} = 0\}
\\
  \ol X_{132}''
     &=&   \{Z \in M_3 \mid z_{12} = 0\}\cup\{Z \in M_3 \mid z_{21} = 0\}
  \ \:=\ \:\{Z \in M_3 \mid z_{12}z_{21} = 0\}
\end{eqnarray*}
and therefore have multidegrees
\begin{eqnarray*}
  [\ol X_{132}']_{\ZZ^{2n}} &=& (x_1-y_1) + (x_2-y_2)
\\{}
  [\ol X_{132}'']_{\ZZ^{2n}}&=& (x_1-y_2) + (x_2-y_1).
\end{eqnarray*}
Either way
calculates $[\ol X_{132}]_{\ZZ^{2n}}$ as in Theorem~\ref{t:pos}.  For
most permutations \mbox{$w \in S_n$}, only antidiagonal degenerations
(such as~$\ol
X_{132}\hspace{-2.5ex}\raisebox{.2ex}{$''$}\hspace{1.4ex}$) can be read
off the minors generating~$I_w$.
\end{example}

Multidegrees are functorial with respect to changes of grading, as the
following proposition says.  It holds for prime monomial quotients
$\Gamma = \kk[\zz]/\!\<z_i \mid i \in D\>$ by normalization, and generally
by Gr\"obner degeneration along with additivity.

\begin{prop} \label{prop:coarsen}
If $\ZZ^d \to \Znp$ is a homomorphism of groups, then any $\ZZ^d$-graded
module $\Gamma$ is also $\Znp$-graded.  Furthermore, $K$-polynomials and
multidegrees specialize naturally:
\begin{numbered}
\item
The $\ZZ^d$-graded $\KK$-polynomial $\KK(\Gamma,\tt)$ maps to the
$\Znp$-graded $K$-polynomial $\KK(\Gamma;\tt')$ under the natural
homomorphism $\ZZ[\ZZ^d] \to \ZZ[\Znp]$ of group rings; and

\item
The $\ZZ^d$-graded multidegree $\cC(\Gamma;\tt)$ maps to the
$\Znp$-graded multidegree $\cC(\Gamma;\tt')$ under the natural
homomorphism $\sym^\spot_\ZZ(\ZZ^d) \to \sym^\spot_\ZZ(\Znp)$.
\end{numbered}
\end{prop}

\begin{example} \label{ex:coarsen}
Changes between the gradings from Example~\ref{ex:exp} go as follows.
\def\mcolempty{\multicolumn{2}{c}{}}
$$
\begin{array}{r|c@{\ \mapsfrom\ }l@{\quad}c@{\ \mapsfrom\ \:}l@{\quad\!\!}
		r@{\ \mapsfrom\ \:}l}
\hbox{change of grading}  &
	\multicolumn{1}{c@{\ \from\ }}{\ZZ} & \ZZ^n&
	\multicolumn{1}{c@{\ \from\ \:}}{\ZZ^n} & \ZZ^{2n}&
	\multicolumn{1}{r@{\ \from\ \:}}{\ZZ^{2n}} & \Znn 
\\[0ex]\hline
\raisebox{-1.3ex}{map on variables in $K$-polynomials} &
	  t & x_i  & x_i  &   x_i   & x_i/y_j & z_{ij}
\\[-1ex]&\mcolempty&  1   &   y_j
\\[0ex]\hline
\raisebox{-1.3ex}{map on variables in multidegrees} &
	  t & x_i  & x_i  &   x_i   & x_i-y_j & z_{ij}
\\[-1ex]&\mcolempty&  0   &   y_j
\end{array}
$$
We often call these maps \bem{specialization}, or \bem{coarsening the
grading}.  Setting all occurrences of $y_j$ to zero in
Example~\ref{ex:s3deg} yields $\ZZ^n$-graded multidegrees, for instance;
compare these to the diagram in Example~\ref{ex:schub}.
\end{example}

The connective tissue in our proof of Theorems~\ref{t:formulae},
\ref{t:gb}, and~\ref{t:mitosis} (Section~\ref{sec:proof}) consists of
the next observation.  It appears in its $\ZZ$-graded form independently
in \cite{MartinPhD} (although Martin applies the ensuing conclusion that
a candidate Gr\"obner basis actually is one~to~a~different~ideal).  It
will be applied with $I' = \IN(I)$ for some ideal~$I$ and some term
order.

\begin{lemma} \label{lemma:IN}
Let $I' \subseteq \kk[z_1,\ldots,z_m]$ be an ideal homogeneous for a
positive $\ZZ^d$-grading.  Suppose that $J$ is an equidimensional
radical ideal contained inside~$I'$.  If the zero schemes of $I'$
and~$J$ have equal multidegrees, then $I' = J$.
\end{lemma}
\begin{proof}
Let $X$ and~$Y$ be the schemes defined by $I'$ and~$J$, respectively.
The multidegree of $\kk[\zz]/J$ equals the sum of the multidegrees of
the components of~$Y$, by additivity.  Since \mbox{$J \subseteq I'$},
each maximal dimensional irreducible component of~$X$ is contained in
some component of~$Y$, and hence is equal to it (and reduced) by
comparing dimensions: equal multidegrees implies equal dimensions by
Theorem~\ref{t:pos}.  Additivity says that the multidegree of~$X$ equals
the sum of multidegrees of components of~$Y$ that happen also to be
components of~$X$.  By hypothesis, the multidegrees of~$X$ and~$Y$
coincide, so the sum of multidegrees of the remaining components of~$Y$
is zero.  This implies that no components remain, by
Theorem~\ref{t:pos}, so $X \supseteq Y$.  Equivalently, $I' \subseteq
J$, whence $I' = J$ by the hypothesis $J \subseteq I'$.
\end{proof}

\section{Subword complexes in Coxeter groups}\label{sec:subword}

This section exploits the properties of reduced words in Coxeter
groups to produce shellings of the initial complex~$\LL_w$ from
Theorem~\ref{t:gb}.  More precisely, we define a new class of
simplicial complexes that generalizes to arbitrary Coxeter groups the
construction in Section~\ref{sec:pipe} of reduced pipe dreams for a
permutation~$w \in S_n$ from the triangular reduced expression
for~$w_0$.  The manner in which subword complexes characterize reduced
pipe dreams is similar in spirit to
\cite{FKyangBax}; however, even for reduced pipe dreams our
topological perspective is~new.

We felt it important to include the Cohen--Macaulayness of the initial
scheme~$\LL_w$ as part of our evidence for the naturality of Gr\"obner
geometry for Schubert polynomials, and the generality of subword
complexes allows our simple proof of their shellability.  However, a
more detailed analysis would take us too far afield, so we have chosen
to develop the theory of subword complexes in Coxeter groups more fully
elsewhere~\cite{KMsubword}.  There, we show that subword complexes are
balls or spheres, and calculate their Hilbert series for applications to
Grothendieck polynomials.  We also comment there on how our forthcoming
Theorem~\ref{t:subword} reflects topologically some of the fundamental
properties of reduced (and nonreduced) expressions in Coxeter groups,
and how Theorem~\ref{t:subword} relates to known results on simplicial
complexes constructed from Bruhat and weak orders.

Let $(\Pi,\Sigma)$ be a Coxeter system, so $\Pi$ is a Coxeter group and
$\Sigma$ is a set of simple reflections, which generate~$\Pi$.  See
\cite{HumphCoxGrps} for background and definitions; the applications to
reduced pipe dreams concern only the case where $\Pi = S_n$ and $\Sigma$
consists of the adjacent transpositions switching $i$ and $i+1$ for $1
\leq i \leq n-1$.

\begin{defn} \label{defn:subword}
A \bem{word} of size~$m$ is an ordered sequence $Q = (\sigma_1, \ldots,
\sigma_m)$ of elements of~$\Sigma$.  An ordered subsequence $P$ of $Q$
is called a \bem{subword} of $Q$.
\begin{numbered}
\item
$P$ \bem{represents} $\pi \in \Pi$ if the ordered product of the simple
reflections in $P$ is a reduced decomposition for $\pi$.
\item
$P$ \bem{contains} $\pi \in \Pi$ if some subsequence of $P$
represents~$\pi$.
\end{numbered}
The \bem{subword complex} $\Delta(Q,\pi)$ is the set of subwords $P
\subseteq Q$ whose complements $Q \minus P$ contain $\pi$.
\end{defn}

Often we write~$Q$ as a string without parentheses or commas, and abuse
notation by saying that $Q$ is a word in~$\Pi$.
Note that $Q$ need not itself be a reduced expression, but the facets of
$\Delta(Q,\pi)$ are the complements of reduced subwords of~$Q$.  The
word~$P$
contains~$\pi$ if and only if the product of~$P$ in the
degenerate Hecke algebra
is~$\geq\!\pi$ in Bruhat order \cite{FKyangBax}.

\begin{example} \label{ex:pentagon}
Let $\Pi = S_4$, and consider the subword complex $\Delta =
\Delta(s_3s_2s_3s_2s_3,1432)$.  Then $\pi = 1432$ has two reduced
expressions, namely $s_3s_2s_3$ and $s_2s_3s_2$.  Labeling the vertices
of a pentagon with the reflections in~$Q = s_3s_2s_3s_2s_3$ (in cyclic
order), we find that the facets of $\Delta$ are the pairs of adjacent
vertices.  Therefore $\Delta$ is the pentagonal boundary.%
\end{example}

\begin{example} \label{ex:complex}
Let $\Pi = S_{2n}$ and let the \bem{square word}
\begin{eqnarray*}
  Q_{n \times n} &=& s_ns_{n-1} \ldots s_2s_1\ s_{n+1}s_n \ldots s_3s_2
  \ \ \ldots\ \ s_{2n-1}s_{2n-2} \ldots s_{n+1}s_n
\end{eqnarray*}
be the ordered list constructed from the pipe dream whose crosses
entirely fill the \mbox{$n \times n$} grid.  Reduced expressions for
permutations $w \in S_n$ never involve reflections~$s_i$ with \mbox{$i
\geq n$}.  Therefore, if $Q_0$~is the triangular long word for $S_n$
(not $S_{2n}$) in Example~\ref{ex:subword}, then $\Delta(Q_{n \times
n},w)$ is the join of $\Delta(Q_0,w)$ with a simplex whose $\binom n2$
vertices correspond to the lower-right triangle of the $n \times n$
grid.  Consequently, the facets of $\Delta(Q_{n \times n},w)$ are
precisely the complements in $[n] \times [n]$ of the reduced pipe dreams
for~$w$, by Lemma~\ref{lemma:subword}.%
\end{example}

The following lemma is immediate from the definitions and the fact that
all reduced expressions for $\pi \in \Pi$ have the same length.

\begin{lemma} \label{lemma:represents}
$\Delta(Q,\pi)$ is a pure simplicial complex whose facets are the
subwords $Q \minus P$ such that $P \subseteq Q$ represents
$\pi$.\hfill$\Box$
\end{lemma}

\begin{defn} \label{defn:vertex}
Let $\Delta$ be a simplicial complex and $F \in \Delta$ a face.
\begin{numbered}
\item
The \bem{deletion} of $F$ from $\Delta$ is $\del(F,\Delta) = \{G \in
\Delta \mid G \cap F = \nothing\}$.
\item
The \bem{link} of $F$ in $\Delta$ is $\link(F,\Delta) = \{G \in \Delta
\mid G \cap F = \nothing$ and $G \cup F \in \Delta\}$.
\end{numbered}
$\Delta$ is \bem{vertex-decomposable} if $\Delta$ is pure and either
(1)~$\Delta = \{\nothing\}$, or (2)~for some vertex $v \in \Delta$, both
$\del(v,\Delta)$ and $\link(v,\Delta)$ are vertex-decomposable.  A
\bem{shelling} of $\Delta$ is an ordered list $F_1, F_2, \ldots, F_t$ of
its facets such that $\bigcup_{j < i} F_j \cap F_i$ is a union of
codimension~1 faces of $F_i$ for each $i \leq t$.  We say $\Delta$ is
\bem{shellable} if it is pure and has a shelling.
\end{defn}
%

Provan and Billera \cite{BilleraProvan} introduced the notion of
vertex-decomposability and proved that it implies shellability (proof:
use induction on the number of vertices by first shelling
$\del(v,\Delta)$ and then shelling the cone from $v$ over
$\link(v,\Delta)$ to get a shelling of~$\Delta$).  It is well-known that
shellability implies Cohen--Macaulayness \cite[Theorem~5.1.13]{BH}.
Here, then, is our central observation concerning subword complexes.

\begin{theorem} \label{t:subword}
Any subword complex $\Delta(Q,\pi)$ is vertex-decomposable.  In
particular, subword complexes are shellable and therefore
Cohen--Macaulay.
\end{theorem}


\begin{proof}
With $Q = (\sigma, \sigma_2, \sigma_3, \ldots, \sigma_m)$, we show
that both the link and the deletion of $\sigma$ from $\Delta(Q,\pi)$
are subword complexes.  By definition, both consist of subwords of $Q'
= (\sigma_2, \ldots , \sigma_m)$.  The link is naturally identified
with the subword complex $\Delta(Q',\pi)$.  For the deletion, there
are two cases.  If $\sigma \pi$ is longer than~$\pi$, then the
deletion of $\sigma$ equals its link because no reduced expression for
$\pi$ begins with~$\sigma$.  On the other hand, when $\sigma \pi$ is
shorter than~$\pi$, the deletion is $\Delta(Q',\sigma\pi)$.%
\end{proof}

\begin{remark} \label{rk:subword}
The vertex decomposition that results for initial ideals of matrix
Schubert varieties has direct analogues in the Gr\"obner degenerations
and formulae for Schubert polynomials.  Consider the sequence
$>_1,>_2,\ldots,>_{n^2}$ of partial term orders, where $>_i$ is
lexicographic in the first $i$ matrix entries snaking from northeast to
southwest one row at a time, and treats all remaining variables equally.
The order $>_{n^2}$ is a total order; this total order is antidiagonal,
and hence degenerates $\ol X_w$ to the subword complex by
Theorem~\ref{t:gb} and Example~\ref{ex:complex}.  Each $>_i$ gives a
degeneration of $\ol X_w$ to a union of components, every one of which
degenerates at $>_{n^2}$ to its own subword complex.

If we study how a component at stage $i$ degenerates into components at
stage $i+1$, by degenerating both using $>_{n^2}$, we recover the vertex
decomposition for the corresponding subword complex.

Note that these components are {\em not}\/ always matrix Schubert
varieties; the set of rank conditions involved does not necessarily
involve only upper-left submatrices. We do not know how general a class
of determinantal ideals can be tackled by partial degeneration of matrix
Schubert varieties, using antidiagonal partial term orders.

However, if we degenerate using the partial order $>_n$ (order just the
first row of variables), then the components {\em are}\/ matrix Schubert
varieties, except that the minors involved are all shifted down one
row. This gives a geometric interpretation of the inductive formula for
Schubert polynomials appearing in Section~1.3 of~\cite{BJS}.
\end{remark}

\part{Applications of the Gr\"obner geometry theorems}\label{part:app}

\section{Positive formulae for Schubert polynomials}\label{sec:bjs}

The original definition of Schubert polynomials by Lascoux and
Sch\"utzenberger via the divided difference recursion involves negation,
so it is quite nonobvious from their formulation that the coefficients of
$\SS_w(\xx)$ are in fact positive.  Although Lascoux and Sch\"utzenberger
did prove positivity using their `transition formula', the first
combinatorial proofs, showing what the coefficients count, appeared in
\cite{BJS, FSnilCoxeter}.  More recently,
\cite{KirMaeno00,BerSot02,KoganThesis} show that the coefficients are
positive for geometric reasons.

Our approach has the advantage that it produces geometrically a uniquely
determined polynomial representative for each Schubert class, and
moreover, that it provides an obvious geometric reason why this
representative has nonnegative coefficients in the variables
$x_1,\ldots,x_n$ (or $\{x_i-y_j\}_{i,j=1}^n$ in the double case).  Only
then do we identify the coefficients as counting known combinatorial
objects; it is coincidence (or naturality of the combinatorics) that our
positive formula for Schubert polynomials agrees with---and provides a
new geometric proof of---the combinatorial formula of Billey, Jockusch,
and Stanley \cite{BJS}.

\begin{thm} \label{thm:positive}
There is a multidegree formula that writes
\begin{eqnarray*}
  \SS_w\:\ =\:\ [\ol X_w] &=& \sum_{L \in \LL_w} [L]
\end{eqnarray*}
as a sum over the facets $L$ of the initial complex $\LL_w$, thereby
expressing the Schubert polynomial $\SS_w(\xx)$ as a positive of
monomials $[L]_{\ZZ^n} = \prod_{(i,j) \in D_L} x_i$ in the variables
$x_1,\ldots,x_n$, and the double Schubert polynomial $\SS_w(\xx,\yy)$ as
a sum of expressions $[L]_{\ZZ^{2n}} = \prod_{(i,j) \in D_L} (x_i -
y_j)$, which are themselves positive in the variables $x_1,\ldots,x_n$
and $-y_1,\ldots,-y_n$.
\end{thm}
\begin{proof}
By Theorem~\ref{t:gb} and degeneration in Theorem~\ref{t:multidegs}, the
multidegrees of~$\ol X_w$ and the zero set $\LL_w$ of~$J_w$ are equal in
any grading.  Since $[\ol X_w] = \SS_w$ by Theorem~\ref{t:formulae}, the
formulae then follow from additivity and normalization in
Theorem~\ref{t:multidegs}, given the ordinary weights in
Example~\ref{ex:coarsen}.%
\end{proof}

\begin{remark} \label{rk:effective}
The version of this positivity in algebraic geometry is the notion of
``effective homology class'', meaning ``representable by a
subscheme''.  On the flag manifold, a homology class is effective
exactly if it is a nonnegative combination of Schubert classes.
(Proof: one direction is a tautology. For the other, if $X$ is a
subscheme of the flag manifold~$\FL$, consider the induced action of
the Borel group $B$ on the Hilbert scheme for $\FL$.  The closure of
the $B$-orbit through the Hilbert point~$X$ will be projective because
the Hilbert scheme is, so Borel's theorem produces a fixed point,
necessarily a union of Schubert varieties, perhaps nonreduced.)  In
particular the classes of monomials in the $x_i$ (the first Chern
classes of the standard line bundles; see Section~\ref{sec:flag}) are
not usually effective.

We work instead on $\mn$, where the standard line bundles become trivial,
but not equivariantly, and a class is effective exactly if it is a
nonnegative combination of monomials
in the equivariant first Chern classes~$x_i$. (Proof: instead of using
$B$ to degenerate a subscheme~$X$ inside~$\mn$, use a $1$-parameter
subgroup of the $n^2$-dimensional torus.  Algebraically, this amounts to
picking a Gr\"obner basis.)
\end{remark}

The next formula was our motivation for relating the antidiagonal
complex $\LL_w$ to the set~$\rc(w)$ of reduced pipe dreams.  Although
formulated here in language based on~\cite{FKyangBax}, its first proof
(in transparently equivalent language) was by Billey, Jockusch, and
Stanley, while Fomin and Stanley shortly thereafter gave a better
combinatorial proof.

\begin{cor}[\cite{BJS,FSnilCoxeter}] \label{cor:BJS}
$\displaystyle \SS_w(\xx)\ = \sum_{D \in \rc(w)} \xx^D$, where
$\displaystyle \xx^D = \prod_{(i,j) \in D} x_i$.
\end{cor}
\begin{proof}
Apply Theorem~\ref{t:gb} to the $\ZZ^n$-graded version of the formula in
Theorem~\ref{thm:positive}.
\end{proof}

\begin{example}
As in the example from the Introduction,
Lemma~\ref{lemma:subword} calculates the multidegree as
$$
\begin{array}{cccccccccc}
  [\ol X_{2143}]
   &=& [L_{11,13}] &+& [L_{11,22}] &+& [L_{11,31}]
\\ &=& \weight(z_{11}z_{13})&+&\weight(z_{11}z_{22})&+&\weight(z_{11}z_{31})
\\ &=& x_1^2&+&x_1x_2&+&x_1x_3
\end{array}
$$
in $\ZZ[x_1,x_2,x_3,x_4]$, for the $\ZZ^n$-grading.%
\end{example}

The double version of Corollary~\ref{cor:BJS} has the same proof, using
the $\ZZ^{2n}$-grading; deriving it directly from reduced pipe dreams
here bypasses the ``double rc-graphs'' of~\cite{BB}.

\begin{cor}[\cite{FKyangBax}] \label{cor:doubSchub}
$\displaystyle \SS_w(\xx,\yy)\ = \sum_{D \in \rc(w)} \prod_{(i,j) \in D}
(x_i - y_j)$.\qed
\end{cor}

Theorem~\ref{thm:positive} and Corollary~\ref{cor:BJS} together are
consequences of Theorem~\ref{t:gb} and the multidegree part of
Theorem~\ref{t:formulae}.  There is a more subtle kind of positivity for
Grothendieck polynomials, due to Fomin and Kirillov
\cite{FKgrothYangBax}, that can be derived from Theorem~\ref{t:gb} and
the Hilbert series part of Theorem~\ref{t:formulae}, along with the
Eagon--Reiner theorem from combinatorial commutative algebra~\cite{ER}.
In fact, this ``positivity'' was our chief evidence leading us to
conjecture the Cohen--Macaulayness of~$J_w$.

More precisely, the work of Fomin and Kirillov implies that for each
$d$ there is a homogeneous polynomial $\GG_w^{(d)}(\xx)$ of degree~$d$
with nonnegative coefficients such that
\begin{eqnarray} \label{eq:alternates}
  \GG_w(\1-\xx) &=& \sum_{d \geq \ell} (-1)^{d-\ell}\GG_w^{(d)}(\xx),
\end{eqnarray}
where $\ell = \length(w)$.  In other words, the coefficients on each
homogeneous piece of $\GG_w(\1-\xx)$ all have the same sign.  On the
other hand, the Eagon--Reiner theorem states:
\begin{quote}
A simplicial complex $\Delta$ is Cohen--Macaulay if and only if the
Alexander dual $J_\Delta^\star$ of its Stanley--Reisner ideal has
\bem{linear free resolution}, meaning that the differential in its
minimal $\ZZ$-graded free resolution over $\kk[\zz]$ can be expressed
using matrices filled with linear~forms.
\end{quote}
The $K$-polynomial of any module with linear resolution alternates as
in~(\ref{eq:alternates}).  But the Alexander inversion formula
\cite{KMsubword} implies that $\GG_w(\1-\xx)$ is the $K$-polynomial of
$J_w^\star$, given that $\GG_w(\xx)$ is the $K$-polynomial of $\cj w$ as
in Theorem~\ref{t:formulae}.  Therefore, $\GG_w(\1-\xx)$ must alternate
as in~(\ref{eq:alternates}), {\em if the Cohen--Macaulayness in
Theorem~\ref{t:gb} holds}.  It would take suspiciously fortuitous
cancelation to have a squarefree monomial ideal $J_w^\star$ whose
$K$-polynomial $\GG_w(\1-\xx)$ behaves like~(\ref{eq:alternates})
without the ideal $J_w^\star$ actually having linear resolution.

In fact, further investigation into the algebraic combinatorics of
subword complexes can identify the coefficients of the homogeneous
pieces of (double) Grothendieck polynomials.  We carry out this program
in \cite{KMsubword}, recovering a formula of Fomin and Kirillov~%
\cite{FKgrothYangBax}.

\section{Degeneracy loci}\label{sec:loci}

We recall here Fulton's theory of degeneracy loci, and explain its
relation to equivariant cohomology.  This was our initial interest in
Gr\"obner geometry of double Schubert polynomials: to get universal
formulae for the cohomology classes of degeneracy loci.  However, since
completing this work, we learned of the papers \cite{FRthomPoly,Kaz97}
taking essentially the same viewpoint, and we refer to them for detail.

Given a flagged vector bundle $E_\spot = (E_1 \into E_2 \into \cdots
\into E_n)$ and a co-flagged vector bundle $F_\spot = (F_n \onto F_{n-1}
\onto \cdots \onto F_1)$ over the same base $X$, a generic map $\sigma :
E_n \to F_n$, and a permutation $w$, define the \bem{degeneracy locus}
$\WW_w$ as the subset
\begin{eqnarray*}
  \WW_w &=& \{x\in X \mid \rank (E_q \to E_n \stackrel \sigma \too F_n
  \onto F_p) \leq \rank(w^T\sub qp) \hbox{ for all } q,p\}.
\end{eqnarray*}
The principal goal in Fulton's paper \cite{FulDegLoc} was to provide
``formulae for degeneracy loci'' as polynomials in the Chern classes of
the vector bundles.  In terms of the Chern roots $\{c_1(E_p/E_{p-1}),
c_1(\ker F_q \to F_{q-1}) \}$, Fulton found that the desired polynomials
were actually the double Schubert polynomials.

It is initially surprising that there is a single formula, for all
$X,E,F$ and not really depending on $\sigma$.  This follows from a
classifying space argument, when $\kk = \CC$, as follows.

The group of automorphisms of a flagged vector space consists of the
invertible lower triangular matrices~$B$, so the classifying space $BB$
of $B$-bundles carries a universal flagged vector bundle. The
classifying space of interest to us is thus $BB \times BB_+$, which
carries a pair of universal vector bundles $\mathcal E$ and $\mathcal
F$, the first flagged and the second co-flagged. We write
$\hhom({\mathcal E},{\mathcal F})$ for the bundle whose fiber at
$(x,y)\in BB\times BB_+$ equals $\Hom({\mathcal E}_x,{\mathcal F}_y)$.

Define the \bem{universal degeneracy locus} $U_w \subseteq
\hhom({\mathcal E},{\mathcal F})$ as the subset
\begin{eqnarray*}
  U_w &=& \{(x,y,\phi) \mid \rank (({\mathcal E}_x)_q \stackrel \phi \too
  ({\mathcal F}_y)_p) \leq \rank(w^T\sub qp) \hbox{ for all } q,p\},
\end{eqnarray*}
where \hbox{$x \in BB$}, \hbox{$y \in BB_+$}, and $\phi : \EE_x \to
\FF_y$.  In other words, the homomorphisms in the fiber of $U_w$ at
$(x,y)$ lie in the corresponding matrix Schubert variety.

The name is justified by the following. Recall that our setup is a space
$X$, a flagged vector bundle $E$ on it, a coflagged vector bundle $F$,
and a `generic' vector bundle map $\sigma : E \to F$; we will soon see
what `generic' means.  Pick a classifying map $\chi : X \to BB \times
BB_+$, which means that $E,F$ are isomorphic to pullbacks of the
universal bundles.  (Classifying maps exist uniquely up to homotopy.)
Over the target we have the universal $\Hom$-bundle $\hhom({\mathcal
E},{\mathcal F})$, and the vector bundle map $\sigma$ is a choice of a
way to factor the map $\chi$ through a map $\tilde\sigma : X \to
\hhom({\mathcal E},{\mathcal F})$.  The degeneracy locus $\WW_w$ is then
$\tilde\sigma^{-1}(U_w)$, and it is natural to request that
$\tilde\sigma$ be transverse to each $U_w$---this will be the notion of
$\sigma$ being generic.

What does this say cohomologically?  The closed subset $U_w$ defines a
class in Borel--Moore homology (and thus ordinary cohomology) of the
bundle $\hhom({\mathcal E},{\mathcal F})$.  If $\sigma$ is generic, then
\begin{eqnarray*}
  [\WW_w] &=& [ \tilde\sigma^{-1}(U_w) ] \ \:=\ \:\tilde\sigma^*([U_w])
\end{eqnarray*}
and this is the sense in which there is a universal formula $[U_w] \in
H^*(\hhom({\mathcal E},{\mathcal F}))$.  The cohomology ring of this
$\Hom$-bundle is the same as that of the base $BB \times BB_+$ (to which
it retracts), namely a polynomial ring in the $2n$ first Chern classes,
so one knows a~priori that the universal formula should be expressible as
a polynomial in these $2n$ variables.

We can rephrase this using Borel's mixing space definition of equivariant
cohomology.  Given a space $S$ carrying an action of a group $G$, and a
contractible space $EG$ upon which $G$ acts freely, the equivariant
cohomology $H^*_G(S)$ of $S$ is defined as
\begin{eqnarray*}
  H^*_G(S) &:=& H^*( (S \times EG)/G),
\end{eqnarray*}
where the quotient is respect to the diagonal action. Note that the Borel
`mixing space' \mbox{$(S\times EG)/G$} is a bundle over $EG/G =: BG$,
with fibers $S$. In particular $H^*_G(S)$ is automatically a module over
$H^*(BG)$, thereby called the `base ring' of $G$-equivariant cohomology.

For us, the relevant group is $B\times B_+$, and we have two spaces $S$:
the space of matrices $\mn$ under left and right multiplication, and
inside it the matrix Schubert variety $\ol X_w$.  Applying the mixing
construction to the pair $\mn \supseteq \ol X_w$, it can be shown that
we recover the bundles $\hhom({\mathcal E},{\mathcal F}) \supseteq U_w$.
As such, the universal formula $[U_w] \in H^*(\hhom({\mathcal
E},{\mathcal F}))$ we seek can be viewed instead as the class defined in
$(B\times B_+)$-equivariant cohomology by $\ol X_w$ inside $\mn$.  As we
prove in Theorem~\ref{t:formulae} (in the setting of multidegrees,
although a direct equivariant cohomological version is possible), these
are the double Schubert polynomials.

The main difference between this mixing space approach and that of
Fulton in \cite{FulDegLoc} is that in the algebraic category, where
Fulton worked, some pairs $(E,F)$ of algebraic vector bundles may have
no {\em algebraic}\/ generic maps~$\sigma$.  The derivation given above
works more generally in the topological category, where no restriction
on $(E,F)$ is necessary.

In addition, we don't even need to know a~priori which polynomials
represent the cohomology classes of matrix Schubert varieties to show
that these classes are the universal degeneracy locus classes.  This
contrasts with methods relying on divided differences.

\section{Schubert classes in flag manifolds}\label{sec:flag}

Having in the main body of the exposition supplanted the topology of the
flag manifold with multigraded commutative algebra, we would like now to
connect back to the topological language.  In particular, we recover a
geometric result from our algebraic treatment of matrix Schubert
varieties: the (double) Grothendieck polynomials represent the
($B_+$-equivariant) $K$-classes of ordinary Schubert varieties in the
flag manifold \cite{LSgrothVarDrap,LascouxGrothFest}.

Our derivation of this result requires no prerequisites concerning the
rationality of the singularities of Schubert varieties: the multidegree
proof of the Hilbert series calculation is based on cohomological
considerations that ignore phenomena at complex codimension~$1$ or more,
and automatically produces the $K$-classes as numerators of Hilbert
series.  The material in this section actually formed the basis for our
original proof of Theorem~\ref{t:gb} over $\kk = \CC$, and therefore of
Theorem~\ref{thm:positive}, before we had available the technology of
multidegrees.

We use standard facts about the flag variety and notions from
(equivariant) algebraic $K$-theory, for which background material can be
found in
\cite{FulIT}.
In particular, we use freely the correspondence between $T$-equivariant
sheaves on~$\mn$ and $\ZZ^n$-graded $\kk[\zz]$-modules, where $T$ is the
torus of diagonal matrices acting by left multiplication on the left.
Under this correspondence, the $K$-polynomial $\KK(\Gamma;\xx)$ equals
the $T$-equivariant $K$-class $[\tilde\Gamma]_T \in K^\circ_T(\mn) \cong
\ZZ[\xx^{\pm1}]$ of the associated sheaf $\tilde\Gamma$ on~$\mn$.

The $K$-cohomology ring $K^\circ(\FL)$ is the quotient of $\ZZ[\xx]$ by
the ideal
\begin{eqnarray*}
  K_n &=& \<e_d(\xx)-{\textstyle{\binom nd}}\mid d\leq n\>,
\end{eqnarray*}
where $e_d$ is the $d^\th$ elementary symmetric function.  These
relations hold in $K^\circ(\FL)$ because the exterior power $\bigwedge^d
\kk^n$ of the trivial rank $n$ bundle is itself trivial of rank $\binom
nd$, and there can be no more relations because $\ZZ[\xx]/K_n$ is an
abelian group of rank $n!\/$.  Indeed, substituting $\tilde x_k = 1 -
x_k$, we find that $\ZZ[\xx]/K_n \cong \ZZ[\tilde\xx]/\tilde K_n$, where
$\tilde K_n = \<e_d(\tilde\xx) \mid d \leq n\>$, and this quotient has
rank $n!$ because it is isomorphic to the familiar cohomology ring
of~$\FL$~\cite{Borel53}.

Thus it makes sense to say that a polynomial in $\ZZ[\xx]$ ``represents
a class'' in $K^\circ(\FL)$.  Lascoux and Sch\"utzenberger, based on
work of Bernstein-Gel$'$fand-Gel$'$fand \cite{BGG} and Demazure
\cite{Dem}, realized that the classes $[\OO_{X_w}] \in \ZZ[\xx]/K_n$ of
(structure sheaves of) Schubert varieties could be represented
independently of~$n$.  To make a precise statement, let $\FLN = \bglN$
be the manifold of flags in $\kk^N$ for $N \geq n$, so $B$ is understood
to consist of $N \times N$ lower triangular matrices.  Let $X_w(N)
\subseteq \FLN$ be the \bem{Schubert variety} for the permutation $w \in
S_n$ considered as an element of $S_N$ that fixes $n+1, \ldots, N$.  In
our conventions, $X_w = B\dom\hspace{-.2ex}(\gln \cap \ol X_w)$, and
similarly for $N \geq n$.

\begin{cor}[\cite{LSgrothVarDrap,LascouxGrothFest}] \label{cor:groth}
The Grothendieck polynomial $\GG_w(\xx)$ represents the $K$-class
$[\OO_{X_w(N)}] \in K^\circ(\FLN)$ for all $N \geq n$.
\end{cor}

This is almost a direct consequence of Theorem~\ref{t:formulae}, but we
do still need a lemma.  Note that $\GG_w(\xx)$ is expressed without
reference to~$N$; here is the reason why.

\begin{lemma} \label{lemma:stable}
The $n$-variable Grothendieck polynomial $\GG_w(\xx)$ equals the
Grothen\-dieck polynomial $\GG_{w_N}(x_1,\ldots,x_N)$, whenever $w_N$
agrees with $w$ on $1,\ldots,n$ and fixes $n+1,\ldots,N$.
\end{lemma}
\begin{proof}
The ideal $I_{w_N}$ in the polynomial ring $\kk[z_{ij} \mid
i,j=1,\ldots,N]$ is extended from the ideal $I_w$ in the multigraded
polynomial subring~$\kk[\zz] = \kk[z_{ij} \mid i,j=1,\ldots,n]$.
Therefore $I_{w_N}$ has the same multigraded Betti numbers as $I_w$, so
their $K$-polynomials, which equal the Grothendieck polynomials $\GG_w$
and $\GG_{w_N}$ by Theorem~\ref{t:formulae}, are equal.
\end{proof}

\begin{proofof}{Corollary~\ref{cor:groth}}
In view of Lemma~\ref{lemma:stable}, we may as well assume $N = n$.  Let
us justify the following diagram:
$$
\begin{array}{ccccc}
	X_w
& &
& &	\ol X_w
\\
	\cap
& &
& &	\cap
\\
\bgln
&\twoheadleftarrow&\gln
&\hookrightarrow&\mn	
\\[5pt]
K^\circ(\bgln)
&\congto&K^\circ_B(\gln)
&\twoheadleftarrow& K^\circ_B(\mn)
\end{array}
$$
Pulling back vector bundles under the quotient map $\bgln
\twoheadleftarrow \gln$ induces the isomorphism $K^\circ(\bgln) \to
K^\circ_B(\gln)$.  The inclusion $\gln \hookrightarrow \mn$ induces a
surjection $K^\circ_B(\gln) \twoheadleftarrow K^\circ_B(\mn)$ because
the classes of (structure sheaves of) algebraic cycles generate both of
the equivariant $K$-homology groups $K_\circ^B(\mn)$ and
$K_\circ^B(\gln)$.

Now let $\wt X_w = \ol X_w \cap \gln$.  Any $B$-equivariant resolution
of $\OO_{\ol X_w} = \ci w$ by vector bundles on $\mn$ pulls back to a
$B$-equivariant resolution $E_\spot$ of $\OO_{\wt X_w}$ on $\gln$.
Viewing a vector bundle on $\gln$ as a geometric object (i.e.\ as the
scheme $E = \SPEC(\sym^\spot\EE^\vee)$ rather than its sheaf of sections
$\EE = \Gamma(E)$), the quotient $B \dom E_\spot$ is a resolution of
$\OO_{X_w}$ by vector bundles on $\bgln$.  Thus $[\OO_{\ol X_w}]_B \in
K^\circ_B(\mn)$ maps to $[\OO_{X_w}] \in K^\circ(\bgln)$.

The corollary follows by identifying the $B$-equivariant class
$[\OO_{\ol X_w}]_B$ as the $T$-equivariant class $[\OO_{\ol X_w}]_T$
under the natural isomorphism $K^\circ_B(\mn) \to K^\circ_T(\mn)$, and
identifying the $T$-equivariant class as the $K$-polynomial $\KK(\ci
w;\xx) = \GG_w(\xx)$ by Theorem~\ref{t:formulae}.
\end{proofof}

\begin{remark} \label{rk:reverse}
The same line of reasoning recovers the double version of
Corollary~\ref{cor:groth}, in which $\GG_w(\xx)$ and $K^\circ(\FLN)$ are
replaced by $\GG_w(\xx,\yy)$ and the equivariant $K$-group
$K^\circ_{B^+}(\FLN)$ for the action of the invertible upper triangular
matrices $B^+$ by inverses on the right.%
\end{remark}

The above proof can be worked in reverse: by assuming
Corollary~\ref{cor:groth} one can then {\em conclude}\/
Theorem~\ref{t:formulae}.  This was in fact the basis for our first
proof of Theorem~\ref{t:formulae}.  However, it requires substantially
more prerequisites (such as rationality of singularities for Schubert
varieties), and is no shorter because it fails to eliminate the
inductive arguments in Part~\ref{part:gb}.

\begin{remark}
There exists technology to assign equivariant cohomology classes to
complex subvarieties of noncompact spaces such as~$\mn$ in the cases that
interest us (see \cite{Kaz97,FRthomPoly}, for instance).  Therefore the
argument for Corollary~\ref{cor:groth} also works when $\GG_w(\xx)$ is
replaced by a Schubert or double Schubert polynomial, and $K^\circ(\FLN)$
is replaced by the appropriate version of cohomology, either $H^*(\FLN)$
or $H^*_{B^+}(\FLN)$.

Just as in Remark~\ref{rk:reverse}, this argument can be reversed: by
assuming the results of \cite{LSpolySchub} that characterize Schubert
polynomials in terms of stability properties, one can then conclude the
multidegree statement in Theorem~\ref{t:formulae}.  Since this part of
Theorem~\ref{t:formulae} is essential to proving the other main theorems
from Part~\ref{part:intro}, what we actually do is give an independent
proof of the multidegree part of Theorem~\ref{t:formulae} in
Section~\ref{sec:multimat}, to avoid issues of direct translation
between equivariant cohomology and multidegrees.
\end{remark}


\begin{remark} \label{rk:quirk}
The substitution $\xx \mapsto \1 - \xx$ in the definition of multidegree
(Section~\ref{sec:mult}) is the change of basis accompanying the
Poincar\'e isomorphism from $K$-cohomology to $K$-homology.  In general
geometric terms, $c_1(L_i) \in H^*(\FL)$ is the cohomology class
Poincar\'e dual to the divisor~$D_i$ of the $i^\th$ standard line
bundle~$L_i$ on~$\FL$.  The exact sequence $0 \to L_i^\vee \to \OO \to
\OO_{D_i} \to 0$ implies that the $K$-homology class $[\OO_{D_i}]$
equals the $K$-cohomology class $1 - [L_i^\vee]$.  Thus $\GG_w(\xx)$
writes $[\OO_{X_w}]$ as a polynomial in the Chern characters $x_i =
e^{c_1(L_i)}$ of the line bundles~$L_i$, whereas $\GG_w(\1-\xx)$ writes
$[\OO_{X_w}]$ as polynomial in the expressions
\begin{eqnarray*}
  1-e^{c_1(L_i)} &=& 1-e^{-c_1(L_i^\vee)}\ \:=\ \:c_1(L_i^\vee) -
  \frac{c_1(L_i^\vee)^2}{2!} + \frac{c_1(L_i^\vee)^3}{3!}  -
  \frac{c_1(L_i^\vee)^4}{4!}  + \cdots,
\end{eqnarray*}
whose lowest degree terms are the first Chern classes $c_1(L_i^\vee)$ of
the dual bundles $L_i^\vee$.  Forgetting the higher degree terms here,
in Definition~\ref{defn:multideg}, and in Lemma~\ref{lemma:schubert}
amounts to taking images in the associated graded ring of
$K_\circ(\fln)$, which is $H^*(\fln)$.  See \cite[Chapter 15]{FulIT} for
details.

It is an often annoying quirk of history that we end up using the same
variable $x_i$ for both $c_1(L_i^\vee) \in H^*$ and $[L_i] \in K^\circ$.
We tolerate (and sometimes even come to appreciate) this confusing abuse
of notation because it can be helpful at times.  In terms of algebra, it
reinterprets the displayed equation as: the lowest degree term in
$1-e^{-x_i}$ is just $x_i$ again.%
\end{remark}

\begin{remark}
Not only do the cohomological and \K-theoretic statements in
Theorem~\ref{t:formulae} descend to the flag manifold~$\FL$, but so also
does the degeneration of Theorem~\ref{t:gb} \cite{KoganM}.  On~$\FL$,
the degeneration can be interpreted in representation theory, where it
explains geometrically the construction of Gel$'$fand--Cetlin bases for
$\gln$ representations \mbox{\cite{GC50,GS83}}.
\end{remark}

\section{Ladder determinantal ideals}\label{sec:ladder}

The importance of Gr\"obner bases in recent work on determinantal ideals
and their relatives, such as their powers and symbolic powers, cannot be
overstated.  They are used in treatments of questions about
Cohen--Macaulayness, rational singularities, multiplicity, dimension,
$a$-invariants, and divisor class groups; see \cite{CGG, StuGBdetRings,
HerzogTrung, ConcaLadDet, MotSoh1, ConcaHerLadRatSing,
BruConKRSandPowers, KP99, GoncMilMixLadDet} for a small sample.  Since
determinantal ideals and their Gr\"obner bases also arise in the study
of (partial) flag varieties and their Schubert varieties (see
\cite{MulayLadders, GonLakSchubToricLadDet, BilLak, GonLakLadDetSchub,
GoncMilMixLadDet}, for instance), it is surprising to us that Gr\"obner
bases for the determinantal ideals defining matrix Schubert varieties
$\ol X_w$ do not seem to be in the literature, even though the ideals
themselves appeared~in~\cite{FulDegLoc}.

Most of the papers above concern a class of determinantal ideals
called `(one-sided) ladder determinantal ideals', about which we now
comment.  Consider a sequence of boxes $(b_1,a_1),\ldots,(b_k,a_k)$ in
the $n \times n$ grid, with
\begin{eqnarray*}
  a_1 \leq a_2 \leq \cdots \leq a_k \quad\hbox{and}\quad
  b_1 \geq b_2 \geq \cdots \geq b_k.
\end{eqnarray*}
Fill the boxes $(b_\ell,a_\ell)$ with nonnegative integers $r_\ell$
satisfying
\begin{equation} \label{eq:ranks}
  0 < a_1 - r_1 < a_2 - r_2 < \cdots < a_k - r_k \quad\hbox{and}\quad
  b_1 - r_1 > b_2 - r_2 > \cdots > b_k - r_k > 0.
\end{equation}
The \bem{ladder determinantal ideal} $I(\ub a, \ub b, \ub r)$ is
generated by the minors of size~$r_\ell$ in the northwest $a_\ell \times
b_\ell$ corner of $Z$ for all $\ell \in 1,\ldots,k$.
Condition~(\ref{eq:ranks}) simply ensures that the vanishing of the
minors of size $r_\ell$ in the northwest $b_\ell \times a_\ell$
submatrix does not imply the vanishing of the minors of size~$r_{\ell'}$
in the northwest $b_{\ell'} \times a_{\ell'}$ submatrix when $\ell \neq
\ell'$.  For example, the ladder determinantal ideals in
\cite{GonLakLadDetSchub} have ranks~$\ub r$ that weakly increase from
southwest to northeast (in our language), while those treated in
\cite{GoncMilMixLadDet} have ranks such that no two labeled boxes lie in
the same row or column (that is, $(b_\ell,a_\ell)_{\ell = 1}^k$ is an
antidiagonal).

Since `ladders' are just another name for `partitions' one might also
like to call these `partition determinantal ideals', but in fact, a
better name is `vexillary determinantal ideals'.  Indeed, Fulton
identified ladder determinantal ideals as Schubert determinantal ideals
$I_w$ for \bem{vexillary permutations} (also known as
\bem{$2143$-avoiding} and \bem{single-shaped} permutations)
\cite[Proposition~9.6]{FulDegLoc}.  Therefore Theorems~\ref{t:formulae}
and~\ref{t:gb} hold in full for ladder determinantal ideals.  Note,
however (as Fulton does), that the probability of a permutation being
vexillary decreases exponentially to zero as $n$ approaches infinity.

Our Gr\"obner bases are new even for the vexillary determinantal ideals
we found in the literature, since previous authors seem to always use
what in our notation are diagonal rather than antidiagonal term orders.
The general phenomenon making the diagonal Gr\"obner basis fail as in
Example~\ref{ex:2143'} is precisely the fact that rank conditions are
``nested'' for every Schubert determinantal ideal that is not vexillary
(this follows from Fulton's essential set characterization
\cite[Section~9]{FulDegLoc}).

The multidegree formula in Theorem~\ref{t:formulae} becomes beautifully
explicit for vexillary ideals.

\begin{cor} \label{cor:ladder}
The $\ZZ^{2n}$-graded multidegree of a ladder determinantal variety is a
\bem{multi-Schur polynomial}, and therefore has an explicit determinantal
expression.
\end{cor}

This corollary is substantially more general than previous $\ZZ$-graded
degree formulae, which held only for special kinds of vexillary ideals,
and were after all only $\ZZ$-graded.  Readers wishing to see the
determinantal expression in its full glory can check \cite{FulDegLoc}
for a brief introduction to multi-Schur polynomials, or
\cite{NoteSchubPoly} for much more.  It would be desirable to make the
Hilbert series in Theorem~\ref{t:gb} just as explicit in closed form,
given that combinatorial formulae are known (see \cite{KMsubword} for
details and references):

\begin{question}
Is there an analogously ``nice'' formula%
	\footnote{In the introductions to \cite{AbhyankarECYT} and its
	second chapter, Abhyankar writes of formulae he first presented
	at a conference at the University of Nice, in France.  Although
	his formulae enumerate certain kinds of tableaux, his results
	were used to obtain formulae for degrees and Hilbert series of
	determinantal ideals.  Since then, some authors have been
	looking for ``nice'' (uncapitalized, and always in quotes)
	formulae for Hilbert series of determinantal ideals; cf.\
	\cite[p.~3]{HerzogTrung} and \cite[p.~55]{AbhyanKulkarHilb}.}
for vexillary double Grothendieck polynomials $\GG_w(\xx,\yy)$, or an
ordinary version for $\GG_w(\xx)$, or even for $\GG_w(t,\ldots,t)$?
\end{question}

The answer is `yes' for $\GG_w(t,\ldots,t)$ in certain vexillary cases;
e.g.\ see \cite{CH94,KP99,Gho02}.

Theorem~\ref{t:gb} provides a new proof that Schubert varieties $X_w
\subseteq \FL$ in the flag manifold are Cohen--Macaulay.
%
%
Instead of giving the quick derivation of this specific consequence, let
us instead mention a more general local equivalence principle between
Schubert varieties and matrix Schubert varieties, special cases of which
have been applied numerous times in the literature.  By a \bem{local
condition}, we mean a condition that holds for a variety whenever it
holds on each subvariety in some open cover.

\begin{thm} \label{thm:equiv}
Let $\mathfrak C$ be a local condition that holds for a variety $X$
whenever it holds for the product of $X$ with any vector space.  Then
$\mathfrak C$ holds for every Schubert variety in every flag variety if
and only if $\mathfrak C$ holds for all matrix Schubert varieties.
\end{thm}
\begin{sketch}
The complete proof of one direction is easy: if $\mathfrak C$ holds for
the matrix Schubert variety $\ol X_w \subseteq \mn$, then it holds for
$\wt X_w = \ol X_w \cap \gln$.  Therefore $\mathfrak C$ holds for the
Schubert variety $X_w \subseteq \bgln$, because $\wt X_w$ is locally
isomorphic to the product of $X_w$ with~$B$, the latter being an open
subset of a vector space.

On the other hand, if $\mathfrak C$ holds for Schubert varieties, then it
holds for matrix Schubert varieties because of the following, whose proof
(which uses Fulton's essential set \cite{FulDegLoc} and would require
introducing a fair amount of notation) we omit.  Given $w \in S_n$,
consider $w$ as an element of $S_{2n}$ fixing $n+1,\ldots,2n$.  The
product $\ol X_w \times \kk^{n^2-n}$ of the matrix Schubert variety $\ol
X_w \subseteq \mn$ with a vector space of dimension $n^2-n$ is isomorphic
to the intersection of the Schubert variety $X_w \subseteq {\mathcal
F}{\ell}_{2n}$ with the opposite big cell in ${\mathcal F}{\ell}_{2n}$.%
\end{sketch}

Since rationality of singularities and normality are among such local
statements, $\ol X_w$ possesses these properties because Schubert
varieties do \cite{ramanathanCM,RamRamNormSchub}.  However, these
statements could just as easily have been derived by Fulton in
\cite{FulDegLoc}, although they are not mentioned there.  Thus we know of
no new results that can be proved using Theorem~\ref{thm:equiv}.

%

\part{Bruhat induction}

\label{part:gb}

\section{Overview}\label{sec:gb}

With motivation coming from the statements of the main results in
Part~\ref{part:intro}, we now introduce the details of \bem{Bruhat
induction} to combinatorial commutative algebra.

Geometric considerations occupy Section~\ref{sec:multimat}, where we
start with large matrix Schubert varieties (associated to shorter
permutations) and chop them rather bluntly with multihomogeneous
functions (certain minors, actually).  This basically yields matrix
Schubert varieties that have dimension one less, which we understand
already by Bruhat induction.  However, the messy hypersurface section
leaves some debris components, which get cleaned up using the technology
of multidegrees (Sections~\ref{sec:mult} and~\ref{sec:pos}).  In
particular, multidegrees allow us to ignore geometric phenomena at
codimension $> 1$.

Beginning in Section~\ref{sec:mutation} we switch to a different track,
namely the combinatorial algebra of antidagonal ideals.  Bruhat
induction manifests itself here via Demazure operators, which we
interpret as combinatorial rules for manipulating Hilbert series
monomial by monomial.  Thus, in Section~\ref{sec:lift} we justify
certain $\Znn$-graded lifts of Demazure operators that take monomials
outside $J_w$ as input and return sums of monomials outside~$J_{ws_i}$.
The resulting Theorem~\ref{thm:ev} is substantially stronger than the
Hilbert series statement required for the main theorems, given that
these operators really do ``lift'' the $\ZZ^{2n}$-graded Demazure
operators to~$\Znn$.  This lifting property is proved in
Section~\ref{sec:coarsen} (Theorem~\ref{thm:induction}) via a certain
combinatorial duality on standard monomials of antidiagonal ideals.

The transition from arbitrary standard monomials to squarefree
monomials, which coincide with faces of the Stanley--Reisner complex,
starts with Section~\ref{sec:Lw}, where a rough interpretation of Bruhat
induction on facets of~$\LL_w$ already proves equidimensionality.  More
detailed (but less algebraic and more pleasantly combinatorial) analysis
in Section~\ref{sec:facets} reduces monomial-by-monomial Demazure
induction to facet-by-facet divided difference Bruhat induction.  It
culminates in mitosis for facets, which is the combinatorial residue of
Bruhat induction for standard monomials of antidiagonal ideals.  Further
use of mitosis in Section~\ref{sec:rp} characterizes the facets of
antidiagonal complexes as reduced pipe dreams.


The final Section~\ref{sec:proof} gathers the main results of Bruhat
induction into a proof of the remaining unproved assertions from
Part~\ref{part:intro} (Theorems~\ref{t:formulae}, \ref{t:gb},
and~\ref{t:mitosis}).  Logically, the remainder of this paper depends
only on the definitions in Part~\ref{part:intro}, and on
Section~\ref{sec:pos}.

\subsection*{Conventions for downward induction on Bruhat order}
We shall repeatedly invoke the hypothesis $\length(ws_i) < \length(w)$.
In terms of permutation matrices, this means that $w^T$ differs from
$(ws_i)^T$ only in rows $i$ and $i+1$, where they look heuristically like
\begin{equation} \label{eq:w}
\begin{array}{c@{\qquad\quad}c}
\begin{array}{lr|c|c|@{}c@{}|c|c|l}
  \multicolumn{6}{c}{}&\multicolumn{1}{c}{
        \begin{array}{@{}c@{}}
                \makebox[0pt]{$\scriptstyle w(i)$}\\
                \downarrow
        \end{array}}&
  \multicolumn{1}{c}{\begin{array}{@{}c@{}}\\\adots\end{array}}
  \\\cline{3-7}
        i  &&   & \phantom{1} & \toplinedots & \phantom{1} & 1
  \\\cline{3-4}\cline{6-7}
        i+1&& 1 &             &              &             &
  \\\cline{3-7}
  \multicolumn{1}{c}{}&
  \multicolumn{1}{c}{\begin{array}{@{}c@{}}\adots\\ \\ \end{array}}&
  \multicolumn{1}{c}{
        \begin{array}{@{}c@{}}
                \uparrow\\
                \makebox[0pt]{$\scriptstyle w(i+1)$}
        \end{array}}
  \\[-.5ex]
  \multicolumn{2}{c}{}&
  \multicolumn{5}{c}{w^T}
\end{array}
&
\begin{array}{lr|c|c|@{}c@{}|c|c|l}
  \multicolumn{2}{c}{}&
  \multicolumn{1}{c}{
        \begin{array}{@{}c@{}}
                \makebox[0pt]{$\scriptstyle w(i+1)$}\\
                \downarrow
        \end{array}}&
  \multicolumn{4}{c}{}&
  \multicolumn{1}{c}{\begin{array}{@{}c@{}}\\\adots\end{array}}
  \\\cline{3-7}
        i  && 1 & \phantom{1} & \toplinedots & \phantom{1} &
  \\\cline{3-4}\cline{6-7}
        i+1&&   &             &              &             & 1
  \\\cline{3-7}
  \multicolumn{1}{c}{}&
  \multicolumn{1}{c}{\begin{array}{@{}c@{}}\adots\\\\\end{array}}&
  \multicolumn{4}{c}{}&\multicolumn{1}{c}{
        \begin{array}{@{}c@{}}
                \uparrow\\
                \makebox[0pt]{$\scriptstyle w(i)$}
        \end{array}}
  \\[-.5ex]
  \multicolumn{2}{c}{}&
  \multicolumn{5}{c}{(ws_i)^T}
\end{array}
\end{array}
\end{equation}
between columns $w(i+1)$ and $w(i)$.  Since reversing the inequality
$\length(ws_i) < \length(w)$ makes so much difference, we always write
the hypothesis this way, for consistency, even though we may actually
{\em use}\/ one of the following equivalent formulations in any given
lemma or proposition.  We hope that collecting this list of standard
statements (``shorter permutation $\iff$ bigger variety'') will prevent
the reader from stumbling on this as many times as we did.  The string of
characters `$\length(ws_i) < \length(w)$' can serve as a visual cue to
this frequent assumption; we shall {\em never}\/ assume the opposite
inequality.

\begin{lemma} \label{cor:bruhat}
The following are equivalent for a permutation $w \in S_n$.
$$
\begin{array}{rl@{\qquad}rl}
  1.&\length(ws_i) < \length(w).    &6.&\dim(\ol X_{ws_i}) > \dim(\ol X_w).
\\2.&\length(ws_i) = \length(w) - 1.&7.&\dim(\ol X_{ws_i}) = \dim(\ol X_w) + 1.
\\3.&w(i) > w(i+1).                 &8.&s_i \ol X_w \neq \ol X_w.
\\4.&ws_i(i) < ws_i(i+1).           &9.&s_i \ol X_{ws_i} = \ol X_{ws_i}.
\\5.&I(\ol X_{ws_i}) \subset I(\ol X_w).&10.&\ol X_{ws_i} \supset \ol X_w.
\end{array}
$$
Here, the transposition $s_i$ acts on the left of $\mn$, switching rows
$i$ and~$i+1$.
\end{lemma}
\begin{proof}
The equivalence of $1$--$4$ comes from Claim~\ref{claim:length}.  The
equivalence of these with $5$--$7$ and~$10$ uses
Proposition~\ref{prop:BwB} and Lemma~\ref{lemma:bruhat}.  Finally,
$8$--$10$ are equivalent by Lemma~\ref{lemma:Z} and its proof.  (We shall
not apply Lemma~\ref{cor:bruhat} until the proof of
Lemma~\ref{lemma:si}.)%
\end{proof}

\section{Multidegrees of matrix Schubert varieties}\label{sec:multimat}

This section provides a proof of the divided difference recursion
satisfied by the multidegrees of matrix Schubert varieties, in
Theorem~\ref{thm:oracle}.  Although many of the preliminary results can
be deduced from a number of places in the literature, notably
\cite{FulDegLoc}, we believe it important to provide proofs (or at least
sketches) so as to make our foundations explicit.

Matrix Schubert varieties are clearly stable under rescaling any row or
column.  Moreover, since we only impose rank conditions on submatrices
that are as far north and west as possible, any operation that adds a
multiple of some row to a row below it (``sweeping downward''), or that
adds a multiple of some column to another column to its right (``sweeping
to the right'') preserves every matrix Schubert variety.

In terms of group theory, let $B$ denote the group of invertible {\em
lower}\/ triangular matrices and $B_+$ the invertible upper triangular
matrices.  The previous paragraph says exactly that each matrix Schubert
variety $\ol X_w$ is preserved by the action%
	\footnote{This is a left group action, in the sense that
	$(b,b_+) \cdot ((b',b'_+) \cdot Z)$ equals
	$((b,b_+)\cdot(b',b'_+))\cdot Z$ instead of
	$((b',b'_+)\cdot(b,b_+))\cdot Z$, even though---in fact
	because---the $b_+$ acts via its inverse \hbox{on the {\em right}}.}
of $B\times B_+$ on $\mn$ in which $(b,b_+) \cdot Z = bZb_+^{-1}$.
Proposition~\ref{prop:BwB} will say more.

The next four results, numbered \ref{lemma:orbits}--\ref{lemma:Z}, are
basically standard facts concerning Bruhat order for~$S_n$, enhanced
slightly for $\mn$ instead of~$\gln$.  They serve as prerequisites for
Proposition~\ref{prop:BwB}, Lemma~\ref{lemma:si}, and
Lemma~\ref{lemma:regular}, which enter at key points in the proof of
Theorem~\ref{thm:oracle}.

Call a matrix $Z \in \mn$ that is zero except for at most one~$1$ in
each row and column a \bem{partial permutation matrix}.  These arise in
the $\mn$ analogue of Bruhat decomposition:

\begin{lemma} \label{lemma:orbits}
In each $B \times B_+$ orbit on $\mn$ lies a unique partial permutation
matrix.
\end{lemma}
\begin{proof}
By doing row and column operations that sweep down and to the right, we
can get from an arbitrary matrix~$Z'$ to a partial permutation
matrix~$Z$.  Such sweeping preserves the ranks of northwest $q \times p$
submatrices, and~$Z$ can be reconstructed uniquely by knowing only
$\rank(Z'\sub qp)$ for $1 \leq q,p \leq n$.
\end{proof}

Define the \bem{length} of a partial permutation matrix~$Z$ as the number
of zeros in $Z$ that lie neither due south nor due east of a~$1$.  In
other words, for every $1$ in~$Z$, cross out all the boxes beneath it in
the same column as well as to its right in the same row, and count the
number of uncrossed-out boxes to get the length of~$Z$.  When $Z = w^T$
is a permutation matrix, $\length(w)$ agrees with the length of~$w^T \in
\mn$.  Write $Z \subseteq Z'$ for partial permutation matrices $Z$
and~$Z'$ if the $1$'s in $Z$ are a subset of the $1$'s in~$Z'$.  Finally,
let $t_{i,i'} \in S_n$ be the transposition switching $i$ and~$i'$.  The
following claim is self-evident.

\begin{claim} \label{claim:length}
Suppose $Z$ is a partial permutation matrix with $1$'s at $(i,j)$ and
$(i',j')$, where $(i,j) \leq (i',j')$.  Switching rows $i$ and $i'$
of~$Z$ creates a partial permutation matrix $Z'$ satisfying\/
$\length(Z') = \length(Z) + 1 + \mbox{}$twice the number of~$1$'s
strictly inside the rectangle enclosed by $(i,j)$ and
$(i',j')$.\hfill$\square$
\end{claim}

\begin{lemma} \label{lemma:bruhat}
Fix a permutation $w \in S_n$ and a partial permutation matrix~$Z$.  If
$Z \in \ol X_w$ and $\length(Z) \geq \length(w)$, then either $Z
\subseteq w^T$, or there is a transposition $t_{i,i'}$ such that $v =
wt_{i,i'}$ satisfies: $Z \in \ol X_v$ and\/ $\length(v) > \length(w)$.
\end{lemma}
\begin{sketch}
Use reasoning similar to the case when $Z$ is a permutation matrix: work
by downward induction on the number of $1$'s shared by $w^T$ and~$Z$,
using Claim~\ref{claim:length}.  We omit the details.
\end{sketch}

\begin{lemma} \label{lemma:Z}
Let $Z$ be a partial permutation matrix with orbit closure $\ol\OO_Z$ in
$\mn$.  If\/ $\length(s_i Z) < \length(Z)$, then
$$
  n^2 - \length(s_i Z)\ \:=\ \:\dim(\ol\OO_{s_i Z})\ \:=\ \:
  \dim(\ol\OO_Z) + 1\ \:=\ \:n^2 - \length(Z) + 1,
$$
and $s_i(\ol \OO_{s_i Z}) = \ol\OO_{s_i Z}$.
\end{lemma}
\begin{proof}
Let $P_i \subseteq \gln$ be the $i^\th$ parabolic subgroup containing
$B$, in which the only nonzero entry outside the lower triangle may lie
at $(i,i+1)$.  Consider the image $Y$ of the multiplication map $P_i
\times \ol\OO_Z \to \mn$ sending $(p,x) \mapsto p \cdot x$.  This map
factors through the quotient $B\dom (P_i \times \ol\OO_Z)$ by the
diagonal action of~$B$, which is an $\ol\OO_Z$-bundle over $P_i/B \cong
\PP^1$ and hence has dimension $\dim(\ol\OO_Z) + 1$.  Thus $\dim(Y) \leq
\dim(\ol\OO_Z) + 1$.

The variety $Y$ is $B \times B_+$-stable by construction.  Since $B
\times B_+$ has only finitely many orbits on~$\mn$, the irreducibility
of~$Y$ implies that $\ol Y$ is an orbit closure $\ol\OO_{Z'}$ for some
partial permutation matrix~$Z'$, by Lemma~\ref{lemma:orbits}.  Clearly
$\ol\OO_Z \subseteq Y$, and $s_i Z \in Y$, so the dimension bound implies
that $Z' = s_i Z$ because $Z \in \ol\OO_{s_i Z}$, as can be checked
directly.  This $\PP^1$-bundle argument also shows that $\ol\OO_{s_i Z} =
\ol Y$ is stable under multiplication by $P_i$ on the left, whence $s_i
\in P_i$ takes $\ol\OO_{s_i Z}$ to itself.

That $\dim(\ol\OO_Z) = n^2 - \length(Z)$ follows by direct calculation
whenever $Z$ has nonzero entries only along the main diagonal, and then
by downward induction on $\length(Z)$.%
\end{proof}

Fulton derived the next result in~\cite{FulDegLoc} from the
corresponding result on flag manifolds.

\begin{prop} \label{prop:BwB}
The matrix Schubert variety $\ol X_w$ is the closure $\ol{B w^T B_+}$ of
the \hbox{$B \times B_+$} orbit on~$\mn$ through the permutation matrix
$w^T$.  Thus $\ol X_w$ is irreducible of dimension $n^2 - \length(w)$,
and $w^T$ is a smooth point of it.
\end{prop}
\begin{proof}
The stability of $\ol X_w$ under $B \times B_+$ means that $\ol X_w$ is a
union of orbits.  By the obvious containment $\ol\OO_{w^T} \subseteq \ol
X_w$ (sweeping down and right preserves northwest ranks) and
Lemma~\ref{lemma:orbits}, it suffices to show that partial permutation
matrices $Z$ lying in $\ol X_w$ lie also in $\ol{B w^T B_+}$.  Standard
arguments analogous to the case where $Z$ is a permutation matrix work
here, using Lemma~\ref{lemma:bruhat}.  The last sentence of the
proposition is standard for orbit closures, except for the dimension
count, which comes from Lemma~\ref{lemma:Z}.
\end{proof}


\begin{lemma} \label{lemma:si}
Let $Z$ be a partial permutation matrix and $w \in S_n$.  If the orbit
closure $\ol\OO_Z$ has codimension~$1$ inside $\ol X_{ws_i}$, then
$\ol\OO_Z$ is mapped to itself by $s_i$ unless~$Z = w$.
\end{lemma}
\begin{proof}
First note that $\length(Z) = \length(ws_i) + 1$ by Lemma~\ref{lemma:Z}
and Proposition~\ref{prop:BwB}.  Using Lemma~\ref{lemma:bruhat} and
Claim~\ref{claim:length}, we find that $Z$ is obtained from $(ws_i)^T$ by
switching some pairs of rows to make partial permutations of strictly
larger length and then deleting some $1$'s.  Since the length of $Z$ is
precisely one less than that of~$ws_i$, we can switch exactly one pair of
rows of~$(ws_i)^T$, or we can delete a single~$1$ from $(ws_i)^T$.

Any $1$ that we delete from $(ws_i)^T$ must have no $1$'s southeast of
it, or else the length increases by more than one.  Thus the~$1$ in
row~$i$ of $(ws_i)^T$ cannot be deleted by statement~4 of
Lemma~\ref{cor:bruhat}, leaving us in the situation of
Lemma~\ref{lemma:Z} with~$Z = w^T$, and completing the case where a 1
has been deleted from $(ws_i)^T$ to get~$Z$.

Suppose now that switching rows $q$ and~$q'$ of $(ws_i)^T$ results in
the matrix $Z = v^T$ for some permutation~$v$, and assume that $s_i(\ol
\OO_Z) \neq \ol\OO_Z$.  Since $s_i(\ol \OO_Z) = \ol\OO_Z$ unless $v$
satisfies $v(i) > v(i+1)$, by Lemma~\ref{lemma:Z}, we find that $v(i) >
v(i+1)$.  At least one of $q$ and~$q'$ must lie in $\{i,i+1\}$ by part~4
of Corollary~\ref{cor:bruhat}, which says that moving neither row~$q$
nor row~$q'$ of~$(ws_i)^T$ leaves $v(i) < v(i+1)$.  On the other hand,
it is impossible for exactly one of $q$ and~$q'$ to lie in $\{i,i+1\}$;
indeed, switching rows $q$ and~$q'$ increases length,
so either the~$1$ at $(i,w(i+1))$ or the~$1$ at $(i+1,w(i))$ would lie
inside the rectangle formed by the switched $1$'s, making $\ol\OO_Z$
have codimension more than~$1$ by Claim~\ref{claim:length} and
Proposition~\ref{prop:BwB}.  Thus $\{q,q'\} = \{i,i+1\}$ and $v = w$,
completing the proof.%
\end{proof}

\begin{lemma} \label{lemma:regular}
If\/ $\length(ws_i) < \length(w)$, and $\mm_{ws_i}$ is the maximal ideal
in the local ring of $(ws_i)^T \in \ol X_{ws_i}$, then the variable
$z_{i+1,w(i+1)}$ maps to a regular parameter in~$\mm_{ws_i}$.  In other
words $z_{i+1,w(i+1)}$ lies in $\mm_{ws_i} \minus \mm_{ws_i}^2$.
\end{lemma}
\begin{proof}
Let $v = ws_i$, and consider the map $B \times B_+ \to \mn$ sending
$(b,b^+) \mapsto b \cdot v^T \cdot b^+$.  The image of this map is
contained in $\ol X_v$ by Proposition~\ref{prop:BwB}, and the identity
$\id := (\id_B,\id_{B_+})$ maps to $v^T$.  The induced map of local rings
the other way thus takes $\mm_v$ to the maximal~ideal
\begin{eqnarray*}
  \mm_\id &:=& \<b_{ii}-1,b^+{}_{\!\!\!\!\!ii}-1 \mid 1 \leq i \leq n\>
  + \<b_{ij},b^+{}_{\!\!\!\!\!ji} \mid i > j\>
\end{eqnarray*}
in the local ring at the identity $\id \in B \times B_+$.  It is enough
to demonstrate that the image of $z_{i+1,w(i+1)}$ lies in $\mm_\id \minus
\mm_\id^2$.

Direct calculation shows that $z_{i+1,w(i+1)}$ maps to 
\begin{eqnarray*}
  b_{i+1,i}b^+{}_{\!\!\!\!\!w(i+1),w(i+1)} + \sum_{q \in Q}
  b_{i+1,q}b^+{}_{\!\!\!\!\!q,w(i+1)} \quad\hbox{where}\quad Q = \{q < i
  \mid w(q) < w(i+1)\}.
\end{eqnarray*}
In particular, all of the summands $b_{i+1,q}b^+{}_{\!\!\!\!\!q,w(i+1)}$
lie in $\mm_\id^2$.  On the other hand, $b^+{}_{\!\!\!\!\!w(i+1),w(i+1)}$
is a unit near the identity, so
$b_{i+1,i}b^+{}_{\!\!\!\!\!w(i+1),w(i+1)}$ lies in $\mm_\id \minus
\mm_\id^2$.
\end{proof}

The argument forming the proof of the previous lemma is an alternative
way to calculate the dimension as in Lemma~\ref{lemma:Z}.

\begin{thm}\label{thm:oracle}
If the permutation $w$ satisfies $\length(ws_i) < \length(w)$, then
\begin{eqnarray*}
  [\ol X_{ws_i} ] &=& \partial_i [\ol X_w]
\end{eqnarray*}
holds for both the $\ZZ^{2n}$-graded and $\ZZ^n$-graded multidegrees.
\end{thm}
\begin{proof}
The proof here works for the $\ZZ^n$-grading as well as the
$\ZZ^{2n}$-grading, simply by ignoring all occurrences of $\yy$, or
setting them to zero.

Let $j = w(i) - 1$, and suppose $\rank(w^T\sub ij) = r-1$.  Then the
permutation matrix $(ws_i)^T$ has $r$ entries equal to $1$ in the
submatrix $(ws_i)^T\sub ij$.  Consider the $r \times r$ minor $\Delta$
using the rows and columns in which $(ws_i)^T\sub ij$ has $1$'s.  Thus
$\Delta$ is not the zero function on $\ol X_{ws_i}$; in fact, $\Delta$
is nonzero everywhere on its interior $B (ws_i)^T B_+$.  Therefore the
subscheme $X_\Delta$ defined by $\Delta$ inside $\ol X_{ws_i}$ is
supported on a union of orbit closures $\ol\OO_Z$ contained in $\ol
X_{ws_i}$ with codimension~$1$.  Now we compare the subscheme $X_\Delta$
to its image $s_i X_\Delta = X_{s_i \Delta}$ under switching rows $i$
\mbox{and~$i+1$}.

\begin{claim} \label{claim:compare}
Every irreducible component of $X_\Delta$ other than~$\ol X_w$ has the
same multiplicity in~$s_i X_\Delta$, and $\ol X_w$ has multiplicity~$1$
in~$X_\Delta$.
\end{claim}
\begin{proof}
Lemma~\ref{lemma:si} says that $s_i$ induces an automorphism of the
local ring at the generic point (i.e.\ the prime ideal) of $\ol\OO_Z$
inside~$\ol X_{ws_i}$, for every irreducible component $\ol\OO_Z$ of
$X_\Delta$ other than~$\ol X_w$.  This automorphism takes $\Delta$
to~$s_i\Delta$, so these two functions have the same multiplicity
along~$\ol\OO_Z$.  The only remaining codimension~$1$ irreducible
component of $X_\Delta$ is~$\ol X_w$, and we shall now verify that the
multiplicity equals~$1$ there.  As a consequence, the multiplicity of
$s_i X_\Delta$ along $s_i \ol X_w$ also equals~$1$.

By Proposition~\ref{prop:BwB}, the local ring of $(ws_i)^T$ in~$\ol
X_{ws_i}$ is regular.  Since $s_i$ is an automorphism of~$\ol X_{ws_i}$,
we find that the local ring of $w^T \in \ol X_{ws_i}$ is also regular.
In a neighborhood of $w^T$, the variables~$z_{qp}$ corresponding to the
locations of the $1$'s in $w^T\sub ij$ are units.  This implies that the
coefficient of $z_{i,w(i+1)}$ in~$\Delta$ is a unit in the local ring of
$w^T \in \ol X_{ws_i}$.  On the other hand, the set of variables in
spots where $w^T$ has zeros generate the maximal ideal in the local ring
at $w^T \in \ol X_{ws_i}$.  Therefore, all terms of~$\Delta$ lie in the
square of this maximal ideal, except for the unit times~$z_{i,w(i+1)}$
term produced above.  Hence, to prove multiplicity one, it is enough to
prove that $z_{i,w(i+1)}$ itself is a regular parameter at $w^T \in \ol
X_{ws_i}$, or equivalently (after applying $s_i$) that $z_{i+1,w(i+1)}$
is a regular parameter at $(ws_i)^T \in \ol X_{ws_i}$.  This is
Lemma~\ref{lemma:regular}.
\end{proof}

Now we use a multidegree trick.  Consider $\Delta$ and $s_i\Delta$ as
elements not in $\kk[\zz]$, but in the ring $\kk[\zz,u]$ with $n^2 + 1$
variables, where the ordinary weight of the new variable~$u$ is $x_i -
x_{i+1}$.  Denote by $\mn \times \AA^1$ the spectrum of~$\kk[\zz,u]$.
Then $\Delta$ and the product $u s_i\Delta$ in $\kk[\zz,u]$ have the
{\em same}\/ ordinary weight $f := f(\xx,\yy)$.  Since the affine
coordinate ring of $\ol X_{ws_i}$ is a domain, neither $\Delta$ nor
$s_i\Delta$ vanishes on $\ol X_{ws_i}$, so we get two short exact
sequences
\begin{equation} \label{eq:Q}
  0\ \to\ \kk[\zz,u]/I(\ol X_{ws_i})\kk[\zz,u](-f)\
  \stackrel\Theta\too\ \kk[\zz,u]/I(\ol X_{ws_i})\kk[\zz,u]\ \too\
  Q(\Theta)\ \to\ 0,
\end{equation}
in which $\Theta$ equals either $\Delta$ or $us_i\Delta$.  The quotients
$Q(\Delta)$ and $Q(us_i\Delta)$ therefore have equal $\ZZ^{2n}$-graded
Hilbert series, and hence equal multidegrees.

Note that $Q(\Delta)$ is the coordinate ring of $X_\Delta \times \AA^1$,
while $Q(us_i\Delta)$ is the coordinate ring of $(s_iX_\Delta \times
\AA^1) \cup (\ol X_{ws_i} \times \{0\})$, the latter component being the
zero scheme of $u$ in $\kk[\zz,u]/I(\ol X_{ws_i})\kk[\zz,u]$.  Breaking
up the multidegrees of $Q(\Delta)$ and $Q(us_i\Delta)$ into sums over
irreducible components by additivity in Theorem~\ref{t:multidegs},
Claim~\ref{claim:compare} says that almost all terms in the equation
\begin{eqnarray*}
  [X_\Delta \times \AA^1] &=& [s_i X_\Delta \times \AA^1] + [\ol X_{ws_i}
  \times \{0\}]
\end{eqnarray*}
cancel, leaving us only with
\begin{eqnarray} \label{eq:cancel}
  [\ol X_w \times \AA^1] &=& [s_i\ol X_w \times \AA^1] + [\ol X_{ws_i}
  \times \{0\}].
\end{eqnarray}

The brackets in these equations denote multidegrees over $\kk[\zz,u]$.
However, the ideals in $\kk[\zz,u]$ of $\ol X_w \times \AA^1$ and
$s_i\ol X_w \times \AA^1$ are extended from the ideals in $\kk[\zz]$ of
$\ol X_w$ and $s_i \ol X_w$.  Therefore their $K$-polynomials agree with
those of $\ol X_w$ and~$s_i \ol X_w$, respectively, whence $[\ol X_w
\times \AA^1] = [\ol X_w]$ and $[s_i\ol X_w \times \AA^1] = [s_i\ol X_w]$
as polynomials in $\xx$ and~$\yy$.  The same argument shows that $[\ol
X_{ws_i} \times \AA^1] = [\ol X_{ws_i}]$.  The coordinate ring of $[\ol
X_{ws_i} \times \{0\}]$, on the other hand, is the right hand term of the
exact sequence that results after replacing $f$ by $x_i - x_{i+1}$ and
$\Theta$ by $u$ in~(\ref{eq:Q}).  We therefore find that
\begin{eqnarray*}
  [\ol X_{ws_i} \times \{0\}] &=& (x_i - x_{i+1})[\ol X_{ws_i} \times
  \AA^1]\ \:=\ \:(x_i - x_{i+1})[\ol X_{ws_i}]
\end{eqnarray*}
as polynomials in $\xx$ and $\yy$.  Substituting back
into~(\ref{eq:cancel}) yields the equation on multidegrees $[\ol X_w] =
[s_i\ol X_w] + (x_i-x_{i+1})[\ol X_{ws_i}]$, which produces the desired
result after moving the $[s_i\ol X_w]$ to the left and dividing through
by $x_i - x_{i+1}$.
\end{proof}

\begin{remark}
This proof, although translated into the language of multigraded
commutative algebra, is actually derived from a standard proof of
divided difference formulae by localization in equivariant cohomology,
when $\kk = \CC$.  The connection is our two functions $\Delta$
and~$s_i\Delta$, which yield a map $\ol X_{ws_i} \to \CC^2$.  The
preimage of one axis is $\ol X_w$ union some junk components, and the
preimage of the other axis is $s_i \ol X_w$ union the same junk
components.  Therefore all of the unwanted (canceling) contributions map
to the point $(0,0) \in \CC^2$.  Essentially, the standard equivariant
localization proof makes the map to $\CC^2$ into a map to $\CC\PP^1$,
thus avoiding the extra components, and pulls back the localization
formula on $\CC\PP^1$ to a formula on whatever $\ol X_{ws_i}$ has become
(a Schubert variety).
\end{remark}

\section{Antidiagonals and mutation}\label{sec:mutation}

In this section we begin investigating the combinatorial properties of
the monomials outside $J_w$ and the antidiagonals generating~$J_w$.  For
the rest of this section, fix a permutation $w$ and a transposition $s_i$
satisfying $\length(ws_i) < \length(w)$.

Define the \bem{rank matrix} $\rk(w)$ to have $(q,p)$ entry equal to
$\rank(w^T\sub qp)$.  There are two standard facts we need concerning
rank matrices, both proved simply by looking at the picture of $(ws_i)^T$
in~(\ref{eq:w}).
\begin{lemma} \label{lemma:rank}
Suppose the $(i,j)$ entry of\/ $\rk(ws_i)$ is $r$.
\begin{numbered}
\item \label{i-1} If $j \geq w(i+1)$ then the $(i-1,j)$ entry of\/
$\rk(ws_i)$ is $r-1$.

\item \label{i+1}
If $j < w(i+1)$ then the $(i+1,j)$ entry of\/ $\rk(ws_i)$ is $r$.
\end{numbered}
\end{lemma}

In what follows, a \bem{rank condition} refers to a statement
requiring ``$\rank(Z\sub qp) \leq r$'' for some $r \geq 0$.  Most
often, $r$ will be either $\rank(w^T\sub qp)$ or $\rank((ws_i)^T\sub
qp)$, thereby making the entries of $\rk(w)$ and $\rk(ws_i)$ into rank
conditions.  We say that a rank condition $\rank(Z\sub qp) \leq r$
\bem{causes} an antidiagonal $a$ of the generic matrix $Z$ if $Z\sub
qp$ contains $a$ and the number of variables in $a$ is strictly larger
than~$r$.  For instance, when the rank condition is in $\rk(w)$, the
antidiagonals it causes include those $a \in J_w$ that are contained
in $Z\sub qp$ but no smaller northwest submatrix.  Although
antidiagonals in $Z$ (that is, antidiagonals of square submatrices of
the generic matrix~$Z$) are by definition monomials, we routinely
identify each antidiagonal with its \bem{support}\/: the subset of the
\hbox{variables dividing it in $\kk[\zz]$.}

\begin{lemma} \label{lemma:i}
Antidiagonals in $J_w \minus J_{ws_i}$ are subsets of $Z\sub i{w(i)}$
and intersect row~$i$.
\end{lemma}
\begin{proof}
If an antidiagonal in $J_w$ is either contained in $Z\sub {i-1}{w(i)}$ or
not contained in $Z\sub i{w(i)}$, then some rank condition causing it is
in both $\rk(w)$ and $\rk(ws_i)$.  Indeed, it is easy to check that the
rank matrices $\rk(ws_i)$ and $\rk(w)$ differ only in row $i$ between the
columns $w(i+1)$ and $w(i)-1$, inclusive.
\end{proof}

Though simple, the next lemma is the key combinatorial observation.  Note
that the permutation $w$ there is arbitrary; in particular, we will
frequently apply the lemma in the context of antidiagonals for a
permutation called $ws_i$.
\begin{lemma} \label{lemma:squish}
Suppose $a \in J_w$ is an antidiagonal and $a' \subset Z$ is another
antidiagonal.
\begin{enumerate}
\item[{\makebox[3.5ex][c]{(W)}}] \label{west} If $a'$ is obtained by
moving west one or more of the variables in $\,a$, then $a' \in J_w$.

\item[{\makebox[3.5ex][c]{(E)}}] \label{east} If $a' \in \kk[\zz]$ is
obtained by moving east any variable {\em except}\/ the northeast one in
$\,a$, then $a' \in J_w$.

\item[{\makebox[3.5ex][c]{(N)}}] \label{north} If $a'$ is obtained by
moving north one or more of the variables in $\,a$, then $a' \in J_w$.

\item[{\makebox[3.5ex][c]{(S)}}] \label{south} If $a' \in \kk[\zz]$ is
obtained by moving south any variable {\em except}\/ the southwest one~in
$\,a$, then $a' \in J_w$.
\end{enumerate}
\end{lemma}
\begin{proof}
Every rank condition causing $a$ also causes all of the antidiagonals
$a'$.
\end{proof}

\begin{example} \label{ex:squish}
Parts~(W) and~(E) of Lemma~\ref{lemma:squish} together imply that the
type of motion depicted in the following diagram preserves the property
of an antidiagonal being in~$J_w$.  The presence of the northeast $*$
justifies moving the southwest $*$ east.
$$
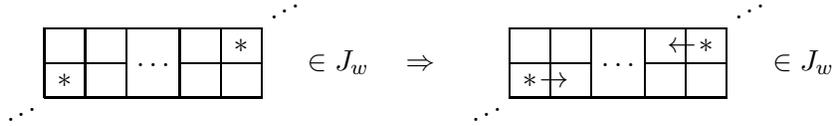
\begin{array}{@{}r@{\;}c@{\;}l@{}}
&&\adots\\&
\begin{array}{@{}l@{}|c|c|@{}c@{}|c|c|@{}}
  \cline{2-6}
        & &             & \toplinedots & \phantom{*} & *
  \\\cline{2-3}\cline{5-6}
        &*& \phantom{*} &              &             &
  \\\cline{2-6}
\end{array}\\
\adots
\end{array}
\in J_w \quad \implies \quad
\begin{array}{@{}r@{\;}c@{\;}l@{}}
&&\adots\\&
\begin{array}{@{}l@{}|c|c|@{}c@{}|c|c|@{}}
  \cline{2-6}
        & & & \toplinedots & \phantom{*}
                           & \makebox[0mm][r]{$\leftarrow$\hspace{1pt}}*
  \\\cline{2-3}\cline{5-6}
        &*\makebox[0mm][l]{\hspace{1pt}$\rightarrow$}
        &\phantom{*} &              & &
  \\\cline{2-6}
\end{array}\\
\adots
\end{array}
\in J_w
$$
The two rows could also be separated by some other rows---possibly
themselves containing elements of the original antidiagonal---as long as
the indicated motion preserves the fact that we have an antidiagonal.
\end{example}

\begin{defn} \label{defn:mu}
Let $\bb$ be an array $\bb = (b_{rs})$ of nonnegative integers, that is,
$\bb$ is an \bem{exponent array} for a monomial $\zz^\bb \in \kk[\zz]$.
Let
\begin{eqnarray*}
  \west_q(\bb) &:=& \min(\{p \mid b_{qp} \neq 0\}\cup\{\infty\})
\end{eqnarray*}
be the column of the leftmost (most ``western'') nonzero entry in row
$q$.  Define the \bem{mutation} of $\bb$ in rows~$i$ and~\mbox{$i+1$} by
\begin{eqnarray*}
  \mu_i(\bb) &:=& \hbox{\rm the exponent array of }
  ({z_{i,p}}/{z_{i+1,p}}) \zz^\bb \hbox{\rm\ for } p = \west_{i+1}(\bb).
\end{eqnarray*}
For ease of notation, we write $\mu_i(\zz^\bb)$ for
$\zz^{\mu_i(\bb)}$.
\end{defn}
If one thinks of the array $\bb$ as a chessboard with some coins stacked
in each square, then $\mu_i$ is performed by taking a coin off the
western stack in row $i+1$ and putting it onto the stack due north of it
in row $i$.

\begin{example} \label{ex:mutate}
Suppose that $\bb$ is the left array in Fig.~\ref{fig:mu}, and that $i =
3$.  We list in Fig.~\ref{fig:mu} (reading left to right as usual) 7
mutations of $\bb$, namely $\bb = (\mu_3)^0(\bb)$ through
$(\mu_3)^6(\bb)$ (after that it involves the dots we left unspecified).
Here, the empty boxes denote entries equal to $0$, and the nonzero
mutated entries at each step are in boxes.  To make things easier to
look at, the entries on or below the main antidiagonal are represented
by dots, each of which may be zero or not (independently of the others).
The $3$ and~$4$ at left are labels for rows $i = 3$ and $i+1 = 4$.
\begin{figure}[t]
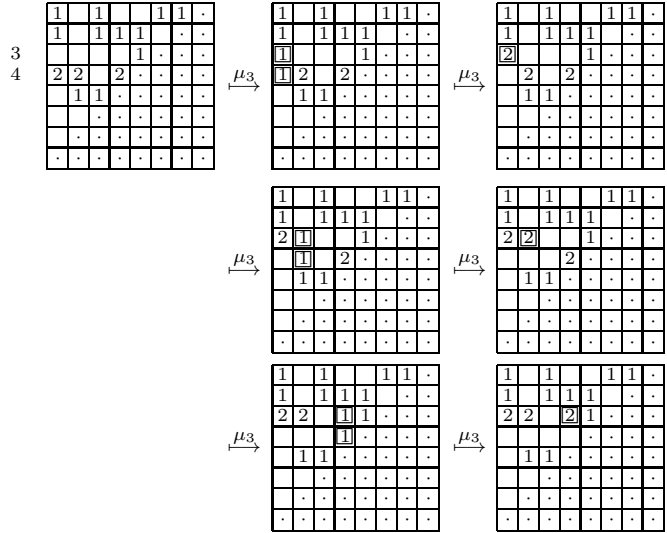

\def\mbf#1{\hbox{{\rlap{\makebox[0ex]{${#1}$}}}%
	{\raisebox{-.55ex}{\normalsize\hspace{-.9ex}$\Box$\hspace{-.85ex}}}}}
$$
\begin{tinyrc}{
\begin{array}{@{}cc@{}}{}\\\\3&\\4\\\\\\\\\\\end{array}
\begin{array}{@{}|@{\:}c@{\:}|@{\:}c@{\:}|@{\:}c@{\:}|@{\:}c@{\:}
  		 |@{\:}c@{\:}|@{\:}c@{\:}|@{\:}c@{\:}|@{\:}c@{\:}|@{}}
\hline    1  &     &  1  &     &     &  1  &  1  &\hspace{1pt}\cdot\hspace{1pt}
\hln      1  &     &  1  &  1  &  1  &     &\cdot&\cdot
\hln         &     &     &     &  1  &\cdot&\cdot&\cdot
\hln      2  &  2  &     &  2  &\cdot&\cdot&\cdot&\cdot
\hln         &  1  &  1  &\cdot&\cdot&\cdot&\cdot&\cdot
\hln         &     &\cdot&\cdot&\cdot&\cdot&\cdot&\cdot
\hln         &\cdot&\cdot&\cdot&\cdot&\cdot&\cdot&\cdot
\hln    \cdot&\cdot&\cdot&\cdot&\cdot&\cdot&\cdot&\cdot
\\\hline
\end{array}
\ \stackrel{\textstyle \mu_3}{\longmapsto}\:
\begin{array}{@{}|@{\:}c@{\:}|@{\:}c@{\:}|@{\:}c@{\:}|@{\:}c@{\:}
  		 |@{\:}c@{\:}|@{\:}c@{\:}|@{\:}c@{\:}|@{\:}c@{\:}|@{}}
\hline    1  &     &  1  &     &     &  1  &  1  &\hspace{1pt}\cdot\hspace{1pt}
\hln      1  &     &  1  &  1  &  1  &     &\cdot&\cdot
\hln    \mbf1&     &     &     &  1  &\cdot&\cdot&\cdot
\hln    \mbf1&  2  &     &  2  &\cdot&\cdot&\cdot&\cdot
\hln         &  1  &  1  &\cdot&\cdot&\cdot&\cdot&\cdot
\hln         &     &\cdot&\cdot&\cdot&\cdot&\cdot&\cdot
\hln         &\cdot&\cdot&\cdot&\cdot&\cdot&\cdot&\cdot
\hln    \cdot&\cdot&\cdot&\cdot&\cdot&\cdot&\cdot&\cdot
\\\hline
\end{array}
\ \stackrel{\textstyle \mu_3}{\longmapsto}\:
\begin{array}{@{}|@{\:}c@{\:}|@{\:}c@{\:}|@{\:}c@{\:}|@{\:}c@{\:}
  		 |@{\:}c@{\:}|@{\:}c@{\:}|@{\:}c@{\:}|@{\:}c@{\:}|@{}}
\hline    1  &     &  1  &     &     &  1  &  1  &\hspace{1pt}\cdot\hspace{1pt}
\hln      1  &     &  1  &  1  &  1  &     &\cdot&\cdot
\hln    \mbf2&     &     &     &  1  &\cdot&\cdot&\cdot
\hln         &  2  &     &  2  &\cdot&\cdot&\cdot&\cdot
\hln         &  1  &  1  &\cdot&\cdot&\cdot&\cdot&\cdot
\hln         &     &\cdot&\cdot&\cdot&\cdot&\cdot&\cdot
\hln         &\cdot&\cdot&\cdot&\cdot&\cdot&\cdot&\cdot
\hln    \cdot&\cdot&\cdot&\cdot&\cdot&\cdot&\cdot&\cdot
\\\hline
\end{array}
}\end{tinyrc}
$$
$$
\begin{tinyrc}{
\phantom{
\begin{array}{@{}cc@{}}{}\\\\3&\\4\\\\\\\\\\\end{array}
\begin{array}{@{}|@{\:}c@{\:}|@{\:}c@{\:}|@{\:}c@{\:}|@{\:}c@{\:}
  		 |@{\:}c@{\:}|@{\:}c@{\:}|@{\:}c@{\:}|@{\:}c@{\:}|@{}}
\hline    1  &     &  1  &     &     &  1  &  1  &\hspace{1pt}\cdot\hspace{1pt}
\hln      1  &     &  1  &  1  &  1  &     &\cdot&\cdot
\hln         &     &     &     &  1  &\cdot&\cdot&\cdot
\hln      2  &  2  &     &  2  &\cdot&\cdot&\cdot&\cdot
\hln         &  1  &  1  &\cdot&\cdot&\cdot&\cdot&\cdot
\hln         &     &\cdot&\cdot&\cdot&\cdot&\cdot&\cdot
\hln         &\cdot&\cdot&\cdot&\cdot&\cdot&\cdot&\cdot
\hln    \cdot&\cdot&\cdot&\cdot&\cdot&\cdot&\cdot&\cdot
\\\hline
\end{array}}
\ \stackrel{\textstyle \mu_3}{\longmapsto}\:
\begin{array}{@{}|@{\:}c@{\:}|@{\:}c@{\:}|@{\:}c@{\:}|@{\:}c@{\:}
  		 |@{\:}c@{\:}|@{\:}c@{\:}|@{\:}c@{\:}|@{\:}c@{\:}|@{}}
\hline    1  &     &  1  &     &     &  1  &  1  &\hspace{1pt}\cdot\hspace{1pt}
\hln      1  &     &  1  &  1  &  1  &     &\cdot&\cdot
\hln      2  &\mbf1&     &     &  1  &\cdot&\cdot&\cdot
\hln         &\mbf1&     &  2  &\cdot&\cdot&\cdot&\cdot
\hln         &  1  &  1  &\cdot&\cdot&\cdot&\cdot&\cdot
\hln         &     &\cdot&\cdot&\cdot&\cdot&\cdot&\cdot
\hln         &\cdot&\cdot&\cdot&\cdot&\cdot&\cdot&\cdot
\hln    \cdot&\cdot&\cdot&\cdot&\cdot&\cdot&\cdot&\cdot
\\\hline
\end{array}
\ \stackrel{\textstyle \mu_3}{\longmapsto}\:
\begin{array}{@{}|@{\:}c@{\:}|@{\:}c@{\:}|@{\:}c@{\:}|@{\:}c@{\:}
  		 |@{\:}c@{\:}|@{\:}c@{\:}|@{\:}c@{\:}|@{\:}c@{\:}|@{}}
\hline    1  &     &  1  &     &     &  1  &  1  &\hspace{1pt}\cdot\hspace{1pt}
\hln      1  &     &  1  &  1  &  1  &     &\cdot&\cdot
\hln      2  &\mbf2&     &     &  1  &\cdot&\cdot&\cdot
\hln         &     &     &  2  &\cdot&\cdot&\cdot&\cdot
\hln         &  1  &  1  &\cdot&\cdot&\cdot&\cdot&\cdot
\hln         &     &\cdot&\cdot&\cdot&\cdot&\cdot&\cdot
\hln         &\cdot&\cdot&\cdot&\cdot&\cdot&\cdot&\cdot
\hln    \cdot&\cdot&\cdot&\cdot&\cdot&\cdot&\cdot&\cdot
\\\hline
\end{array}
}\end{tinyrc}
$$
$$
\begin{tinyrc}{
\phantom{
\begin{array}{@{}cc@{}}{}\\\\3&\\4\\\\\\\\\\\end{array}
\begin{array}{@{}|@{\:}c@{\:}|@{\:}c@{\:}|@{\:}c@{\:}|@{\:}c@{\:}
  		 |@{\:}c@{\:}|@{\:}c@{\:}|@{\:}c@{\:}|@{\:}c@{\:}|@{}}
\hline    1  &     &  1  &     &     &  1  &  1  &\hspace{1pt}\cdot\hspace{1pt}
\hln      1  &     &  1  &  1  &  1  &     &\cdot&\cdot
\hln         &     &     &     &  1  &\cdot&\cdot&\cdot
\hln      2  &  2  &     &  2  &\cdot&\cdot&\cdot&\cdot
\hln         &  1  &  1  &\cdot&\cdot&\cdot&\cdot&\cdot
\hln         &     &\cdot&\cdot&\cdot&\cdot&\cdot&\cdot
\hln         &\cdot&\cdot&\cdot&\cdot&\cdot&\cdot&\cdot
\hln    \cdot&\cdot&\cdot&\cdot&\cdot&\cdot&\cdot&\cdot
\\\hline
\end{array}}
\ \stackrel{\textstyle \mu_3}{\longmapsto}\:
\begin{array}{@{}|@{\:}c@{\:}|@{\:}c@{\:}|@{\:}c@{\:}|@{\:}c@{\:}
  		 |@{\:}c@{\:}|@{\:}c@{\:}|@{\:}c@{\:}|@{\:}c@{\:}|@{}}
\hline    1  &     &  1  &     &     &  1  &  1  &\hspace{1pt}\cdot\hspace{1pt}
\hln      1  &     &  1  &  1  &  1  &     &\cdot&\cdot
\hln      2  &  2  &     &\mbf1&  1  &\cdot&\cdot&\cdot
\hln         &     &     &\mbf1&\cdot&\cdot&\cdot&\cdot
\hln         &  1  &  1  &\cdot&\cdot&\cdot&\cdot&\cdot
\hln         &     &\cdot&\cdot&\cdot&\cdot&\cdot&\cdot
\hln         &\cdot&\cdot&\cdot&\cdot&\cdot&\cdot&\cdot
\hln    \cdot&\cdot&\cdot&\cdot&\cdot&\cdot&\cdot&\cdot
\\\hline
\end{array}
\ \stackrel{\textstyle \mu_3}{\longmapsto}\:
\begin{array}{@{}|@{\:}c@{\:}|@{\:}c@{\:}|@{\:}c@{\:}|@{\:}c@{\:}
  		 |@{\:}c@{\:}|@{\:}c@{\:}|@{\:}c@{\:}|@{\:}c@{\:}|@{}}
\hline    1  &     &  1  &     &     &  1  &  1  &\hspace{1pt}\cdot\hspace{1pt}
\hln      1  &     &  1  &  1  &  1  &     &\cdot&\cdot
\hln      2  &  2  &     &\mbf2&  1  &\cdot&\cdot&\cdot
\hln         &     &     &     &\cdot&\cdot&\cdot&\cdot
\hln         &  1  &  1  &\cdot&\cdot&\cdot&\cdot&\cdot
\hln         &     &\cdot&\cdot&\cdot&\cdot&\cdot&\cdot
\hln         &\cdot&\cdot&\cdot&\cdot&\cdot&\cdot&\cdot
\hln    \cdot&\cdot&\cdot&\cdot&\cdot&\cdot&\cdot&\cdot
\\\hline
\end{array}
}\end{tinyrc}
$$
\caption{Mutation} \label{fig:mu}
\end{figure}
\end{example}

Having chosen our permutation $w$ and row index $i$, various entries of
a given $\bb \in \Znn$ play special roles.  To begin with, we call the
union of rows $i$ and $i+1$ the \bem{gene}%
        \footnote{All of the unusual terminology in what follows comes
        from genetics.  Superficically, our diagrams with two rows of
        boxes look like geneticists' schematic diagrams of the DNA
        double helix; but there is a much more apt analogy that will
        become clear only in Section~\ref{sec:coarsen}, where the
        biological meanings of the terms can be found in another
        footnote.}
of $\bb$.  For exponent arrays $\bb$ such that $\zz^\bb \not\in J_w$, the
spot in row $i$ and column
\begin{eqnarray} \label{eq:start}
  \start_i(\bb) &:=& \min\{p \mid z_{ip}\zz^\bb \not\in J_w\}
\end{eqnarray}
is called the \bem{start codon} of $\bb$.  The minimum defining
$\start_i(\bb)$ is taken over a nonempty set because $J_w \subseteq
J_{w_0} = I_{w_0} = \<z_{qp} \mid q+p \leq n\>$, so that $z_{in}\zz^\bb$
remains outside of $J_w$.  Of course, $z_{ip}\zz^\bb \not\in J_w$
whenever $z_{ip}$ divides $\zz^\bb \not\in J_w$ because $J_w$ is
generated by squarefree monomials.  Therefore,
\begin{eqnarray} \label{eq:leq}
  \start_i(\bb) &\leq& \west_i(\bb).
\end{eqnarray}
For completeness, set $\start_i(\bb)=0$ if $\zz^\bb \in J_w$.

Also of special importance is the \bem{promoter} $\prom(\bb)$, consisting
of the rectangular \mbox{$2 \times (\start_i(\bb)-1)$} array of locations
in the gene of $\bb$ that are strictly west of $\start_i(\bb)$.  Again,
we omit the explicit reference to $i$ and $w$ in the notation because
these are fixed for the discussion.  The sum of all entries in the
promoter of $\bb$ is
\begin{eqnarray} \label{eq:loose}
  |\prom(\bb)| &=& \sum_{j < \start_i(\bb)} b_{i+1,j}.
\end{eqnarray}

\begin{example} \label{ex:prom}
Let $\bb$ be the left array in Fig.~\ref{fig:mu}, $i = 3$, and $w =
13865742$, the permutation displayed in Example~\ref{ex:intro}.  Then
the gene of $\bb$ consists of rows $i = 3$ and $i+1 = 4$, and we claim
$\start_i(\bb) = 5$.

To begin with, we have $6$ choices for an antidiagonal $a \in J_w$
dividing $z_{31}\zz^\bb$: we must have $z_{31} \in a$, but other than
that we are free to choose one element of $\{z_{23}, z_{24}, z_{25}\}$
and one element of $\{z_{16},z_{17}\}$.  (This gives an example of the
$a$ produced in the first paragraph of the proof of
Lemma~\ref{lemma:outside}, below.)  Even more varied choices are
available for $z_{32}\zz^\bb$, such as $z_{41}z_{32}z_{23}$ or
$z_{41}z_{32}z_{13}$.  We can similarly find lots of antidiagonals in
$J_w$ dividing $z_{33}\zz^\bb$, and $z_{34}\zz^\bb$.  On the other hand,
$z_{35}$ already divides $\zz^\bb$, and one can verify that $\zz^\bb$ is
not divisible by the antidiagonals of any of the $2 \times 2$ or $3
\times 3$ minors defining $I_w$ (see Example~\ref{ex:intro}).  Therefore
$z_{35}\zz^\bb \not\in J_w$, so $\start_i(\bb) = 5$.

The promoter $\prom(\bb)$ consists of the $2 \times 4$ block
\begin{rcgraph}
  \begin{array}{l|c|c|c|c|}
          \cline{2-5}  3\  &     &     &     &
        \\\cline{2-5}  4\  &  2  &  2  &\ \, &  2
        \\\cline{2-5}
  \end{array}
\end{rcgraph}
at the western end.  In particular, $|\prom(\bb)| = 6$.

Nothing in this example depends on the values chosen for the dots on or
below the main antidiagonal.%
\end{example}

\section{Lifting Demazure operators}\label{sec:lift}

Now we need to understand the Hilbert series of $\cj w$ for varying $w$.
Since $J_w$ is a monomial ideal, its $\Znn$-graded Hilbert series $H(\cj
w; \zz)$ is simply the sum of all monomials outside~$J_w$.  Using the
combinatorics of the previous section, we construct operators $\ddem iw$
defined on monomials and taking the power series $H(\cj w; \zz)$ to
$H(\cj{ws_i}; \zz)$ whenever $\length(ws_i) < \length(w)$.  In other
words, the sum of all monomial outside $J_{ws_i}$ is obtained from the
sum of monomials outside $J_w$ by replacing $\zz^\bb \not\in J_w$ with
$\ddem iw(\zz^\bb)$.  It is worth keeping in mind that we shall
eventually show (in Section~\ref{sec:coarsen}) how $\ddem iw$ refines
the usual $\ZZ^n$-graded Demazure operator $\dem i$, when these
operators are applied to the variously graded Hilbert series of $\cj w$.

Again, fix for the duration of this section a permutation $w$ and a
transposition $s_i$ satisfying $\length(ws_i) < \length(w)$ .

\begin{defn} \label{defn:ddem}
The \bem{lifted Demazure operator} corresponding to $w$ and $i$ is a map
of abelian groups $\ddem iw : \ZZ[[\zz]] \too \ZZ[[\zz]]$ determined by
its action on monomials:
\begin{eqnarray*}
  \ddem iw(\zz^\bb) &:=& \sum_{d=0}^{|\prom(\bb)|} {\mu_i}^d(\zz^\bb).
\end{eqnarray*}
Here, ${\mu_i}^d$ means take the result of applying $\mu_i$ a total of
$d$ times, and ${\mu_i}^0(\bb) = \bb$.
\end{defn}

\begin{example} \label{ex:ddem}
If $\bb$ is the array in Examples~\ref{ex:mutate} and~\ref{ex:prom},
then $\ddem 3w(\zz^\bb)$ is the sum of the $7$ monomials whose exponent
arrays are displayed in Fig.~\ref{fig:mu}.%
\end{example}

Observe that $\ddem iw$ replaces each monomial by a homogeneous
polynomial of the same total degree, so the result of applying $\ddem
iw$ to a power series is actually a power~series.

In preparation for Theorem~\ref{thm:ev}, we need a few lemmas detailing
the effects of mutation on monomials and their genes.  The first of
these implies that $\ddem iw$ takes monomials outside $J_w$ to sums of
monomials outside $J_{ws_i}$, given that ${\mu_i}^0(\bb) = \bb$.

\begin{lemma} \label{lemma:outside}
If $\zz^\bb \not\in J_w$ and $1 \leq d \leq |\prom(\bb)|$ then
${\mu_i}^d(\zz^\bb) \in J_w \minus J_{ws_i}$.
\end{lemma}
\begin{proof}
We may as well assume $|\prom(\bb)| \geq 1$, or else the statement is
vacuous.  By definition of $\prom(\bb)$ and $\start_i(\bb)$, some
antidiagonal $a \in J_w$ divides $z_{ip}\zz^\bb$, where here (and for
the remainder of this proof) $p = \west_{i+1}(\bb)$.  Since $a$ doesn't
divide $\zz^\bb$, we find that $z_{ip} \in a$, whence $a$ cannot
intersect row $i+1$, which is zero to the west of~$z_{ip}$.  Thus $a$
also divides $\mu_i(\zz^\bb)$, and hence ${\mu_i}^d(\zz^\bb)$ for all
$d$ (including $d > |\prom(\bb)|$, but we won't need this).

It remains to show that ${\mu_i}^d(\zz^\bb) \not\in J_{ws_i}$ when $d
\leq |\prom(\bb)|$.  Let's start with $d \leq b_{i+1,p}$.  Any
antidiagonal $a$ dividing ${\mu_i}^d(\zz^\bb)$ does not continue
southwest of~$z_{ip}$; this is by Lemma~\ref{lemma:squish}(S) and the
fact that $\zz^\bb \not\in J_w$ (we could move $z_{ip}$ south).  Suppose
for contradiction that $a \in J_{ws_i}$, and consider the smallest
northwest submatrix $Z\sub i{j(a)}$ containing $a$.  If $j(a) \geq
w(i+1)$ then the antidiagonal $a' = a/z_{ip}$ obtained by omitting
$z_{ip}$ from $a$ is still in $J_{ws_i}$, being caused by the entry of
$\rk(ws_i)$ at $(i-1,j(a))$ as per Lemma~\ref{lemma:rank}.\ref{i-1}.  On
the other hand, if $j(a) < w(i+1)$, then $a'' = (z_{i+1,p}/z_{ip})a$ is
still in $J_{ws_i}$, being caused by the entry of $\rk(ws_i)$ at
$(i+1,j(a))$ as per Lemma~\ref{lemma:rank}.\ref{i+1}.  Since both $a'$
and $a''$ divide $\zz^\bb$ by construction, we find that $\zz^\bb \in
J_{ws_i} \subset J_w$, the desired contradiction.  It follows that
${\mu_i}^d(\zz^\bb) \not\in J_{ws_i}$ for $d \leq b_{i+1,p}$.

Assuming the result for $d \leq \sum_{j=p}^{p'} b_{i+1,j}$, where $p' <
\start_i(\bb) - 1$, we now demonstrate the result for $d \leq
\sum_{j=p}^{p'+1} b_{i+1,j}$.  Again, any antidiagonal $a \in J_w$
dividing ${\mu_i}^d(\zz^\bb)$ must end at row $i$, for the same reason
as in the previous paragraph.  But now if $a \in J_{ws_i}$, then moving
its southwest variable to $z_{ip}$ creates an antidiagonal that is in
$J_{ws_i}$ (by Lemma~\ref{lemma:squish}(W)) and divides
$\mu_i(\zz^\bb)$, which we have seen is impossible.%
\end{proof}

Now we show that mutation of monomials outside $J_w$ cannot produce the
same monomial more than once, as long we stop after $|\prom|$ many
steps.
\begin{lemma} \label{lemma:unique}
Suppose $\zz^\bb, \zz^{\bb'} \not\in J_w$ and that $d,d' \in \ZZ$
satisfy $1 \leq d \leq |\prom(\bb)|$ and $1 \leq d' \leq |\prom(\bb')|$.
If $\bb \neq \bb'$ then ${\mu_i}^d(\bb) \neq {\mu_i}^{d'}(\bb')$.
\end{lemma}
\begin{proof}
The inequality $d \leq |\prom(\bb)|$ guarantees that the mutations of
$\bb$ only alter the promoter of $\bb$, which is west of $\west_i(\bb)$
by~(\ref{eq:leq}).  Therefore, assuming (by switching $\bb$ and $\bb'$
if necessary) that $\west_i(\bb') \leq \west_i(\bb)$, we reduce to the
case where $\bb$ and $\bb'$ differ only in their genes, in columns
strictly west of $\west_i(\bb)$.

Let $\cc = {\mu_i}^d(\bb)$ and $\cc' = {\mu_i}^{d'}(\bb')$.  Mutating
preserves the sums
$$
  b_{i+1,j} = c_{ij} + c_{i+1,j} \quad {\rm and} \quad c'_{ij} +
  c'_{i+1,j} = b_{ij} + b_{i+1,j}
$$
for $j < \west_i(\bb)$, and we may as well assume these are equal for
every $j$, or else $\cc \neq \cc'$ is clear.  The westernmost column
where $\bb$ and $\bb'$ disagree is now necessarily $p = \west_i(\bb')$.
It follows that $z_{ip}\zz^\bb \not\in J_w$, because $\bb$ agrees with
$\bb'$ strictly to the west of column~$p$ as well as strictly to the
north of row $i$, and any antidiagonal $a \in J_w$ dividing
$z_{ip}\zz^\bb$ must be contained in this region (since it contains
$z_{ip}$).  In particular, $\start_i(\bb) \leq p$.  We conclude that
mutating $\bb$ and $\bb'$ fewer than $|\prom(\bb)|$ or $|\prom(\bb')|$
times cannot alter the column $p$ where $\bb$ and $\bb'$ differ.  Thus
$\cc$ differs from $\cc'$ in column$~p$.
\end{proof}

\begin{example}
If we apply $\mu_i$ more than $|\prom(\bb)|$ times to some array $\bb$,
it is possible to reach ${\mu_i}^{d'}(\bb')$ for some $\bb' \neq \bb$
and $d' \leq |\prom(\bb')|$.  Take $\bb$, $i$, and $w$ as in
Examples~\ref{ex:mutate}, \ref{ex:prom}, and~\ref{ex:ddem}, and set the
dot in $\bb$ at position $(4,5)$ equal to $3$.  If $\zz^{\bb'} =
(z_{35}/z_{45})\zz^\bb$, then we have $|\prom(\bb)| = |\prom(\bb')| =
6$, but the entries of $\bb$ and $\bb'$ in column $5$ of their genes are
$\aoverb 13$ and $\aoverb 22$, respectively.  Mutating $\bb$ and $\bb'$
up to $6$ times yields $7$ rrays each, all distinct because of the
$\aoverb 13$ and $\aoverb 22$ in column $5$.  However, mutating $\bb$ to
an $8^\th$ array $(\mu_3)^7(\bb)$ changes the $\aoverb 13$ to $\aoverb
22$, and outputs $(\mu_3)^6(\bb')$.
\end{example}

If $\cc$ is the result of applying $\mu_i$ to $\bb$ some number of
times, we can recover $\bb$ from $\cc$ by \bem{reverting} certain
entries of $\cc$ from row $i$ back to row $i+1$.  Formally, reverting an
entry $c_{ij}$ of $\cc$ means making a new array that agrees with $\cc$
except at $(i,j)$ and $(i+1,j)$.  In those spots, the new array has
$(i,j)$ entry $0$ and $(i+1,j)$ entry $c_{ij} + c_{i+1,j}$.  (In terms
of the stacks-of-coins picture, we revert only entire stacks of coins,
not single coins.)  Even if we are just given $\cc$ without knowing
$\bb$, we still have a criterion to determine when a certain reversion
of $\cc$ yields a monomial $\zz^\bb \not\in J_w$.

\begin{claim} \label{claim:min}
Suppose $\zz^\cc \in J_w \minus J_{ws_i}$.  If\/ $\bb$ is obtained from
$\cc$ by reverting all entries of $\cc$ in row $i$ that are west of or
at column $\west_{i+1}(\cc)$, then $\zz^\bb \not\in J_w$.
\end{claim}
\begin{proof}
Suppose $\zz^\bb \in J_w$, and let us try to produce an antidiagonal
witness $a \in J_w$ dividing~it.  Either $a$ ends at row $i$, or not.
In the first case, $a$ divides $\zz^\cc$, because the nonzero entries in
row $i$ of $\bb$ are the same as the corresponding entries of $\cc$.
Thus we can replace $a$ by the result $a'$ of tacking on $z_{i+1,p}$ to
$a$, where $p = \west_{i+1}(\cc)$.  This new $a'$ is in $J_w$
because $a \in J_w$ divides $a'$.  It follows from Lemma~\ref{lemma:i}
that $a' \in J_{ws_i}$.  Furthermore, $a'$ divides $\zz^\cc$ by
construction, and thus contradicts our assumption that $\zz^\cc \not\in
J_{ws_i}$.  Therefore, we may assume for the remainder of the proof of
this lemma that $a$ does not end at row $i$, so $a \in J_{ws_i}$ by
Lemma~\ref{lemma:i}.

We now prove that $\zz^{\bb} \not\in J_w$ by showing that if $a \in
J_{ws_i}$ and $a$ divides $\zz^{\bb}$, then from $a$ we can synthesize
$a' \in J_{ws_i}$ dividing $\zz^\cc \in J_{ws_i}$, again contradicting
our running assumption $\zz^\cc \in J_w \minus J_{ws_i}$.  There are
three possibilities (an illustration for~(\ref{ii}) is described in
Example~\ref{ex:ii}):
\begin{romanlist}
\item \label{i} The antidiagonal $a \in J_{ws_i}$ intersects row $i$ but
does not end there.

\item \label{ii} The antidiagonal $a \in J_{ws_i}$ skips row $i$ but
intersects row $i+1$.

\item \label{iii} The antidiagonal $a \in J_{ws_i}$ skips both row $i$
as well as row $i+1$.
\end{romanlist}
In case~(\ref{i}) either $a$ already divides $\zz^\cc$ or we can move
east the row $i+1$ variable in $a$, into the location
$(i+1,\west_{i+1}(\cc))$.  The resulting antidiagonal $a'$ divides $\cc$
by construction and is in $J_{ws_i}$ by Lemma~\ref{lemma:squish}(E).  In
case~(\ref{iii}) the antidiagonal already divides $\zz^\cc$ because
$\bb$ agrees with $\cc$ outside of their genes.

This leaves case~(\ref{ii}).  If $a$ does not already divide $\zz^\cc$,
then the intersection $z_{i+1,j}$ of $a$ with row $i+1$ is strictly west
of $\west_{i+1}(\cc)$.  The antidiagonal $a' = (z_{ij}/z_{i+1,j})a$ then
divides $\zz^\cc$ by construction, and is in $J_{ws_i}$ by
Lemma~\ref{lemma:squish}(N).
\end{proof}

\begin{example} \label{ex:ii}
Here is an instance of what occurs in case~(\ref{ii}) from the proof of
Claim~\ref{claim:min}.  Let $\bb$ and $\cc$ be the first and last arrays
from Example~\ref{ex:mutate}, and consider what happens when we fiddle
with their $(5,1)$ entries.  The antidiagonals $z_{51}z_{42}z_{23}$ and
$z_{51}z_{42}z_{24} \in J_{ws_i}$ both divide $z_{51}\zz^\bb$.  Using
Lemma~\ref{lemma:squish}(N) we can move the $z_{42}$ north to $z_{32}$
to get $a' \in \{z_{51}z_{32}z_{23}, z_{51}z_{32}z_{24}\}$ in $J_{ws_i}$
dividing $z_{51}\zz^\cc$.  (It almost goes without saying, of course,
that $z_{51}\zz^\cc$ is no longer in $J_w \minus J_{ws_i}$, so it does
not satisfy the hypothesis of Claim~\ref{claim:min}; we were, after all,
looking for a contradiction.)
\end{example}

Any array $\cc$ whose row $i$ begins west of its row $i+1$ can be
expressed as a mutation of {\em some}\/ array $\bb$.  By
Claim~\ref{claim:min}, we even know how to make sure $\zz^\bb \not\in
J_w$ whenever $\zz^\cc \in J_w \minus J_{ws_i}$.  But we also want each
$\zz^\cc \in J_w \minus J_{ws_i}$ to appear in $\ddem iw(\zz^\bb)$ for
some $\zz^\bb \not\in J_w$, and this involves making sure
$\start_i(\bb)$ is not too far west.

\begin{example} \label{ex:revert}
If $\west_i(\cc)$ is sufficiently smaller than $\west_{i+1}(\cc)$, then
it might be hard to determine which entries in row $i$ of $\cc$ to
revert while still assuring that $\zz^\cc$ appears in $\ddem
iw(\zz^\bb)$.  For example, let $\cc$ be the last array in
Example~\ref{ex:mutate}, that is, $\cc = (\mu_3)^6(\bb)$.  Suppose
further that the dot at $(4,5)$,
is really blank---i.e.\ zero.  Without a~priori knowing $\bb$, how are
we to know {\em not}\/ to revert the $1$ in position $(3,5)$?  Well,
suppose we did revert this entry, along with all of the entries west of
it in row $3$.  Then we would end up with an array $\bb'$ such that
$\zz^{\bb'} \not\in J_w$ all right, as per Claim~\ref{claim:min}, but
also such that $z_{35}\zz^{\bb'} \not\in J_w$.  This latter condition is
intolerable, since $5 \geq \start_i(\bb)$ implies that our original
$\zz^\cc$ will not end up in the sum $\ddem 3w(\zz^{\bb'})$.

Thus the problem with trying to revert the $1$ in position $(3,5)$ is
that it's too far east.  On the other hand, we might also try reverting
only the row $3$ entries in columns $1$ and~$2$, but with dire
consequences: we end up with an array $\bb''$ such that $\zz^{\bb''}$ is
divisible by $z_{42}z_{34}z_{25} \in J_w$.  (This is an example of the
antidiagonal $a$ to be produced after the displayed equation in the
proof of Lemma~\ref{lemma:revert}.)  We are left with only one choice:
revert the boldface~${\mathbf 2}$ in position $(3,4)$ and all of its
more western brethren.%
\end{example}

In general, as in the previous example, the \bem{critical column} for
$\zz^\cc \in J_w \minus J_{ws_i}$ is
\begin{eqnarray*}
  \crit(\cc) &:=& \min(p \leq \west_{i+1}(\cc)\mid z_{ip}\hbox{ divides }
  \zz^\cc \hbox{ and } z_{i+1,p}\zz^\cc \not\in J_{ws_i}).
\end{eqnarray*}

\begin{claim} \label{claim:crit}
If $\zz^\cc \in J_w \minus J_{ws_i}$ then:
\begin{numbered}
\item \label{crit1}
the set used to define $\crit(\cc)$ is nonempty;

\item \label{crit2}
reverting $c_{i,\crit(\cc)}$ creates an array $\cc'$ such that
$\zz^{\cc'} \not\in J_{ws_i}$; and

\item \label{crit3}
if $\west_i(\cc) < \crit(\cc)$, then the monomial $\zz^{\cc'}$ from
statement~\ref{crit2} remains in $J_w$.
\end{numbered}
\end{claim}
\begin{proof}
Claim~\ref{claim:min} implies $\west_i(\cc) \leq \west_{i+1}(\cc)$, so
$p' = \max(p \leq \west_{i+1}(\cc) \mid c_{ip} \neq 0)$ is well-defined.
If $a$ is an antidiagonal dividing the monomial whose exponent array is
the result of reverting $c_{ip'}$, then $a$ divides either $\zz^\cc$ or
the monomial $\zz^\bb$ from Claim~\ref{claim:min}, and neither of these
is in $J_{ws_i}$.  Thus $a \not\in J_{ws_i}$ and statemen~\ref{crit1} is
proved.  Part~\ref{crit2} is by definition, and statement~\ref{crit3}
follows from it by Lemmas~\ref{lemma:i} and~\ref{lemma:squish}(W).
\end{proof}

\begin{lemma} \label{lemma:revert}
Suppose $\zz^\cc \in J_w \minus J_{ws_i}$ and that~$\bb$ is obtained
by reverting all row $i$ entries of $\cc$ west of or at $\crit(\cc)$.
Then $\zz^\bb \not\in J_w$, and $\crit(\cc) < \start_i(\bb)$.
\end{lemma}
\begin{proof}
The proof that this $\zz^\bb$ is not in $J_w$ has two cases.  In the
first case we have $\crit(\cc) = \west_{i+1}(\cc)$, and
Claim~\ref{claim:min} immediately implies the result.  In the second
case we have $\crit(\cc) < \west_{i+1}(\cc)$, and we can apply
Claim~\ref{claim:min} to the monomial $\zz^{\cc'} \in J_w \minus
J_{ws_i}$ from Claim~\ref{claim:crit}.

Now we need to show $z_{ip}\zz^\bb \in J_w$ for two kinds of $p$: for $p
\leq \west_i(\cc)$ and $\west_i(\cc) < p \leq \crit(\cc)$.  (Of course,
when $\west_i(\cc) = \crit(\cc)$ the second of these cases is vacuous.)
The case $p \leq \west_i(\cc)$ is a little easier, so we treat it first.

There is some antidiagonal in $J_w$ ending on row $i$ and dividing
$\zz^\cc$, by Lemma~\ref{lemma:i}.  When $p \leq \west_i(\cc)$, we get
the desired result by appealing to Lemma~\ref{lemma:squish}(W).

Next we treat $\west_i(\cc) < p \leq \crit(\cc)$.  These inequalities
mean precisely that
\begin{eqnarray*}
  j &=& \max\{p' < p \mid c_{ip'} \neq 0\}
\end{eqnarray*}
is well-defined, and that $z_{i+1,j}\zz^\cc \in J_{ws_i}$.  Any
antidiagonal $a \in J_{ws_i}$ dividing $z_{i+1,j}\zz^\cc$ must contain
$z_{i+1,j}$ because $a$ does not divide $\zz^\cc$, and the fact that
$\zz^\bb \not\in J_w$ implies that $a$ also does not divide $\zz^\bb$.
It follows that $a$ intersects row $i$ at some spot in which $\cc$ is
nonzero strictly east of column $j$.  This spot is necessarily east of
or at $(i,p)$ by construction.  Without changing whether $a \in
J_{ws_i}$, Lemma~\ref{lemma:squish}(W) says that we may assume $a$
contains $z_{ip}$ itself.  This $a$ divides $z_{ip}\zz^\bb$, whence
$z_{ip}\zz^\bb \in J_{ws_i}$.
\end{proof}

The next theorem is the main result of Section~\ref{sec:lift}.  It
pinpoints, at the level of individual standard monomials, the relation
between $J_w$ and~$J_{ws_i}$.
\begin{thm} \label{thm:ev}
$_{\!}H(\cj{ws_i}; \zz) = \ddem iw H(\cj w; \zz)$ if\/ $\length(ws_i) <
\length(w)$.
\end{thm}
\begin{proof}
We need the sum $H(\cj{ws_i}; \zz)$ of monomials outside $J_{ws_i}$ to be
obtained from the sum of monomials outside $J_w$ by replacing $\zz^\bb
\not\in J_w$ with $\ddem iw(\zz^\bb)$.  We know by
Lemma~\ref{lemma:outside} that $\ddem iw H(\cj w)$ is a sum of monomials
outside $J_{ws_i}$.  Furthermore, no monomial $\zz^\cc$ is repeated in
this sum: if $\zz^\cc \not\in J_w$ appears in $\ddem iw(\bb)$, then $\bb$
must equal $\cc = {\mu_i}^0(\bb)$ by Lemma~\ref{lemma:outside}; and if
$\zz^\cc \in J_w$ then Lemma~\ref{lemma:unique} applies.

It remains to demonstrate that each monomial $\zz^\cc \not\in J_{ws_i}$
is equal to ${\mu_i}^d(\zz^\bb)$ for some monomial $\zz^\bb \not\in J_w$
and $d \leq |\prom(\bb)|$.  This is easy if $\zz^\cc$ is not even in
$J_w$: we take $\zz^\bb = {\mu_i}^0(\zz^\bb) = \zz^\cc$.  Since we can
now assume $\zz^\cc \in J_w \minus J_{ws_i}$, the result follows from
Lemma~\ref{lemma:revert}, once we notice that the inequality $\crit(\cc)
< \start_i(\bb)$ there is equivalent to the inequality $d \leq
|\prom(\bb)|$.
\end{proof}

\smallskip

\section{Coarsening the grading}\label{sec:coarsen}

As in Section~\ref{sec:lift}, fix a permutation $w$ and an index $i$
such that $\length(ws_i) < \length(w)$.  Our goal in this section is to
prove (in Theorem~\ref{thm:induction}) that the set of $\ZZ^n$-graded
Hilbert series $H(\cj w; \xx)$ for varying $w$ is closed under Demazure
operators.  The idea is to combine lifted Demazure operators $\ddem iw$
with the specialization $\XX : \ZZ[[\zz]] \to \ZZ[[\xx]]$ sending
$z_{qp} \mapsto x_q$; see Example~\ref{ex:coarsen}.  We present
Lemma~\ref{lemma:involution} and the proof
of~\mbox{Proposition~\ref{prop:lift}} in ``single'' language, for ease
of notation, but indicate at the end of the proof of
Proposition~\ref{prop:lift} which changes of notation make the arguments
work for the $\ZZ^{2n}$-graded Hilbert series $H(\cj w;\xx,\yy)$, with
the specialization $\XX_\YY : \ZZ[[\zz]] \to \ZZ[[\xx,\yy^{-\1}]]$
sending $z_{qp} \mapsto x_q/y_p$.

At the outset, we could hope that $\XX \circ \ddem iw = \dem i \circ
\XX$ monomial by monomial.  However, although this works in some cases
(see~(\ref{eq:fixed}), below) {\em it fails in general}.  The next
lemma will be used to take care of the general case.  Its proof is
somewhat involved (but fun) and irrelevant to its application, so we
postpone the proof until after
Theorem~\ref{thm:induction}.  Denote by $\std(J_w)$ the set of
\bem{standard exponent arrays\/}: the exponents on monomials not in
$J_w$.

\begin{lemma} \label{lemma:involution}
There is an involution $\tau : \std(J_w) \to \std(J_w)$ such that
$\tau^2 = 1$ and:
\begin{numbered}
\item
$\tau\bb$ agrees with $\bb$ outside their genes;

\item
$\prom(\tau\bb) = \prom(\bb)$;

\item \label{item:3}
if $\XX(\zz^\bb) = x_{i+1}^\ell\xx^\aa$ with $\ell = |\prom(\bb)|$, then
$\XX(\zz^{\tau\bb}) = x_{i+1}^\ell s_i(\xx^\aa)$; and

\item \label{item:4}
$\tau$ preserves column sums.  In other words, if\/ $\bb' = \tau\bb$,
then $\sum_q b_{qp} = \sum_q b'{}_{\!qp}$ for any fixed column
index~$p$.
\end{numbered}
In particular, $\XX(\ddem iw\zz^{\tau\bb}) = \dem i(x_{i+1}^\ell)
(s_i\xx^\aa)$.
\end{lemma}

\begin{remark} \label{rk:schur}
The squarefree monomials outside $J_w$ for a {Grassmannian
permutation}~$v$ (that is, a permutation having a unique descent) are in
natural bijection with the semistandard Young tableaux of the
appropriate shape and content.  (This follows from
Definition~\ref{defn:rp} and Theorem~\ref{thm:rp}, below, along with the
bijection in \cite{KoganThesis} between reduced pipe dreams and
semistandard Young tableaux.)  Under this natural bijection, intron
mutation reduces to an operation that arises in a well-known
combinatorial proof of the symmetry of the Schur function associated
to~$v$, which equals~$\SS_v$.
\end{remark}

Our next result justifies the term `lifted Demazure operator' for $\ddem
iw$.

\begin{prop} \label{prop:lift}
Specializing $\zz$ to $\xx$ in $\ddem iw H(\cj w; \zz)$ yields $\dem i
H(\cj w; \xx)$.  More generally, specializing $z_{qp}$ to $x_q/y_p$ in
$\ddem iw H(\cj w; \zz)$ yields $\dem i H(\cj w; \xx,\yy)$.
\end{prop}
\begin{proof}
Suppose $\zz^\bb \not\in J_w$ specializes to $\XX(\zz^\bb) =
x_{i+1}^\ell\xx^\aa$, where $\ell = |\prom(\bb)|$.  The definition of
$\ddem iw\zz^\bb$ implies that
$$
\begin{array}{r@{\ \:=\ \:}l}
\XX(\ddem iw \zz^\bb)
        &\sum_{d=0}^\ell x_i^d x_{i+1}^{\ell-d}\xx^\aa
\\[5pt] &\frac{x_{i+1}^{\ell+1} - x_i^{\ell+1}}{x_{i+1}-x_i}\xx^\aa
\\[5pt] &\dem i(x_{i+1}^\ell)\xx^\aa.
\end{array}
$$
If it happens that $s_i\xx^\aa = \xx^\aa$, so $\xx^\aa$ is symmetric in
$x_i$ and $x_{i+1}$, then
\begin{equation} \label{eq:fixed}
  \XX(\ddem iw\zz^\bb)\ \:=\ \:\dem i(x_{i+1}^\ell) \xx^\aa\ \:=\ \:
  \dem i(x_{i+1}^\ell \xx^\aa)\ \:=\ \:\dem i\XX(\zz^\bb).
\end{equation}
Of course, there will in general be lots of $\zz^\bb \not\in J_w$ whose
$\xx^\aa$ is not fixed by~$s_i$.  We overcome this difficulty using
Lemma~\ref{lemma:involution}, which says how to pair each $\zz^\bb
\not\in J_w$ with a partner so that their corresponding $\XX \circ \ddem
iw$ sums add up nicely.  Using the notation of the Lemma, notice that if
$\tau\bb = \bb$, then $s_i\xx^\aa = \xx^\aa$ and $\XX(\ddem iw\zz^\bb) =
\dem i\XX(\zz^\bb)$, as in~(\ref{eq:fixed}).  On the other hand, if
$\tau\bb \neq \bb$, then the Lemma~implies
$$
\begin{array}{r@{\ \:=\ \:}l}
\XX(\ddem iw(\zz^\bb + \zz^{\tau\bb}))
        &\dem i(x_{i+1}^\ell)(\xx^\aa + s_i\xx^\aa)
\\[5pt] &\dem i(x_{i+1}^\ell(\xx^\aa + s_i\xx^\aa))
\\[5pt] &\ddem iw(\XX(\zz^\bb + \zz^{\tau\bb}))
\end{array}
$$
because $\xx^\aa + s_i\xx^\aa$ is symmetric in $x_i$ and $x_{i+1}$.
This proves the $\ZZ^n$-graded statement.

The $\ZZ^{2n}$-graded version of the argument works mutatis mutandis by
the preservation of column sums under mutation (Definition~\ref{defn:mu}
and statement~\ref{item:4} of Lemma~\ref{lemma:involution}), which
allows us to replace $\XX$ by $\XX_\YY$ and $\xx^\aa$ by a monomial in
the $\xx$ variables and the inverses of the $\yy$ variables.
\end{proof}

Now we come to a result that will be crucial in proving the main
theorems of Part~\ref{part:intro}.

\begin{thm} \label{thm:induction}
$H(\cj{ws_i}; \xx) = \dem i H(\cj w; \xx)$ if\/ $\length(ws_i) <
\length(w)$.  More generally, $H(\cj{ws_i}; \xx,\yy) = \dem i H(\cj w;
\xx,\yy)$ if\/ $\length(ws_i) < \length(w)$.
\end{thm}
\begin{proof}
Theorem~\ref{thm:ev} and Proposition~\ref{prop:lift}.
\end{proof}

Before constructing this magic involution~$\tau$, we introduce some
necessary notation and provide examples.  Recall that the union of rows
$i$ and $i+1$ is the \textbf{gene} of $\bb$ (we view the row index~$i$
as being fixed for the discussion).  Order the boxes in columns east of
$\start_i(\bb)$ in the gene of $\bb$ as in the diagram, using the
notation $\start_i(\bb)$ from~(\ref{eq:start})
in Section~\ref{sec:mutation}:
$$
\begin{array}{l|c|c|c|c|c|c}
\multicolumn{1}{l}{\mbox{}}\\[-4ex]
\multicolumn{2}{l}{\mbox{}}&
  \multicolumn{1}{c}{
        \begin{array}{@{}c@{}}
                \makebox[0pt]{$\scriptstyle \start_i(\bb)$}\\
                \downarrow
        \end{array}}
                                        \\[.2ex]\cline{2-7}
\petit{i}  &\toplinedots&  1  &  3  &  5  &  7  &\toplinedots\\\cline{3-6}
\petit{i+1}&            &  2  &  4  &  6  &  8  &            \\\cline{2-7}
\end{array}
$$
Now define five different kinds of blocks in the gene of $\bb$, called
the promoter, the start codon, exons, introns, and the stop codon.%
        \footnote{All of these are terms from genetics.  The DNA
        sequence for a single gene is not necessarily contiguous.
        Instead, it sometimes comes in blocks called \bem{exons}.  The
        intervening DNA sequences whose data are excised are called
        \bem{introns} (note that the structure of the gene of an
        exponent array is determined by its exons, not its introns).
        The \bem{promoter} is a medium-length region somewhat before the
        gene that signals the transcriptase enzyme where to attach to
        the DNA, so that it may begin transcribing the DNA into RNA.
        The \bem{start codon} is a short sequence signaling the
        beginning of the actual gene; the \bem{stop codon} is a similar
        sequence signaling the end of the gene.}
In the following, $k,\ell \in \NN$.
\begin{itemize}
\item
\bem{promoter}: the rectangle consisting of unnumbered boxes at the left
end

\item
\bem{start codon}: the box numbered~$1$, which lies at
$(i,\start_i(\bb))$

\item
\bem{stop codon}: the last numbered box, which lies at $(i+1,n)$

\item
\bem{exon}: any sequence $2k,\ldots,2\ell+1$ (with $k \leq \ell$) of
consecutive boxes satisfying:
\begin{numbered}
\item
the entries of~$\bb$ in the boxes corresponding to $2k+1,\ldots,2\ell$
are all zero;

\item
either box $2k+1$ is the start codon, or box~$2k$ has a nonzero entry
in~$\bb$;~and

\item
either box $2\ell$ is the stop codon, or box $2\ell+1$ has a nonzero
entry in~$\bb$
\end{numbered}

\item
\bem{intron}: any rectangle of consecutive boxes $2\ell+1,\ldots,2k$
(with $\ell < k$) satisfying:
\begin{numbered}
\item
the rectangle contains no exons;

\item
box $2\ell+1$ is either the start codon or the last box in an exon; and

\item
box $2k$ is either the stop codon or the fisrt box in an exon
\end{numbered}
\end{itemize}
Roughly speaking, the nonzero entries in gene$(\bb)$ are parititioned
into the promoter and introns, the latter being contiguous rectangles
having nonzero entries in their northwest and southeast corners.  Exons
connect adjacent introns via bridges of zeros.

\begin{example} \label{ex:gene}
Suppose we are given a permutation $w$, an array $\bb$ such that
$\zz^\bb \not\in J_w$, and a row index $i$ such that $\start_i(\bb) = 4$
and $\bb$ has the gene in Figure~\ref{fig:intron}.
\begin{figure}[t]
\begin{rcgraph}
\begin{array}{rrc@{\ }c@{\hspace{-3.6ex}}c@{\hspace{-3.6ex}}c@{\hspace{-3.6ex}}
		    c@{\hspace{-3.6ex}}c@{\hspace{-3.6ex}}c@{\hspace{-3.6ex}}c}
&
{\rm promoter\colon}
&
\begin{array}{|c|c|c|}
  \hline           & \ \,&
  \\\hline      2  &     &  3
  \\\hline
\end{array}\!\!
\\
\\
{\rm gene\ of\ }\bb\colon\hspace{-5ex}
&
\multicolumn{9}{r}{
  \begin{array}{l|c|c|c|c|c|c|c|c|c|c|c|c|c|c|c|c|c|c|c|c|}
  \multicolumn{1}{l}{\mbox{}}\\[-4ex]
  \multicolumn{4}{l}{\mbox{}}&
    \multicolumn{1}{c}{
          \begin{array}{@{}c@{}}
                  \makebox[0pt]{start codon}\\
                  \downarrow
          \end{array}}
                                          \\\cline{2-21}
  i  &  &\ \,&  & 6&  &\ \,&  &  & 4&\ \,& 3& 8& 6&\ \,& 2&  &  &\ \,&\ \,& 5
  \\\cline{2-21}
  i+1& 2&    & 3& 1& 4&    & 5& 3& 7&    &  &  &  &    & 5& 1& 4&    &    &
  \\\cline{2-21}
    \multicolumn{20}{c}{}&
    \multicolumn{1}{c}{
          \begin{array}{@{}c@{}}
                  \uparrow\\
                  \makebox[0pt]{stop codon}
          \end{array}}
  \end{array}\:
}
\\
\\
&
{\rm exons\colon}
&&&
\raisebox{.15ex}{%
$
\begin{array}{|c|c|}
  \cline{2-2}
        \multicolumn{1}{c|}{} & 4
  \\\hline                 3
  \\\cline{1-1}
\end{array}
$}
&&
\raisebox{-.15ex}{%
$
\begin{array}{|c|c|c|}
  \cline{2-3}
        \multicolumn{1}{c|}{} &\ \,& 3 
  \\\hline                 7  &
  \\\cline{1-2}
\end{array}
$}
&&
\begin{array}{|c|c|c|c|}
  \cline{2-4}
        \multicolumn{1}{c|}{} &\ \,&    & 5
  \\\hline                 4  &    &\ \,
  \\\cline{1-3}
\end{array}
\\[3ex]
&
{\rm introns\colon}
&
&
\begin{array}{@{}|c|c|c|c|c|@{}}
  \hline        6  &     &    &     &
  \\\hline      1  &  4  &\ \,&  5  &  3  
  \\\hline
\end{array}
&&
\begin{array}{@{}|c|@{}}
  \hline        4  
  \\\hline      7  
  \\\hline
\end{array}
&&
\begin{array}{@{}|c|c|c|c|c|c|c|@{}}
  \hline        3  &  8  &  6  &\ \,&  2  &     &   
  \\\hline         &     &     &    &  5  &  1  &  4
  \\\hline
\end{array}
&&
\begin{array}{@{}|c|@{}}
  \hline        5  
  \\\hline         
  \\\hline
\end{array}
\end{array}
\end{rcgraph}
\medskip
For ease of comparison, we dissect here the mutated gene $\tau\bb$
of~$\bb$.\hfill\mbox{}
\medskip
\begin{rcgraph}
\begin{array}{rrc@{\ }c@{\hspace{-3.6ex}}c@{\hspace{-3.6ex}}c@{\hspace{-3.6ex}}
		    c@{\hspace{-3.6ex}}c@{\hspace{-3.6ex}}c@{\hspace{-3.6ex}}c}
&
{\rm promoter\colon}
&
\begin{array}{|c|c|c|}
  \hline           & \ \,&
  \\\hline      2  &     &  3
  \\\hline
\end{array}\!\!
\\
\\
{\rm gene\ of\ }\tau\bb\colon\hspace{-5ex}
&
\multicolumn{9}{r}{
  \begin{array}{l|c|c|c|c|c|c|c|c|c|c|c|c|c|c|c|c|c|c|c|c|}
  \multicolumn{1}{l}{\mbox{}}\\[-4ex]
  \multicolumn{4}{l}{\mbox{}}&
    \multicolumn{1}{c}{
          \begin{array}{@{}c@{}}
                  \makebox[0pt]{start codon}\\
                  \downarrow
          \end{array}}
                                          \\\cline{2-21}
  i  &  &\ \,&  & 7& 4&\ \,& 1&  & 7&\ \,& 3& 7&  &\ \,&  &  &  &\ \,&\ \,& 1
  \\\cline{2-21}
  i+1& 2&    & 3&  &  &    & 4& 3& 4&    &  & 1& 6&    & 7& 1& 4&    &    & 4
  \\\cline{2-21}
    \multicolumn{20}{c}{}&
    \multicolumn{1}{c}{
          \begin{array}{@{}c@{}}
                  \uparrow\\
                  \makebox[0pt]{stop codon}
          \end{array}}
  \end{array}\:
}
\\
\\
&
{\rm exons\colon}
&&&
\raisebox{.15ex}{%
$
\begin{array}{|c|c|}
  \cline{2-2}
        \multicolumn{1}{c|}{} & 7
  \\\hline                 3
  \\\cline{1-1}
\end{array}
$}
&&
\raisebox{-.15ex}{%
$
\begin{array}{|c|c|c|}
  \cline{2-3}
        \multicolumn{1}{c|}{} &\ \,& 3 
  \\\hline                 4  &
  \\\cline{1-2}
\end{array}
$}
&&
\begin{array}{|c|c|c|c|}
  \cline{2-4}
        \multicolumn{1}{c|}{} &\ \,&    & 1
  \\\hline                 4  &    &\ \,
  \\\cline{1-3}
\end{array}
\\[3ex]
&
{\rm introns\colon}
&
&
\begin{array}{|c|c|c|c|c|}
  \hline        7  &  4  &    &  1  &
  \\\hline         &     &\ \,&  4  &  3  
  \\\hline
\end{array}
&&
\begin{array}{|c|}
  \hline        7  
  \\\hline      4  
  \\\hline
\end{array}
&&
\begin{array}{|c|c|c|c|c|c|c|}
  \hline        3  &  7  &     &\ \,&     &     &   
  \\\hline         &  1  &  6  &    &  7  &  1  &  4
  \\\hline
\end{array}
&&
\begin{array}{|c|}
  \hline        1  
  \\\hline      4  
  \\\hline
\end{array}
\end{array}
\end{rcgraph}
\caption{Intron mutation}\label{fig:intron}
\end{figure}
The gene of $\bb$ breaks up into promoter, start codon, exons, introns,
and stop codon as indicated.  We shall say something more about the
mutated gene $\tau\bb$ in Example~\ref{ex:tau}.
\end{example}

If $\cc$ is an array having two rows filled with nonnegative
integers, then let
$\ol\cc$ be the rectangle obtained by rotating $\cc$ through an angle of
$180^\circ$.  For purposes of applying the mutation operator~$\mu_i$
(Definition~\ref{defn:mu}), we identify the rows of an intron $\cc$ as
rows $i$ and $i+1$ in a gene, and we view $\cc$ as an $n \times n$ array
that happens to be zero outside of its $2 \times k$ rectangle.

\begin{defn}[Intron mutation] \label{defn:mutation}
Let $c_i$ and $c_{i+1}$ be the sums of the entries in the top and bottom
nonzero rows of an intron~$\cc$, and set $d = |c_i-c_{i+1}|$.  Then
\begin{eqnarray*}
\tau\cc &=& \left\{\begin{array}{@{}l@{\:\ \textrm{if}\ }l}
			\ol{{\mu_i}^d(\ol\cc)} & c_i > c_{i+1} \\
			{\mu_i}^d(\cc) & c_i < c_{i+1}
		   \end{array}\right.
\end{eqnarray*}
is the \bem{mutation} of~$\cc$.  Define the \bem{intron mutation}
$\tau\bb$ of an exponent array $\bb$
by
\begin{itemize}
\item
adding~$1$ to the start and stop codons of~$\bb$;

\item
mutating every intron in the gene of the resulting exponent array; and
then
\item
subtracting~$1$ from the boxes that were the start and stop codons
of~$\bb$.
\end{itemize}
\end{defn}

Intron mutation pushes the entries of each intron either upward from
left to right or downward from right to left---whichever initially
brings the row sums in that particular intron closer to agreement.

\begin{example} \label{ex:tau}
Although the ``look'' of $\bb$ in Example~\ref{ex:gene} completely
changes when it is mutated into $\tau\bb$, the columns of $\tau\bb$
containing a nonzero entry are exactly the same as those in $\bb$, and
the column sums are preserved.  Note that mutating the gene of $\tau\bb$
yields back the gene of $\bb$, {\em as long as $\zz^\cc \not\in J_w$ and
the location of the start codon has not changed}.  The proof of
Lemma~\ref{lemma:involution} shows why $\tau$ always works this~way.%
\end{example}

\begin{lemma} \label{lemma:exons}
Intron mutation outputs an exponent array (that is, the entries are
nonnegative).  Assume, for the purpose of defining exons in $\tau\bb$,
that the start codon of $\tau\bb$ lies at the same location as the start
codon in~$\bb$.  The boxes occupied by exons of\/ $\tau\bb$ thus defined
coincide with the boxes occupied by exons of\/ $\bb$ itself.
\end{lemma}
\begin{proof}
The definitions ensure that any intron not containing the start or stop
codon has nonzero northwest and southeast corners.  After adding $1$ to
the start and stop codons, every intron has this property.  Mutation of
such an intron leaves strictly positive entries in the northwest and
southeast corners (this is crucial---it explains why we have to add and
subtract the $1$'s from the codons), so subtracting~$1$ preserves
nonnegativity.  Furthermore, intron mutation does not introduce any new
exons, because the nonzero entries in an intron both before and after
mutation follow a snake pattern that drops from row~$i$ to row~$i+1$
precisely once.%
\end{proof}

\begin{proofof}{Lemma~\ref{lemma:involution}}
First we show that $\tau\bb \in \std(J_w)$, or equivalently that
$\mbox{$\zz^{\tau\bb} \in J_w$} \implies \zz^\bb \in J_w$.  Observe
that $\zz^\bb \in J_w$ if and only if $z_{ip}z_{i+1,n}\zz^\bb \in
J_w$, where $p = \start_i(\bb)$, by definition of $\start_i(\bb)$ and
the fact that $z_{i+1,n}$ is a nonzerodivisor modulo $J_w$ for all
$w$. Therefore, it suffices to demonstrate how an antidiagonal $a \in
J_w$ dividing $\zz^{\tau\bb}$ gives rise to a possibly different
antidiagonal $a' \in J_w$ dividing $z_{ip}z_{i+1,n}\zz^\bb$, where $p
= \start_i(\bb)$.  There are five cases:
\begin{romanlist}
\item $a$ intersects neither row $i$ nor row $i+1$;

\item the southwest variable in $a$ is in row $i$;

\item $a$ intersects row $i$ and continues south, but skips row $i+1$;

\item $a$ intersects both row $i$ and row $i+1$; or

\item $a$ skips row $i$ but intersects row $i+1$.
\end{romanlist}
In each of these cases, $a'$ is constructed as follows.  Outside the
gene of the generic matrix~$Z$, the new $a'$ will agree with $a$ in all
five cases, since $\bb$ and $\tau\bb$ agree outside of their genes.
Inside their genes, we may need some adjustments.
\begin{romanlist}
\item Leave $a'=a$ as is.

\item Move the variable in row $i$ west to $z_{ip}$, using
Lemma~\ref{lemma:squish}(W).

\item \label{item:iii}
The gene of $\bb$ has nonzero entries in precisely the same columns as
the gene of $\tau\bb$, by definition. Either $a$ already divides
$z_{ip}\zz^\bb$, or moving the variable in row $i$ due south to row
$i+1$ yields $a'$ by Lemma~\ref{lemma:squish}(S).

\item
Use Example~\ref{ex:squish} to make $a'$ contain the nonzero entries
in some exon of~$\bb$ (see Lemma~\ref{lemma:exons}).

\item
Same as~(\ref{item:iii}), except that either $a$ already divides
$z_{i+1,n}\zz^\bb$ or Lemma~\ref{lemma:squish}(N) says we can move the
variable due north from row $i+1$ to row $i$.
\end{romanlist}

Now that we know $\tau\bb \in \std(J_w)$, we find that
\begin{eqnarray*}
  \start_i(\tau\bb) &=& \start_i(\bb).
\end{eqnarray*}
Indeed, when $j \leq \start_i(\bb)$, we have $z_{ij}\zz^{\tau\bb} \in
J_w$ if and only if $z_{ij}\zz^\bb \in J_w$, because any antidiagonal
containing $z_{ij}$ interacts with $\tau\bb$ north of row $i$ and west
of column $\start_i(\bb)$, where $\tau\bb$ agrees with $\bb$.  It
follows that $\prom(\tau\bb) = \prom(\bb)$, so exons of $\tau\bb$ occupy
the same boxes as those of~$\bb$ by Lemma~\ref{lemma:exons}.
We conclude that $\tau\bb$ also has introns in the same boxes as the
introns of~$\bb$.  The statement $\tau^2 =$ identity holds because the
partitions of the genes of $\bb$ and~$\tau\bb$ into promoter and introns
are the same, and mutation on each of these blocks in the partition has
order $1$ or $2$.  Part~\ref{item:3} follows intron by intron, except
that the added and subtracted $1$'s in the first and last introns
cancel.%
\end{proofof}

\section{Equidimensionality}\label{sec:Lw}

We demonstrate here that the facets of $\LL_w$ all have the same
dimension.

The monomials $\zz^\bb$ that are nonzero in $\cj w$ are the so-called
\bem{standard monomials} for $J_w$, and are precisely those with support
sets
\begin{eqnarray*}
  \supp(\zz^\bb) &:=& \{z_{qp} \in Z \mid z_{qp} \hbox{ divides }
  \zz^\bb\}
\end{eqnarray*}
in the complex $\LL_w$.  The maximal support sets of standard monomials
are the facets of $\LL_w$.  The following lemma says that maximal
support monomials for $J_{ws_i}$ can only be mutations of maximal
support monomials for~$J_w$.

\begin{lemma} \label{lemma:notafacet}
If $\zz^\bb$ divides $\zz^\cc \not\in J_w$ and $d \leq |\prom(\bb)|$,
then ${\mu_i}^d(\zz^\bb)$ divides ${\mu_i}^e(\zz^\cc)$ for some $e \leq
|\prom(\cc)|$.  If $\supp(\zz^\bb) \in \LL_w$ is not a facet and $d \leq
|\prom(\bb)|$, then $\supp({\mu_i}^d(\zz^\bb)) \in \LL_{ws_i}$ is not a
facet.
\end{lemma}
\begin{proof}
If $\zz^\bb$ divides $\zz^\cc$, then $\start_i(\bb) \leq \start_i(\cc)$
by definition ($z_{ip}\zz^\bb \in J_w \implies z_{ip}\zz^\cc \in J_w$).
Therefore, if the $d^\th$ mutation of $\bb$ is the $k^\th$ occurring in
column $j < \start_i(\bb)$, we can choose $e$ so that the $e^\th$
mutation of $\cc$ is also the $k^\th$ occurring in column~$j$.  If, in
addition, $z_{qp} \in \supp(\zz^\cc) \minus \supp(\zz^\bb)$, then either
$z_{qp}$ or $z_{q-1,p}$ ends up in $\supp({\mu_i}^e(\zz^\cc)) \minus
\supp({\mu_i}^d(\zz^\bb))$, depending on whether or not $(q,p) \in
\prom(\cc)$ and $p < j$.  Note that ${\mu_i}^d(\zz^\bb)$ and
${\mu_i}^e(\zz^\cc)$ are not in $J_{ws_i}$ by Theorem~\ref{thm:ev}, so
their supports are in $\LL_{ws_i}$.%
\end{proof}

Recall that $D_0$ is the pipe dream with crosses in the strict
upper-left triangle, that is, all locations $(q,p)$ such that $q+p \leq
n$.  Number the crosses in $D_0$ as follows, where $N = \binom n2$:
\def\mcn#1{\multicolumn{#1}{@{}c@{}}{}}
\def\phaN{\phantom{8}}
\def\drop{\begin{array}{@{}c@{}}\vdots\\[-5ex]\end{array}}
\begin{rcgraph}
\begin{array}{|c|c|c|c|c|c|c|}
\hline
\!\!N\!\!&       &\!\!10\!\!& 6  &  3  &  1  &\ \:
\\\cline{1-1}\cline{3-6}
  \phaN &        &  9  &  5  &  2  &\mcn1&
\\\cline{3-5}
  \drop & \cdots &  8  &  4  &\mcn2&
\\\cline{3-4}
        &        &  7  &\mcn3&
\\\cline{3-3}
        &\multicolumn{1}{|c@{}}{\raisebox{1ex}{$\!\!\!\!\adots$}}&\mcn4&
\end{array}
\end{rcgraph}
For notation, given an $n \times n$ exponent array~$\bb$, let $D(\bb) =
[n]^2 \minus \supp(\zz^\bb)$.

\begin{lemma} \label{lemma:alpha}
Consider the following condition on a permutation~$w$: there is a fixed
$\alpha \geq 1$ such that for every facet $L$ of~$\LL_w$, the associated
pipe dream $D_L$ has crosses in boxes marked $\geq \alpha$ and an elbow
joint in the box due south of the box marked~$\alpha$.  Given this
condition, it follows that if $\alpha$ sits in row~$i$, then for every
facet $L'$ of~$\LL_{ws_i}$, the associated pipe dream $D_{L'}$ has
crosses in boxes marked $\geq \alpha + 1$ and an elbow joint at
$\alpha$.  Moreover, if $\alpha$ sits in column~$j$ and $\supp(\zz^\bb)$
is a facet of~$\LL_w$, then $\start_i(\bb) > j$.
\end{lemma}
\begin{proof}
If $\alpha$ sits in column~$j$, then every variable other than
$z_{i+1,j}$ in $Z\sub {i+1}j$ lies in~$J_w$, by the hypothesis
on~$\alpha$.  It follows that $w(i+1) = j$ and $\length(ws_i) <
\length(w)$, because $\rank(w^T\sub {i+1}j) = 1$ while $\rank(w^T\sub qp)
= 0$ whenever $z_{i+1,j} \neq z_{qp} \in Z\sub {i+1}j$.

By Lemma~\ref{lemma:notafacet}, every facet $L' \in \LL_{ws_i}$ can be
expressed as $\supp({\mu_i}^d(\zz^\bb))$ for some monomial $\zz^\bb$
such that $\supp(\zz^\bb) \in \LL_w$ is a facet.  The maximality of
$\supp(\zz^\bb)$ implies $\start_i(\bb) = \west_i(\bb)$; but
$\west_i(\bb) > j$ because $D(\bb)$ has crosses in $Z\sub ij$.  The
result follows because all mutations of $\zz^\bb$ (except $\zz^\bb$
itself) are therefore divisible by~$z_{ij}$.%
\end{proof}

In the hypothesis of Lemma~\ref{lemma:alpha}, note that the box due
south of $\alpha$ will not be marked $\alpha-1$ when $\alpha-1$ lies in
the top row.  The next result allows induction on~$\alpha$.

\begin{lemma} \label{lemma:n!}
Given a permutation $v \in S_n$ with $v \neq w_0$, there exists an
integer $\alpha \geq 1$ such that $w$ satsifies the hypothesis of
Lemma~\ref{lemma:alpha} and $v = ws_i$.  In particular, there exists a
unique $\alpha \in \{1,\ldots,\binom n2\}$ such that for every facet $L$
of~$\LL_v$, the associated pipe dream $D_L$ has crosses in boxes marked
$\geq \alpha+1$ and an elbow joint at~$\alpha$.
\end{lemma}
\begin{proof}
For each $\alpha \geq 0$, let $\nu(\alpha)$ be the number of
permutations $ws_i$ satisfying the conclusion of Lemma~\ref{lemma:alpha}
for that choice of~$\alpha$, where we declare $\nu(0) = 1$ to take care
of~$ws_i = w_0$.  Lemma~\ref{lemma:alpha} implies by induction
on~$\alpha$ that
\begin{itemize}
\item
$\nu(\alpha) \geq \nu(\alpha-1)$ if $\alpha-1$ does not lie in the top
row; and

\item
$\nu(\alpha) \geq \sum_{\beta \leq \alpha} \nu(\beta)$ if $\alpha-1$
lies in the top row.
\end{itemize}
Let $\nu_j = \sum \nu(\beta)$ be the sum over all $\beta$ in column~$j$,
and set $\sigma_{j+1} = \sum_{j' \geq j+1} \nu_{j'}$.  The itemized
claims above imply that if $\alpha$ lies in column~$j$, then
$\nu(\alpha) \geq \sigma_{j+1}$.  Assume by downward induction on~$j$
that $\sigma_{j+1} \geq (n-j)!$.  Then $\nu_j \geq (n-j)(n-j)!$ because
there are $n-j$ marked boxes in column~$j$.  It follows that $\sigma_j
\geq (n-j)! + (n-j)(n-j)! = (n-j+1)!$.  In particular, $\sigma_1 \geq
n!$, whence $\sigma_1 = n!$.%
\end{proof}

If we knew a~priori that $J_w$ were the initial ideal of $I(\ol X_w)$,
then the following Proposition would follow from \cite{KSprimeIdeals}
(except for embedded components).  Since we are working in the opposite
direction, we have to prove it ourselves.

\begin{prop} \label{prop:pure}
The simplicial complex $\LL_w$ is pure, each facet having $\dim\ol X_w =
n^2 - \length(w)$ vertices; i.e.\ $\spec(\cj w)$ is equidimensional of
dimension $\dim\ol X_w$.
\end{prop}
\begin{proof}
The result is obvious for $\LL_{w_0}$.  Taking $w$ and~$ws_i$ as in
Lemma~\ref{lemma:alpha} by Lemma~\ref{lemma:n!}, we prove the result
for~$ws_i$ by assuming it for~$w$.  Theorem~\ref{thm:ev} implies that
the facets of $\LL_{ws_i}$ are supports of monomials
${\mu_i}^d(\zz^\bb)$ for $\zz^\bb \not\in J_w$ and $d \leq
|\prom(\bb)|$.  By Lemma~\ref{lemma:notafacet}, we may restrict our
attention to square monomials $\zz^{2\bb}$ with maximal support.

Repeated mutation increases the support size of a monomial by~$0$ or~$1$
over the original.  Hence it suffices to show that if the cardinalities
of $\supp({\mu_i}^d(\zz^{2\bb}))$ and $\supp(\zz^{2\bb})$ are equal for
some $d \leq |\prom(2\bb)|$, then $\supp({\mu_i}^d(\zz^{2\bb}))
\subsetneq \supp({\mu_i}^e(\zz^{2\bb}))$ for some $e \leq
|\prom(2\bb)|$.  By Lemma~\ref{lemma:alpha} we may take $e = 1$ if $d =
0$.  If $d > 0$, then mutation will have just barely pushed a row $i+1$
entry of $2\bb$ up to row~$i$, and we may take $e = d-1$.  Containment
uses that every entry $b_{ij}$ is at least~$2$; strict containment is
automatic, by checking cardinalities.%
\end{proof}

\begin{example} \label{ex:434}
The ideal $J_{1432}$ is generated by the antidiagonals of the five $2
\times 2$ minors contained in the union of the northwest $2 \times 3$
and $3 \times 2$ submatrices of $(z_{ij})$:
\begin{eqnarray*}
J_{1432}
& = & \<z_{12}z_{21}, z_{13}z_{21}, z_{13}z_{22}, z_{12}z_{31},
        z_{22}z_{31}\>\\
& = & \<z_{12},z_{13},z_{22}\> \cap \<z_{12},z_{21},z_{22}\> \cap
        \<z_{21},z_{22},z_{31}\> \cap \<z_{13},z_{21},z_{31}\> \cap
        \<z_{12},z_{13},z_{31}\>.
\end{eqnarray*}
$\LL_{1432}$ is the join of a pentagon with a simplex having $11$
vertices $\{z_{11}\} \cup \{z_{rs} \mid r+s \geq 5\}$ (note $n=4$ here).
Each facet of $\LL_w$ therefore has $13 = 4^2 - \length(1432)$
vertices.%
\end{example}

\section{Mitosis on facets}\label{sec:facets}

This section and the next translate the combinatorics of lifted Demazure
operators into a language compatible with reduced pipe dreams.  The
present section concerns the relation between mutation and mitosis.  We
again think of the row index~$i$ as being fixed, as we did in
Sections~\ref{sec:mutation}--\ref{sec:coarsen}.

Recall the definition of pipe dream from Section~\ref{sec:pipe}.  The
similarity between $\start_i(D)$ for pipe dreams~$D$ introduced there
and $\start_i(\bb)$ for arrays in Section~\ref{sec:lift} is apparent.
It is explained precisely in the following lemma, whose proof is
immediate from the definitions.

\begin{lemma} \label{lemma:max}
If $\zz^\bb \not\in J_w$ has maximal support and $D(\bb) = [n]^2 \minus
\supp(\zz^\bb)$, then
$$
  \start_i(\bb)\:\ =\:\ \west_i(\bb)\:\ =\:\ \start_i(D(\bb)).
$$
\end{lemma}

\begin{excise}{
For precision, we record here (perhaps laboriously) in symbols the
definition of mitosis from Section~\ref{sec:alg}.
\begin{defn} \label{d:mitosis}
Given a pipe dream $D \subseteq [n]^2$ and a row index~$i$, construct a
set $\mitosis_i(D)$ of pipe dreams contained in $[n]^2$ as follows.  Set
\begin{eqnarray*}
\JJ(D) &=& \{j < \start_i(D) \mid (i+1,j) \not\in D\};\\
\JJ_{i,\leq p} &=& \{(i,j) \mid j \in \JJ(D) \hbox{ and } j \leq p\}
\hbox{ for } p \in \JJ(D);\\
\JJ_{i+1,< p} &=& \{(i+1,j) \mid j \in \JJ(D) \hbox{ and } j < p\}
\hbox{ for } p \in \JJ(D).
\end{eqnarray*}
Define the \bem{offspring} $D_p(i) = D \cup \JJ_{i+1,< p} \minus
\JJ_{i,\leq p}$ and $\mitosis_i(D) = \{D_p(i) \mid p \in \JJ(D)\}$.
Given a set ${\mathcal P}$ of pipe dreams, write $\mitosis_i({\mathcal
P}) = \bigcup_{D \in {\mathcal P}} \mitosis_i({\mathcal D})$.%
\end{defn}
}\end{excise}%

Proposition~\ref{prop:offspring} will present a verbal description of
mitosis that is different from the one in Section~\ref{sec:alg}.  The
new description is a little more algorithmic, using a certain local
transformation on pipe dreams that was discovered by Bergeron and
Billey.

\begin{defn}[\cite{BB}] \label{defn:chute}
A \bem{chutable rectangle} is a connected $2 \times k$ rectangle $C$
inside a pipe dream $D$ such that $k \geq 2$ and all but the following 3
locations in $C$ are crosses: the northwest, southwest, and southeast
corners.  Applying a \bem{chute move}%
        \footnote{The transpose of a chute move is called a \bem{ladder
        move} in \cite{BB}.}
to $D$ is accomplished by placing a \cross in the southwest corner of a
chutable rectangle $C$ and removing the \cross from the northeast corner
of the same~$C$.%
\end{defn}

Heuristically, a chute move therefore looks like:
\begin{rcgraph}
\begin{array}{@{}r|c|c|c|@{}c@{}|c|c|c|l@{}}
  \multicolumn{7}{c}{}&\multicolumn{1}{c}{
                                          \phantom{\!+\!}}&
  \multicolumn{1}{c}{\begin{array}{@{}c@{}}\\\adots\end{array}}
  \\\cline{2-8}
           &\cdot&\!+\!&\!+\!& \toplinedots &\!+\!&\!+\!&\!+\!
  \\\cline{2-4}\cline{6-8}
           &\cdot&\!+\!&\!+\!&              &\!+\!&\!+\!&\cdot
  \\\cline{2-8}
  \multicolumn{1}{c}{\begin{array}{@{}c@{}}\adots\\ \\ \end{array}}&
  \multicolumn{1}{c}{\phantom{\!+\!}}
\end{array}
\quad\stackrel{\rm chute}\rightsquigarrow\quad
%
\begin{array}{@{}r|c|c|c|@{}c@{}|c|c|c|l@{}}
  \multicolumn{7}{c}{}&\multicolumn{1}{c}{
                                          \phantom{\!+\!}}&
  \multicolumn{1}{c}{\begin{array}{@{}c@{}}\\\adots\end{array}}
  \\\cline{2-8}
           &\cdot&\!+\!&\!+\!& \toplinedots &\!+\!&\!+\!&\cdot
  \\\cline{2-4}\cline{6-8}
           &\!+\!&\!+\!&\!+\!&              &\!+\!&\!+\!&\cdot
  \\\cline{2-8}
  \multicolumn{1}{c}{\begin{array}{@{}c@{}}\adots\\ \\ \end{array}}&
  \multicolumn{1}{c}{\phantom{\!+\!}}
\end{array}
\end{rcgraph}

\begin{prop} \label{prop:offspring}
Let $D$ be a pipe dream, and suppose $j$ is the smallest column index
such that $(i+1,j) \not\in D$ and $(i,p) \in D$ for all $p \leq j$.
Then $D_p(i) \in \mitosis_i(D)$ is obtained from $D$ by
\begin{numbered}
\item
removing $(i,j)$, and then
\item
performing chute moves from row~$i$ to row~$i+1$, each one as far west
as possible, so that $(i,p)$ is the last {\rm \cross$\!$} removed.
\end{numbered}
\end{prop}
\begin{proof}
Immediate from Definitions~\ref{defn:chute} and~\ref{defn:mitosis}.
\end{proof}

\begin{example} \label{ex:offspring}
The offspring in Example~\ref{ex:mitosis} are listed in the order they
are born via the algorithmic `chute' form of mitosis in
Proposition~\ref{prop:offspring}, with $i = 3$.%
\begin{excise}{
  The dot in each pipe dream represents the last \cross removed.%
  \begin{figure}[t]
  \begin{rcgraph}
  \begin{array}{ccc}
  \begin{array}{@{}cc@{}}{}\\\\3&\\4\\\\\\\\\\\end{array}
  \begin{array}{|c|c|c|c|c|c|c|c|}
  \hline       &\!+\!&     &\!+\!&\!+\!&     &\ \, &\ \,
  \\\hline     &\!+\!&     &     &     &\!+\!&     &
  \\\hline\!+\!&\!+\!&\!+\!&\!+\!&     &     &     &
  \\\hline     &     &\!+\!&     &     &     &     &
  \\\hline\!+\!&     &     &     &     &     &     &
  \\\hline\!+\!&\!+\!&     &     &     &     &     &
  \\\hline\!+\!&     &     &     &     &     &     &
  \\\hline     &     &     &     &     &     &     &
  \\\hline
  \end{array}
  &\longmapsto&
  \begin{array}{|c|c|c|c|c|c|c|c|}
  \hline       &\!+\!&     &\!+\!&\!+\!&     &\ \, &\ \,
  \\\hline     &\!+\!&     &     &     &\!+\!&     &
  \\\hline\cdot&\!+\!&\!+\!&\!+\!&     &     &     &
  \\\hline     &     &\!+\!&     &     &     &     &
  \\\hline\!+\!&     &     &     &     &     &     &
  \\\hline\!+\!&\!+\!&     &     &     &     &     &
  \\\hline\!+\!&     &     &     &     &     &     &
  \\\hline     &     &     &     &     &     &     &
  \\\hline
  \end{array}
  \\\\[-8pt]
  \quad\ D
  &       &
  \downarrow
  \\[4pt]
  &       &
  \begin{array}{|c|c|c|c|c|c|c|c|}
  \hline       &\!+\!&     &\!+\!&\!+\!&     &\ \, &\ \,
  \\\hline     &\!+\!&     &     &     &\!+\!&     &
  \\\hline     &\cdot&\!+\!&\!+\!&     &     &     &
  \\\hline\!+\!&     &\!+\!&     &     &     &     &
  \\\hline\!+\!&     &     &     &     &     &     &
  \\\hline\!+\!&\!+\!&     &     &     &     &     &
  \\\hline\!+\!&     &     &     &     &     &     &
  \\\hline     &     &     &     &     &     &     &
  \\\hline
  \end{array}
  \\\\[-8pt]
  &       &
  \downarrow
  \\[4pt]
  &       &
  \begin{array}{|c|c|c|c|c|c|c|c|}
  \hline       &\!+\!&     &\!+\!&\!+\!&     &\ \, &\ \,
  \\\hline     &\!+\!&     &     &     &\!+\!&     &
  \\\hline     &     &\!+\!&\cdot&     &     &     &
  \\\hline\!+\!&\!+\!&\!+\!&     &     &     &     &
  \\\hline\!+\!&     &     &     &     &     &     &
  \\\hline\!+\!&\!+\!&     &     &     &     &     &
  \\\hline\!+\!&     &     &     &     &     &     &
  \\\hline     &     &     &     &     &     &     &
  \\\hline
  \end{array}
  \end{array}
  \end{rcgraph}
  \caption{Mitosis} \label{fig:mitosis}
  \end{figure}%
}\end{excise}
\end{example}

Proposition~\ref{prop:offspring} for mitosis has the following analogue
for mutation.

\begin{claim} \label{claim:mu}
Suppose that $\length(ws_i) < \length(w)$, and let $\zz^\bb \not\in J_w$
be a squarefree monomial of maximal support.  If $D = D(\bb)$ then
\begin{eqnarray*}
  |\prom(\bb)| &=& |\JJ(D)|.
\end{eqnarray*}
If $0 \leq d < |\prom(\bb)|$, then $D(\mu_i^{2d+1}(2\bb))$ is obtained
from $D$ by
\begin{numbered}
\item
removing $(i,j)$, where $j$ is as in Proposition~\ref{prop:offspring},
and then
\item
performing $d$ chute moves from row~$i$ to row~$i+1$, each as far west
as possible.
\end{numbered}
\end{claim}
\begin{proof}
By Lemma~\ref{lemma:max}, the columns in $\JJ(D)$ are in bijection with
the nonzero entries in the promoter of $\bb$, each of which is a~$1$ in
row~$i+1$.  The final statement follows easily from the definitions.
\end{proof}

\begin{example} \label{ex:mitosis'}
The left array $\bb$ in Fig.~\ref{fig:mu} has maximal support among
exponent arrays on monomials not in $J_w$, where $w = 13865742$ as in
Examples~\ref{ex:intro}, \ref{ex:mutate}, and~\ref{ex:prom}.
Substituting \cross for each blank space and then removing the numbers
and dots yields the left pipe dream Example~\ref{ex:mitosis}.  Applying
the same makeover to the middle column of Fig.~\ref{fig:mu} results in
the offspring $D(\mu_3^1(\bb))$, $D(\mu_3^3(\bb))$, and
$D(\mu_3^5(\bb))$.%
\end{example}

\begin{lemma} \label{lemma:agree}
If\/ $\length(ws_i) < \length(w)$ and $\zz^\bb \not\in J_w$ is squarefree
of maximal support, then $\mitosis_i(D(\bb)) = \{D(\mu_i^{2d+1}(2\bb))
\mid 0 \leq d < |\prom(\bb)|\}$.
\end{lemma}
\begin{proof}
Compare Claim~\ref{claim:mu} with Proposition~\ref{prop:offspring}.
\end{proof}

Theorem~\ref{thm:alg} will conclude the translation of mutation on
monomials into mitosis on facets, but the translation requires an
intermediate result.

\begin{lemma} \label{lemma:facets}
If\/ $\length(ws_i) < \length(w)$, then $\{D_L \mid L$ is a facet of
$\LL_{ws_i}\}$ is the set of pipe dreams $D(\mu_i^{2d+1}(2\bb))$ such
that $\zz^\bb \not\in J_w$ is squarefree of maximal support and $0 \leq d
< |\prom(\bb)|$.
\end{lemma}
\begin{proof}
By Theorem~\ref{thm:ev} and Lemma~\ref{lemma:notafacet}, every facet of
$\LL_{ws_i}$ is the support of a mutation $\mu_i^d(\zz^{2\bb})$ for some
monomial $\zz^{2\bb} \not\in J_w$ and $d \leq |\prom(2\bb)|$.
Furthermore, it is clear from Lemma~\ref{lemma:notafacet} and the
definition of mutation that we may assume $\zz^\bb$ is a squarefree
monomial (so the entries of $2\bb$ are all~$0$ or~$2$).  In this case,
the supports of odd mutations $\mu_i^{2d+1}(\zz^{2\bb})$ for $0 \leq d <
|\prom(\bb)|$ are facets of $\LL_{ws_i}$ by Proposition~\ref{prop:pure}
because they each have cardinality $n^2 - \length(w) + 1 = n^2 -
\length(ws_i)$, while the supports of even mutations
$\mu_i^{2d}(\zz^{2\bb})$ for $d \leq |\prom(\bb)|$ are not facets, each
having cardinality $n^2 - \length(w)$.
\end{proof}

\begin{thm} \label{thm:alg}
If\/ $\length(ws_i) < \length(w)$, then
\begin{eqnarray*}
  \{D_L \mid L \hbox{ is a facet of } \LL_{ws_i}\} &=& \mitosis_i(\{D_L
  \mid L \hbox{ is a facet of } \LL_w\}).
\end{eqnarray*}
Moreover, $\mitosis_i(D_L) \cap \mitosis_i(D_{L'}) = \nothing$ if $L \neq
L'$ are facets of $\LL_w$.
\end{thm}
\begin{proof}
The displayed equation is a consequence of Lemma~\ref{lemma:facets} and
Lemma~\ref{lemma:agree}, so we concentrate on the final statement.  Let
$\zz^\bb \not\in J_w$ be a squarefree monomial of maximal support $L =
\supp(\bb)$, and let $0 \leq d < |\prom(\bb)|$.  The entries of the array
$\mu_i^{2d+1}(2\bb)$ are all either~$0$ or~$2$, except for precisely two
$1$'s, both in the same column~$p$ (the boldface entries in
Fig.~\ref{fig:mu}, middle column).  By Lemma~\ref{lemma:max},~$p$ is the
westernmost column of $D(\mu_i^{2d+1}(2\bb))$ in which neither row~$i$
nor row~$i+1$ has a cross.

Now suppose $\zz^{\bb'} \not\in J_w$ is another squarefree monomial of
maximal support $L'$, and let $0 \leq d' < |\prom(\bb')|$.  If
$D(\mu_i^{2d+1}(2\bb)) = D(\mu_i^{2d'+1}(2\bb'))$, then the argument in
the first paragraph of the proof implies that $\mu_i^{2d+1}(2\bb) =
\mu_i^{2d'+1}(2\bb')$, since they have the same entries equal to $1$ as
well as the same support, and all of their other nonzero entries
equal~$2$.  We conclude that $\bb = \bb'$ by Lemma~\ref{lemma:unique}.
Using Lemma~\ref{lemma:facets}, we have proved that $\mitosis_i(D_L)
\cap \mitosis_i(D_{L'}) \neq \nothing$ implies $L = L'$.
\end{proof}

\begin{example} \label{ex:remove}
The pipe dream $D$ in Example~\ref{ex:mitosis'} and the left side of
Example~\ref{ex:mitosis} is $D_L$ for a facet of $\LL_{13865742}$.  By
Theorem~\ref{thm:alg}, the three pipe dreams at the right of
Example~\ref{ex:mitosis} can be expressed as $D_{L'}$ for facets $L' \in
\LL_{13685742}$, where $13685742 = 13865742 \cdot s_3$.%
\end{example}

\section{Facets and reduced pipe dreams}\label{sec:rp}

To~derive the connection between the initial complex~$\LL_w$ and reduced
pipe dreams, in Theorem~\ref{thm:rp}, we need the next result, whose
proof connects chuting with antidiagonals.

\begin{lemma} \label{lemma:chuting}
The set $\{D_L \mid L \in \facets(\LL_w)\}$ is closed under chute moves.
\end{lemma}
\begin{proof}
A pipe dream $D$ is equal to $D_L$ for some (not necessarily maximal)
$L \in \LL_w$ if and only if $D$ meets every antidiagonal in $J_w$,
which by definition of $\LL_w$ equals $\bigcap_{L \in \LL_w} \<z_{qp}
\mid (q,p) \in D_L\>$.  Suppose that $C$ is a chutable rectangle in
$D_L$ for $L \in \LL_w$.  It is enough to show that the intersection
$a \cap D_L$ of any antidiagonal $a \in J_w$ with $D_L$ does not
consist entirely of the single cross in the northeast corner of~$C$,
unless $a$ also contains the southwest corner of~$C$.  Indeed, the
purity of $\LL_w$ (Proposition~\ref{prop:pure}) will then imply that
chuting $D_L$ in $C$ yields $D_{L'}$ for some {\em facet}\/ $L'$
whenever $L \in \facets(\LL_w)$.

To prove the claim concerning $a \cap D_L$, we may assume $a$ contains
the cross in the northeast corner $(q,p)$ of $C$, but not the cross in
the southwest corner of~$C$, and split into cases:
\begin{romanlist}
\item
$a$ does not continue south of row~$q$.

\item
$a$ continues south of row~$q$ but skips row~$q+1$.

\item
$a$ intersects row~$q+1$, but strictly east of the southwest corner of
$C$.

\item
$a$ intersects row~$q+1$, but strictly west of the southwest corner of
$C$.
\end{romanlist}
Letting $(q+1,t)$ be the southwest corner of $C$, construct new
antidiagonals $a'$ that are in $J_w$ (and hence intersect $D_L$) by
replacing the cross at $(q,p)$ with a cross at:
\begin{romanlist}
\item
$(q,t)$, using Lemma~\ref{lemma:squish}(W);

\item
$(q+1,p)$, using Lemma~\ref{lemma:squish}(S);

\item
$(q,p)$, so $a = a'$ trivially; or

\item
$(q,t)$, using Lemma~\ref{lemma:squish}(W).
\end{romanlist}
Observe that in case~(iii), $a$ already shares a box in row~$q+1$ where
$D_L$ has a cross.  Each of the other antidiagonals $a'$ intersects both
$a$ and $D_L$ in some box that is not $(q,p)$, since the location of $a'
\minus a$ has been constructed not to be a cross in $D_L$.%
\end{proof}

Lemma~\ref{lemma:subword} implies the following criterion for when
removing a~\cross from a pipe dream $D \in \rc(w)$ yields a pipe dream
in~$\rc(ws_i)$.  Specifically, it concerns the removal of a cross at
$(i,j)$ from configurations that look like
$$
\begin{array}{lccccccccc}
&\perm 1{}&\perm{}{}&\perm{}{}&\perm{}{}&
   \perm{\hskip-11.6pt \textstyle \cdots}{}&\perm{}{}&\perm{}{}&\perm j\
\\[3pt]
\petit{i}  &\+ &  \+ &  \+ &  \+ &  \+ &  \+ &  \+ &  \+ &\toplinedots
\\[-1pt]
\petit{i+1}&\+ &  \+ &  \+ &  \+ &  \+ &  \+ &  \+ &  \jr
\end{array}
\begin{array}{c}\\ = \\\end{array}\
\begin{array}{l|c|c|c|c|c|c|c|c|c}
\multicolumn{1}{c}{}
        &\multicolumn{1}{c}{\petit 1}
                &\multicolumn{6}{c}{\cdots}
                        &\multicolumn{1}{c}{\petit j}
                                \\[3pt]\cline{2-10}
\petit{i}  &\!+\!&\!+\!&\!+\!&\!+\!&\!+\!&\!+\!&\!+\!&\!+\!&\toplinedots
                                                            \\\cline{2-9}
\petit{i+1}&\!+\!&\!+\!&\!+\!&\!+\!&\!+\!&\!+\!&\!+\!&\cdot&\\\cline{2-10}
\end{array}
$$
at the west end of rows~$i$ and~$i+1$ in~$D$.

\begin{lemma} \label{lemma:rp}
Let $D \in \rc(w)$ and $j$ be a fixed column index with $(i+1,j) \not\in
D$, but $(i,p) \in D$ for all $p \leq j$, and $(i+1,p) \in D$ for all $p
< j$.  Then $\length(ws_i) < \length(w)$, and if $D' = D \minus (i,j)$
then $D' \in \rc(ws_i)$.
\end{lemma}
\begin{proof}
Removing $(i,j)$ only switches the exit points of the two pipes starting
in rows $i$ and $i+1$, so the pipe starting in row $k$ of $D'$ exits out
of column $ws_i(k)$ for each $k$.  The result follows from
Lemma~\ref{lemma:subword}.
\end{proof}

The connection beween the complexes $\LL_w$ and reduced pipe dreams
requires certain facts proved by Bergeron and Billey \cite{BB}.  The
next lemma consists mostly of the combinatorial parts (a),~(b), and~(c)
of \cite[Theorem~3.7]{BB}, their main result.  Its proof there relies
exclusively on elementary properties of reduced pipe dreams.

\begin{lemma}[\cite{BB}] \label{lemma:BB} \mbox{}
\begin{numbered}
\item \label{closed}
The set $\rc(w)$ of reduced pipe dreams for~$w$ is closed under chute
operations.

\item \label{unique}
There is a unique \bem{top reduced pipe dream} for~$w$ such that every
cross not in the first row has a cross due north of it.

\item
Every reduced pipe dream for~$w$ can be obtained by applying a sequence
of chute moves to the top reduced pipe dream for~$w$.
\end{numbered}
\end{lemma}

\begin{lemma} \label{lemma:top}
The top reduced pipe dream for~$w$ is $D_L$ for a facet $L \in \LL_w$.
\end{lemma}
\begin{proof}
The unique reduced pipe dream for $w_0$, whose crosses lie at $\{(q,p)
\mid q+p \leq n\}$, is also $D_L$ for the unique facet of $\LL_{w_0}$.
By Lemma~\ref{lemma:n!}, take $w$ and~$ws_i$ as in
Lemma~\ref{lemma:alpha} and assume the result for~$w$.  By
Lemma~\ref{lemma:n!} again, the top pipe dream $D \in \rc(w)$ satisfies
the conditions of Lemma~\ref{lemma:rp}, where $\alpha$ lies at
$(i+1,j)$.  Now combine Lemma~\ref{lemma:rp} with the last sentence of
Lemma~\ref{lemma:alpha} to prove the desired result for~$ws_i$.%
\end{proof}

\begin{thm} \label{thm:rp}
$\rc(w) = \{D_L \mid L$ is a facet of $\LL_w\}$, where $D_L = [n]^2
\minus L$.  In other words, reduced pipe dreams for~$w$ are complements
of maximal supports of monomials $\not\in J_w$.
\end{thm}
\begin{proof}
Lemma~\ref{lemma:top} and Lemma~\ref{lemma:chuting} imply that $\rc(w)
\subseteq \{D_L \mid L \in \facets(\LL_w)\}$, given
Lemma~\ref{lemma:BB}.  Since the opposite containment $\{D_L \mid L \in
\facets(\LL_v)\} \subseteq \rc(v)$ is obvious for $v = w_0$, it suffices
to prove it for $v = ws_i$ by assuming it for $v = w$.

A pair $(i+1,j)$ as in Lemma~\ref{lemma:rp} exists in $D$ if and only if
the set $\JJ(D)$ in Definition~\ref{defn:mitosis} is nonempty, which
occurs if and only if $\mitosis_i(D) \neq \nothing$.  In this case, the
first offspring of~$D$ under $\mitosis_i$ (as in
Proposition~\ref{prop:offspring}) is $D' = D \minus (i,j) \in
\rc(ws_i)$.  The desired containment follows from Theorem~\ref{thm:alg}
and statement~\ref{closed} of Lemma~\ref{lemma:BB}.
\end{proof}

\section{Proof of Theorems~\ref{t:formulae}, \ref{t:gb}, 
and~\ref{t:mitosis}}\label{sec:proof}

In this section, we tie up loose ends, completing the proofs of
Theorems~\ref{t:formulae}, \ref{t:gb}, and~\ref{t:mitosis}.  All
statements in these theorems are straightforward for the long
permutation $w = w_0$, so we may prove the rest by Bruhat
induction---that is, by downward induction on~$\length(w)$.

The multidegrees of $\{\ol X_w\}_{w \in S_n}$ of matrix Schubert
varieties satisfy the divided difference recursion defining double
Schubert polynomials by Theorem~\ref{thm:oracle}.  Separately, the
Hilbert series of $\{\cj{w}\}_{w \in S_n}$ satisfy the Demazure
recursion defining the double Grothendieck polynomials by
Theorem~\ref{thm:induction}.  It follows by Lemma~\ref{lemma:schubert}
that the multidegrees $\{[\LL_w]\}_{w \in S_n}$ satisfy the divided
difference recursion.

The subspace arrangement $\LL_w$ is equidimensional by
Proposition~\ref{prop:pure}, and the multidegree of~$\ol X_w$ equals
that of $\kk[\zz]/\IN(I(\ol X_w))$ under any term order, by degeneration
in Theorem~\ref{t:multidegs}.  Moreover, if the term order is
antidiagonal, then $J_w \subseteq \IN(I_w) \subseteq \IN(I(\ol X_w))$,
whence $J_w = \IN(I_w) = \IN(I(\ol X_w))$ by Lemma~\ref{lemma:IN} with
$I' = \IN(I(\ol X_w))$ and $J = J_w$.  The primality of~$I_w$ in
Theorem~\ref{t:formulae} follows because $\IN(I_w) = \IN(I(\ol X_w))$
and $I_w \subseteq I(\ol X_w)$.  The Gr\"obner basis statement in
Theorem~\ref{t:gb} is an immediate consequence.  Hilbert series are
preserved under taking initial ideals, so the \K-polynomials and
multidegrees of $\{\ci w\}_{w \in S_n}$ satisfy the Demazure and divided
difference recursions, proving Theorem~\ref{t:formulae}.

The last sentence of Theorem~\ref{t:gb} is Theorem~\ref{thm:rp}.  Using
that, we conclude the remaining statement in Theorem~\ref{t:gb}, namely
shellability and Cohen--Macaulayness, by applying
Theorem~\ref{t:subword} to Example~\ref{ex:complex}.  Mitosis in
Theorem~\ref{t:mitosis} comes from applying Theorem~\ref{thm:rp} to the
mitosis in Theorem~\ref{thm:alg}.%

\raggedbottom
\def\cprime{$'$}
\providecommand{\bysame}{\leavevmode\hbox to3em{\hrulefill}\thinspace}


\end{document}